\documentclass[12pt, reqno]{amsart}

\usepackage[margin=1in]{geometry} 

\usepackage{times}
\usepackage[english]{babel}
\usepackage[utf8]{inputenc}
\usepackage{csquotes}
\usepackage{cancel}
\usepackage{hyperref}
\usepackage{amsmath,amssymb,amsthm,bbm}
\usepackage{amsfonts}
\usepackage{mathtools}
\usepackage{tabulary}
\usepackage{bm}
\usepackage{graphicx, color}
\usepackage{svg}
\usepackage{enumitem}
\usepackage{tikz-cd}
\usepackage{physics}
\usepackage{tikz}
\usepackage{todonotes}
\usepackage[outline]{contour} 
\usetikzlibrary{patterns,snakes}
\usetikzlibrary{arrows.meta} 
\contourlength{0.4pt}

\allowdisplaybreaks

\makeatletter
\@namedef{subjclassname@2020}{\textup{2020} Mathematics Subject Classification}
\makeatother

\tikzstyle{mydashed}=[black,dashed,line width=0.25,dash pattern=on 2.2pt off 2.2pt]
\tikzstyle{axis}=[->,thick] 
\tikzstyle{ell}=[{Latex[length=3.3,width=2.2]}-{Latex[length=3.3,width=2.2]},line width=0.3]
\tikzstyle{dx}=[-{Latex[length=3.3,width=2.2]},black,line width=0.3]
\tikzstyle{ground}=[preaction={fill,top color=black!10,bottom color=black!5,shading angle=20},
                    fill,pattern=north east lines,draw=none,minimum width=0.3,minimum height=0.6]
\tikzstyle{mass}=[line width=0.8,black,fill=black,rounded corners=1,
                  top color=black,bottom color=black,shading angle=20]
\tikzstyle{spring}=[line width=0.8,black,snake=coil,segment amplitude=5,segment length=5,line cap=round]
\tikzset{>=latex} 
\tikzstyle{force}=[->,myred,very thick,line cap=round]

\usetikzlibrary{decorations.pathmorphing}
\usepackage{parskip}
\numberwithin{equation}{section}

\newcommand{\numberset}{\mathbb}
\newcommand{\N}{\numberset{N}}
\newcommand{\R}{\numberset{R}}
\newcommand{\C}{\numberset{C}}

\newcommand{\E}{\mathcal{E}}
\newcommand{\D}{\mathbb{D}}

\DeclareMathOperator*{\argmin}{arg\,min}
\newtheorem{theorem}{Theorem}[section]
\newtheorem{definition}{Definition}[section]
\newtheorem{proposition}{Proposition}[section]
\newtheorem{corollary}{Corollary}[section]
\newtheorem{lemma}{Lemma}[section]
\theoremstyle{definition}

\newcommand{\EEE}{\color{black}}

\title[Finite-strain Poynting-Thomson model]{Finite-strain Poynting-Thomson model: \\ existence and
  linearization}
 
\author[A. Chiesa]{Andrea Chiesa}
	\address[Andrea Chiesa]{University of
		Vienna, Faculty of Mathematics and  Vienna School of Mathematics,
                Oskar-Morgenstern-Platz 1, A-1090 Vienna, Austria.}
              \email{andrea.chiesa@univie.ac.at}
              \urladdr{https://www.mat.univie.ac.at/$\sim$achiesa/}
               
\author[M. Kru\v{z}\'{\i}k]{Martin Kru\v z\`\i k}
\address[Martin Kru\v{z}\'ik]{Czech Academy of Sciences, Institute of Information Theory and Automation, Pod vod\' arenskou ve\v z\' \i\ 4, 182 08, Prague 8, Czech Republic and Faculty of Civil Engineering, Czech Technical University, Th\'{a}kurova 7, 166 29, Prague 6, Czech Republic.}
\email{kruzik@utia.cas.cz}
	\urladdr{http://staff.utia.cas.cz/kruzik/}

\author[U. Stefanelli]{Ulisse Stefanelli} 
	\address[Ulisse Stefanelli]{University of
		Vienna, Faculty of Mathematics,
                Oskar-Morgenstern-Platz 1, A-1090 Vienna, Austria, 
		University of Vienna, Vienna Research Platform on Accelerating
		Photoreaction Discovery, W\"ahringerstra\ss e 17, 1090
                Wien, Austria, and Istituto di	Matematica Applicata e Tecnologie Informatiche {\it E. Magenes}, via
		Ferrata 1, I-27100 Pavia, Italy.}
	\email{ulisse.stefanelli@univie.ac.at}
	\urladdr{http://www.mat.univie.ac.at/$\sim$stefanelli}

\subjclass[2020]{49J20, 
74D10}
\keywords{Poynting-Thomson model, variational approach, existence, linearization.}

\begin{document}

	\begin{abstract}
     We analyse  the finite-strain   Poynting-Thomson viscoelastic
    model.  In its linearized small-deformation limit, this
    corresponds to the  serial connection of an elastic
    spring and a Kelvin-Voigt viscoelastic element.  In the
    finite-strain case, the 
    total
    deformation of the body results from the composition of two maps,
    describing  the deformation of the viscoelastic
    element  and   the  elastic
    one,  respectively.    We prove the existence of  suitably  weak solutions by a
    time-discretization approach based on incremental
    minimization. Moreover,  we prove a rigorous linearization
    result, showing that the corresponding small-strain model is
    indeed recovered
    in the small-loading limit.  

	\end{abstract}

        \maketitle
        

	\section{Introduction}
	\label{sec:intro}

 Viscoelastic solids appear ubiquitously in
applications. Polymers, rubber, biomaterials, wood, clay, and soft solids,
 including  metals at close-to-melting temperatures, behave
viscoelastically.
 The   mechanical response  of viscoelastic solids  is governed by the interplay between elastic
and viscous dynamics: by applying stresses both strains and
strain rates ensue \cite{Nahn}. This is at the basis of different
effects, from viscoelastic creep, to viscous relaxation, to
rate-dependence in material response, to dissipation of mechanical
energy \cite{Shaw}.

The modelization of viscoelastic solid response dates back to the early days
of Continuum Mechanics. In the linearized, infinitesimal-strain setting of the standard-solid rheology, two
basic models  are the {\it Maxwell} and the {\it Kelvin-Voigt} one, where an
elastic spring is connected to a viscous dashpot in series or in
parallel, respectively.  These models offer only a simplified description
of actual viscoelastic  behavior.  
More accurate descriptions  necessarily call  for  more complex models. 
 A
first option in this direction is the {\it Poynting-Thomson} model,
resulting from the combination in series of   an elastic and a
Kelvin-Voigt component, see Figure \ref{fig:Poynting Thomson}.  
A second option would be the {\it Zener} model, which
consists of an elastic and a Maxwell element in parallel. Note
however, that Poynting-Thomson and Zener can be proved to be
equivalent in the linearized setting, see \cite[Remark 6.5.4]{Kruzik Roubicek}.

The aim of this paper is to investigate the Poynting-Thomson
model in the finite-strain setting. From the modeling
viewpoint, extending the model beyond the small-strain case is crucial, for
viscoelastic materials commonly experience large deformations. In
fact, finite-strain versions of the Poynting-Thomson
model have already been considered. The reader is referred to
\cite{PoyntingZerner}, where a comparison between Poynting-Thomson and
Zener models at finite strains is discussed, and to 
\cite{Meo}, focusing on the  anisothermal  version the 
Poynting-Thomson  model.

To the best of our knowledge, mathematical results on the
finite-strain Poynting-Thomson model are still not available. The
focus of this paper is to fill this gap by presenting
\begin{itemize}
\item an {\bf existence theory} for solutions  of the
  finite-strain Poynting-Thomson model,  as well as a convergence
  result for time-discretizations  (Theorem \ref{thm: existence
  limit}); \vspace{2mm}
  \item a rigorous
  {\bf linearization result}, proving that finite-strain solutions converge
  (up to subsequences) to solutions of the linearized system in the
  limit of small loadings and, correspondingly, small strains  (Theorem \ref{thm limit lin}). 
\end{itemize}\vspace{3mm}

  Our analysis is variational in nature.  The convergence result 
 provides a rigorous
   counterpart  to the classical heuristic arguments based on Taylor
  expansions \cite{PoyntingZerner}.

 We postpone to Section \ref{sec:model} both the detailed
discussion of the model and a first presentation of our main
results. We anticipate however here that the theory requires no
second-gradient terms but rather relies on a decomposition of the
total deformation in terms of an elastic and a viscous deformation,
see \eqref{composition} below. Correspondingly, the variational
formulation of the problem features both Lagrangian and Eulerian
terms. Note moreover that the viscous dissipation is  here  assumed to be $p_\psi$-homogeneous,
with superlinear homogeneity $p_\psi\geq 2$.

Our notion of solution (see Definition \ref{def:approximable sol}) hinges on the validity
of an energy inequality, an elastic semistability inequality, and an
approximability property via time-discrete problems. Albeit very weak,
this notion replicates the important features of viscoelastic
evolution, including elastic equilibrium, energy dissipation, and viscous relaxation.

Before moving on, let us put our results in context
with respect to the available literature. In the purely PDE setting,
existence results for viscoelastic  dissipative  systems are
classical. The reader is referred to the recent monograph \cite{Kruzik
Roubicek} for a comprehensive collection of references. As it is well
known, the PDE
setting is local in nature and, as such, does not allow considering
global constraints such as  injectivity of deformations, i.e., 
noninterpenetration   of matter. Variational theories for viscoelastic
evolution offer a remedy in this respect. By making use of the
underlying gradient-flow structure of viscoelastic evolution,
existence results for variational
solutions have been obtained in the one-dimensional \cite{MOS} and in the
multi-dimensional case \cite{Friedrichs}. The latter paper also delivers a
rigorous evolutive $\Gamma$-convergence linearization result.  See
also \cite{KR} for the case of self-contact and \cite{BFK,MR20} for some
extension to nonisothermal situations. With respect to these
contributions, we deal here with an internal-variable formulation,
where the elastic variable does not dissipate. From the technical
viewpoint,  
the novelty of our approach resides in avoiding the
second-gradient theory by virtue of the composition assumption
\eqref{composition}.  This impacts on the functional setting, as well
as on the required mathematical  techniques.

In the different but related frame of activated inelastic
deformations,  the 
closest contributions to ours are \cite{MielkeRossiSavare} and
\cite{Roger}, both dealing with rate-dependent viscoplasticity
($p_\psi>1$) under the multiplicative-decomposition setting. In both papers, existence of solutions is discussed, by
taking into account additional gradient-type terms for the viscous
strain. In particular, the full gradient is considered in \cite{MielkeRossiSavare}, whereas in
\cite{Roger} only its curl is penalized. The approach in \cite{Roger}
analogous to ours in terms of solution notion, despite the differences in the model.
In contrast with
these papers, viscous evolution is here not activated. In addition,
 by
not considering here additional gradient terms, we avoid 
introducing   a
second length scale in the model and thus tackle so-called {\it simple materials}.
Moreover, we investigate  here linearization, which was not discussed
in \cite{MielkeRossiSavare,Roger}. 

In the  fully rate-independent setting $p_\psi=1$ of activated
  elastoplasticity, the
papers \cite{Melching,Stefanelli1} and \cite{MielkeStefanelli}, contribute
an existence and linearization theory  which is  parallel
the current viscoelastic one. More precisely,
\cite{Melching,Stefanelli1} deal with a decomposition of deformations
in the same spirit of \eqref{composition} below,  
avoiding the use of second gradients, whereas \cite{MielkeStefanelli}
features no gradients, but is a pure convergence result, in a setting where
existence is not known. With respect to these contributions, the
superlinear,  non activated  nature of the
dissipation of the  viscous setting calls for using a different set of analytical
tools from gradient-flow theory \cite{MRS08b}.
Note that, also in the  rate-independent  setting, by including a gradient
term of the plastic strain, hence resorting to so-called {\it
  strain-gradient} finite plasticity, one obtains stronger results. In particular, the existence of
energetic solutions in strain-gradient finite plasticity is in \cite{MM09} and the linearization in some
symmetrized case is in \cite{GS2}. Under the
mere penalization of the curl of the gradient of the plastic strain,
existence of incremental solutions is proved in \cite{Mielke-Mueller} and
linearization is in \cite{curl}.

The paper is organized as follows.  In Section \ref{sec:model},  we  provide an illustration of the
finite-strain Poynting-Thomson model under consideration, as well as 
an introduction  to  our main results.
  Some preliminary material
and comment on the functional setting  is provided in Section
\ref{sec:pre}.   In particular, we discuss the
set of admissible deformations in Subsection \ref{sec:admin}. In
Subsections \ref{sec:hp existence} and \ref{sec:hp linearization} we list and comment the assumptions,
whereas the
statements of our main results, Theorem \ref{thm: existence limit} and
Theorem \ref{thm limit lin} are presented in Subsections \ref{sec:res existence} and \ref{sec:res linearization}, respectively. The solvability of the time-discrete incremental
problems is discussed in Section \ref{sec:time discr},
whereas the proofs of  Theorems  \ref{thm: existence limit} and
\ref{thm limit lin} are given in Sections \ref{sec:Proof of Existence} and \ref{sec:Proof of Linearization},
respectively.

\section{The finite-strain Poynting-Thomson model}\label{sec:model}

  %

 \begin{figure}[h!]
  \centering
  \includegraphics[trim={4cm 20cm 10cm 4cm},clip,width=0.55\textwidth]{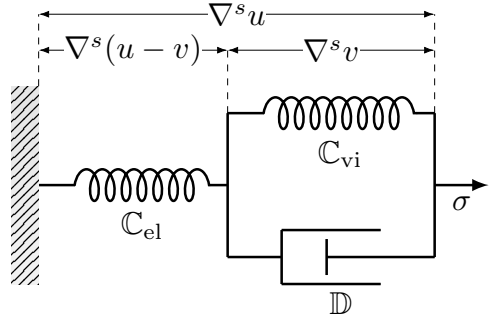}

  \caption{ The Poynting-Thomson rheological model (linearized
    setting). }
    \label{fig:Poynting Thomson}
  \end{figure}

In order to illustrate our results, we start by recalling the
classical Poynting-Thomson in the linearized setting of infinitesimal
strains. By indicating
by $u:\Omega \to \R^d$ the {\it infinitesimal displacement} from the
{\it reference configuration} $\Omega \subset \R^d$, the total
strain $\nabla^s u$ (here $\nabla^s$ denotes the symmetrized gradient
 $\nabla^su = (\nabla u + \nabla u^\top)/2$)  is
additively decomposed in its elastic and its viscous parts as
$\nabla^s u = \C_{\rm el}^{-1} \sigma + \nabla^s v$. In its
quasistatic approximation, the evolution of the body results from the
combination of the equilibrium system and the constitutive relation, namely,
\begin{align*}
-\operatorname{div}\left(\C_{\rm el} \nabla^s(u-v)\right)&=f  \quad
  \text{ in } \Omega \times (0,T),   
  \\
\mathbb{D} \nabla^s \dot{v}+(\C_{\rm vi}{+} \C_{\rm el})\nabla^s v  
&=\C_{\rm el} \nabla^s u    \quad  \text{ in } \Omega \times (0,T), 
\end{align*}
where $f$ stands for a given body force and $\dot v$ denotes the time-derivative of $v$. The reader is referred to the
  monographs \cite{Fabrizio92,Kruzik Roubicek,Visintin} for a
  comprehensive collection of analytical results. Let us remark
  that in this paper we will specifically consider the case of incompressible viscosity, i.e., in the linearized setting $\operatorname{tr}v=0$. Hence, the evolution of the system considered is actually determined by the following equations
  \begin{align}
    -\operatorname{div}\left(\C_{\rm el} \nabla^s(u-v)\right)&=f  \quad
      \text{ in } \Omega \times (0,T),   \label{linear PDE1}  \\
    \mathbb{D} \nabla^s \dot{v}+(\C_{\rm vi}{+} \C_{\rm el})\nabla^s v  
    &=\operatorname{dev}(\C_{\rm el} \nabla^s u)    \quad  \text{ in } \Omega \times (0,T) \label{linear PDE2},
    \end{align} 
    where $\operatorname{dev}$ denotes the deviatoric part of a
    tensor. Restricting to the incompressible case  would call for
    accordingly specifying the rheological diagram from Figure
    \ref{fig:Poynting Thomson} by distinguishing the volumetric and
    the deviatoric components.  \EEE

   In the finite-strain Poynting-Thomson model
  \cite{PoyntingZerner,Meo}, the   state of the viscoelastic system is specified in terms of its
{\it deformation} $y:\Omega\rightarrow \R^d$.  As it is common in finite-strain theories \cite{Lee69}, the
deformation gradient $\nabla y$ is {\it multiplicatively decomposed}
as $\nabla y=F_{\rm el}F_{\rm vi},$
where $F_{\rm el} $ and $F_{\rm vi}$ are the {\it elastic} and {\it
  viscous strain} tensors, representing the elastic and viscous
response of the medium, respectively. 

 A distinctive feature of
our approach is that we assume the viscous strain to
be {\it compatible}: we identify $F_{\rm vi}$  with the
gradient $\nabla y_{\rm vi}$ of a
 {\it viscous deformation} $y_{\rm vi}:\Omega\rightarrow y_{\rm
  vi}(\Omega)\subset\R^d$,  mapping  the {\it reference}
configuration $\Omega$ to  the  {\it intermediate} one 
$y_{\rm vi}(\Omega)$. Correspondingly, the elastic strain is
compatible as well and there   exists an elastic deformation
$y_{\rm el}: y_{\rm vi}(\Omega)\rightarrow\R^d$   with $F_{\rm
  el}=\nabla y_{\rm el}$ mapping the
intermediate configuration to the {\it actual} one. As such, the
multiplicative decomposition $\nabla y = F_{\rm el}F_{\rm vi}$ ensues
from an application of the classical chain rule to the composition 
\begin{equation}
   y:=y_{\text{\rm el}}\circ y_{\text{\rm vi}}:\Omega \rightarrow
   \R^d. \label{composition}
\end{equation}
 Moving from this position, the state of the medium is described
by the pair $( y_{\rm vi},y_{\rm el})$, effectively distinguishing viscous and
elastic responses.

Before moving on, let us stress that the compatibility assumption on
$F_{\rm vi}$, whence the composition assumption \eqref{composition},
realistically describes a variety of viscoelastic evolution
settings and refer to \cite{Melching,Stefanelli1} for some parallel theory in the frame of
  finite-strain plasticity.  In particular, position
  \eqref{composition}  is flexible enough to cover both limiting
  cases of a purely  elastic ($y_{\rm vi}={\rm id}$)  and of a  plain   Kelvin-Voigt 
  ($y_{\rm el}={\rm id}$)  materials. In the
  linearized setting, these would formally correspond to the cases ${\mathbb
    C}_{\rm vi}\to \infty$ and ${\mathbb
    C}_{\rm el}\to \infty$, respectively. Let us note that by
  choosing ${\mathbb
    C}_{\rm vi}=0$ the linearized system \eqref{linear
    PDE1}-\eqref{linear PDE2} reduces to the {\it Maxwell} fluidic
  rheological model. By assuming
  \eqref{composition} we exclude the onset of defects, such as
  dislocations and disclinations. Albeit this could limit the
  application of the theory in some specific cases, it is to 
  remark that viscous materials are often amorphous, so that the
  relevance of strictly crystallographic descriptions may be
  questionable. From the more analytical viewpoint, assumption \eqref{composition}  allows us to present a
comprehensive mathematical theory within the setting of so-called {\it
  simple materials}, i.e., without resorting to second-gradient
theories. The alternative path of including second-order deformation
gradients, is also viable and, as far as existence is concerned, has
been  considered  in \cite{MielkeRossiSavare}  in the
activated case of viscoplasticity.   

A first consequence of the composition  \eqref{composition} is that the
elastic deformation $y_{\rm el}$ is defined on
the a-priori unknown intermediate configuration $ y_{\rm vi}(\Omega)$,
making the analysis delicate. In particular,  
the variational description of the viscoelastic behavior  results
in a mixed Lagrangian-Eulerian variational problem.  This mixed
nature of the problem will be tamed by means of change-of-variables
techniques, which in turn ask for some specification on the class of
admissible intermediate configurations.  Let us anticipate that
 $y_{\rm vi}$ will
be  required  to
be an incompressible  ($\det \nabla y_{\rm vi}=1$)  homeomorphism
throughout. We refer to \cite{ Haupt-Lion, Wijaya} for models of incompressible viscoelastic solids.  As it is mentioned in \cite{Devedran}, incompressibility 
is a somewhat standard assumption in the setting of 
biological applications. See also \cite{Berjamin} for modeling of
brain tissues. Our interest in the incompressible case is also
motivated by the prospects of devising a sound existence
theory. Assuming incompressibility has the net effect of simplifying
change-of-variable formulas, ultimately allowing the mathematical
treatment.

The  {\it stored energy} of the medium  is assumed to be of
the form
\begin{equation}
   {\mathcal W}(y_{\text{\rm el}},y_{\text{\rm vi}}):=  \int_{y_{\text{\rm vi}}(\Omega)} W_{\text{\rm el}}(\nabla y_{\text{\rm el}}(\xi) )\,{\rm d}\xi + \int_{\Omega}  W_{\text{\rm vi}}(\nabla y_{\text{\rm vi}}(X))\,{\rm d}X. \label{stored}
  \end{equation}
  Here and in the following, we indicate  by $X$  the variable in the reference configuration
 $\Omega$  and  by $\xi$  the variable in the intermediate configuration
 $y_{\rm vi}(\Omega)$. The first integral above corresponds
 to the {\it stored elastic energy} and the given function $W_{\rm
   el}$ is the  stored elastic energy density. Its argument $\nabla
 y_{\text{\rm el}}(\xi) $ can be equivalently rewritten in Lagrangian
 variables as the usual product $\nabla y(X) \nabla y_{\rm
   vi}^{-1}(X)$.  By comparing these two expressions, the
 advantage of working in Eulerian variables is apparent, for $ \nabla
 y_{\text{\rm el}}(\xi) $ is linear in $y_{\text{\rm el}}$. 
The function $W_{\text{\rm vi}}$ is the {\it stored
   viscous energy density} instead and the corresponding integral is
 Lagrangian.

The {\it instantaneous dissipation}  of the system is given by
\begin{equation}
    \Psi(y_{\text{\rm vi}},\dot y_{\text{\rm vi}}):= \int_{\Omega} \psi(\nabla \dot y_{\text{\rm vi}}(\nabla y_{\text{\rm vi}})^{-1})\,{\rm d}X \label{dissi}
\end{equation}
where $\psi(\cdot)$ models the instantaneous dissipation density and is
assumed to be $p_{\psi}$-positively homogeneous  for $p_\psi\geq 2$.

 By formally taking variations of the above introduced
functionals, we obtain the quasistatic equilibrium system 
\cite{mielke2015rate} 
\begin{align*}
  &-\operatorname{div} {\rm D} W_{\rm el}\left(\nabla {y_{\rm
        el}}\right)=f\circ y_{\rm vi}^{-1} \quad  \text{ in } y_{\rm
    vi}(\Omega) \times (0,T) \\[1.5mm]
  &{\rm D} W_{\rm el}\left(\nabla {y_{\rm el}}(y_{\rm vi}\right):{\rm
  D}^2 y_{\rm el}\left(y_{\rm vi})\right)+\operatorname{div} {\rm D}
  W_{\rm vi}\left(\nabla {y_{\rm vi}}\right)  \\
  &\quad \quad\quad-
  \nabla {y_{\rm el}}(y_{\rm vi})^\top f
  =-\operatorname{div}\left({\rm D} \psi\left(\nabla {\dot{y}_{\rm
  vi}}(\nabla {y_{\rm vi}})\right)^{-1}\right)(\nabla {y_{\rm
  vi}})^{-\top}) \quad  \text{ in } \Omega  \times (0,T) .
\end{align*}
 The highly nonlinear character of this system, combined with the
absence of higher-order gradients in the viscous variable, forces us
to consider a suitable weak-solution notion. 

Inspired by \cite[Def. 2.12]{DalMasoLazzaroni} and
\cite[Def. 2.2]{Roger},   in our first main result, Theorem \ref{thm:
  existence limit}, we prove  
the existence  of  {\it approximable solutions} (see
Definition \ref{def:approximable sol}). 
These are everywhere defined  trajectories  $(y_{\rm el}, y_{\rm
  vi}):[0,T]\rightarrow W^{1,p_{\rm el}}(y_{\rm
  vi}(\Omega);\R^d)\times W^{1,p_{\rm vi}}(\Omega;\R^d)$  starting
from some given initial datum $(y_{\text{\rm el},0},y_{\text{\rm
    vi},0})$ and satisfying for all $t\in [0,T]$

\begin{align}
  &\text{\it Energy inequality:}\nonumber\\[2mm]
  &\nonumber
  \int_{y_{\text{\rm vi}}(t,\Omega)} W_{\text{\rm el}}(\nabla
    y_{\text{\rm el}}(t,\xi) )\,{\rm d}\xi + \int_{\Omega}
    W_{\text{\rm vi}}(\nabla y_{\text{\rm vi}}(t,X))\,{\rm d}X-
    \int_{\Omega} f(t,X)\cdot  y_{\text{\rm el}}(t,  y_{\text{\rm
    vi}}(t,X) ) \, {\rm d}X \\
  &\qquad \nonumber  + 
  p_{\psi}\int_{0}^{t}\!\!\int_{\Omega} \psi(\nabla \dot y_{\text{\rm
    vi}}(s,X)(\nabla y_{\text{\rm vi}}(s,X) )^{-1}) \, {\rm d}X\, {\rm d}s\\
  &\quad  \leq \int_{y_{\text{\rm vi},0}(\Omega)} W_{\text{\rm el}}(\nabla
    y_{\text{\rm el},0}(\xi) )\,{\rm d}\xi + \int_{\Omega}
    W_{\text{\rm vi}}(\nabla y_{\text{\rm vi},0}( X))\,{\rm d}X-
    \int_{\Omega} f(0,X)\cdot  y_{\text{\rm el},0}( y_{\text{\rm
    vi},0}(X) ) \, {\rm d}X \nonumber\\
  &\qquad  {}-
  \int_{0}^{t}\int_{\Omega} \partial_s f(s,X)\cdot  y_{\text{\rm el}}(s,  y_{\text{\rm
    vi}}(s,X) ) \, {\rm d}X\, {\rm d}s \label{eneq}\\[4mm]
  &\text{\it Semistability condition:}\nonumber\\[3mm]
  & \int_{y_{\text{\rm vi}}(t,\Omega)} W_{\text{\rm el}}(\nabla
    y_{\text{\rm el}}(t,\xi) )\,{\rm d}\xi  -
    \int_{\Omega} f(t,X)\cdot  y_{\text{\rm el}}(t,  y_{\text{\rm
    vi}}(t,X) ) \, {\rm d}X \nonumber\\
  & \quad  \leq  \int_{y_{\text{\rm vi}}(t,\Omega)} W_{\text{\rm el}}(\nabla
    \tilde{y}_{\text{\rm el}}(\xi) )\,{\rm d}\xi  -
    \int_{\Omega} f(t,X)\cdot  \tilde{y}_{\text{\rm el}} (y_{\text{\rm
    vi}}(t,X) ) \, {\rm d}X \nonumber \\[2mm]
  &\qquad \qquad \qquad
   \forall \tilde{y}_{\text{\rm el}}
  \text { with }\left(\tilde{y}_{\text{\rm el}}, y_{\text{\rm
    vi}}(t,\cdot )\right) \in \mathcal{A} \label{elst}
\end{align}
 where $ \mathcal{A} $ is the set of admissible deformations,
introduced in Section \ref{sec:admin} below. 
The first line of inequality \eqref{eneq} corresponds to the {\it
  total energy} of the medium at time $t$ and state $(y_{\rm
  el}(t,\cdot), y_{\rm
  vi}(t,\cdot))$. In particular, the term $-\int_\Omega f \cdot (y_{\rm
  el}\circ y_{\rm
  vi})\, {\rm d} X$ is the work of the (external) force $f$ (later, a
boundary traction will be considered, as well). 
Solutions $t
\mapsto (y_{\rm el}(t), y_{\rm vi}(t))$ are moreover required to be
{\it approximable}, namely, to ensue as  limit of time
discretizations.  In this respect,  we
consider the {\it incremental minimization problems}, for $i=1,\dots,N$, 
\begin{align}
  &\min_{(y_{\text{\rm el}},y_{\text{\rm vi}})\in \mathcal{A}} \Bigg\{    
    \int_{y_{\text{\rm vi}}(\Omega)} W_{\text{\rm el}}(\nabla
    y_{\text{\rm el}}( \xi) )\,{\rm d}\xi + \int_{\Omega}
    W_{\text{\rm vi}}(\nabla y_{\text{\rm vi}}( X))\,{\rm d}X-
    \int_{\Omega} f(i\tau,X)\cdot  y_{\text{\rm el}}( y_{\text{\rm
  vi}}(X) ) \, {\rm d}X \nonumber\\
 &\qquad \qquad \qquad   +  \tau  \int_{\Omega} \psi\left(\frac{\nabla  y_{\text{\rm
    vi}}(X)- \nabla  y_{\text{\rm
    vi}}^{i-1} (X)}{\tau}(\nabla y_{\text{\rm vi}}^{i-1} (X) )^{-1}\right) \, {\rm d}X\,\Bigg\} \quad \text{ for } y_{\rm vi}^{i-1} \text{ given}\label{discr scheme 1}
\end{align}
on a given uniform time-partition $\{0=t_0<t_1<\dots<t_N=T\}$, where
the set of admissible states $\mathcal{A}$ is defined in Subsection
\ref{sec:admin} below. 

 The   notion of approximable
solution is capable of reproducing the main features of viscoelastic
evolution.
First of all,  the semistability condition
\eqref{elst} implies  that $y_{\rm el}$ solves the
elastic equilibrium at all times, given the viscous-state
evolution. Correspondingly, the description of the purely elastic response of the
material is complete.
 Secondly, the
energy inequality \eqref{eneq} is sharp, in the sense
that it may indeed hold as equality in specific smooth situations. In
other words, all dissipative contributions are correctly taken into
account in  \eqref{eneq}.  Note in this respect the presence of
the factor $p_\psi$ multiplying the dissipation term in  \eqref{eneq}.
Eventually,
the approximation property ensures that viscous evolution actually
occurs, even in absence of applied loads. We give an illustration of
this fact in  Section \ref{sec:res existence} below, see  Figure
\ref{fig:example1}. 

Under
suitable assumptions, the incremental minimization problems \eqref{discr scheme 1} are proved to admit
solutions in Proposition \ref{thm: existence} below. These
time-discrete solutions  fulfill  a discrete energy inequality and a
discrete semistability inequality. The
existence of approximable solutions (Theorem \ref{thm: existence
  limit}) follows by passing to the limit in the time-discrete
problems. 
In order to pass from the  time-discrete  to the 
time-continuous   energy
inequality  \eqref{eneq},  lower semicontinuity of the energy and dissipation
functionals is necessary, which translates in our setting in 
asking for the 
polyconvexity of the respective densities.  In order to obtain the
specific form \eqref{eneq} we need to resort to the notion of
De~Giorgi variational interpolant \cite[Def. 3.2.1, p.~66]{AGS}, and adapt this tool from its
original metric-space application to the current one.

For  establishing the
elastic semistability \eqref{elst},  a suitable recovery-sequence
 construction is required. This calls for the extension of the 
elastic deformations from the intermediate configurations to the whole
$\R^d$.  The possibility of performing this extension requires
some regularity of the boundary of the intermediate configurations,
which we ask to be   
{\it Jones domains} (see Definition~\ref{def:Jones}).

%

 The second main focus of the paper is on the 
rigorous linearization of the system through evolutionary
$\Gamma$-convergence \cite{MRS08} in the case of quadratic
dissipations, namely   for  $p_\psi= 2$.  Moving from the seminal
paper \cite{DalMaso}  in the stationary, hyperelastic case, the
application of $\Gamma$-convergence to inelastic evolutive problems
has been started in \cite{MielkeStefanelli} and has been applied to
different settings.  In particular, linearization in the incompressible
case has been discussed in \cite{Jesenko,Mainini1,Mainini2}. 
The goal is to provide a rigorous formalization of  heuristic
Taylor expansion arguments  which  for the finite-strain
Poynting-Thomson model were already presented in
\cite{PoyntingZerner}. At first, let us review this heuristic by 
assuming sufficient regularity of all ingredients.
Consider the functions $u,\,v,\,w$ defined as
  \begin{equation*}
   u:=\frac{y-\operatorname{id}_{\Omega}}{\varepsilon }, \quad  v:=\frac{y_{\text{\rm vi}}-\operatorname{id}_{\Omega}}{\varepsilon }, \quad \text{ and } \quad w:=\frac{y_{\text{\rm el}}-\operatorname{id}_{y_{\text{\rm vi}}(\Omega)}}{\varepsilon },
  \end{equation*}
  so that  $u,\,v,\,w$ actually correspond to the $\varepsilon$-rescaled
  displacements of $y, \, y_{\rm vi},\, y_{\rm el}$, respectively. 
  To compute the linearization it will be more convenient to  work
  with 
  the pair $(u,v)$  corresponding to the total and viscous
  deformations $(y,y_{\rm vi})$. In particular, we replace
  $\nabla y=I{+}\varepsilon \nabla u$ and $\nabla y_{\rm
    vi}=I{+}\varepsilon \nabla v$ in the stored
  energy and $\nabla \dot y_{\rm
    vi}= \varepsilon \nabla \dot v$ in the dissipation. 
   By formally Taylor expanding  the (rescaled) energy terms
  and  taking  $\varepsilon\rightarrow 0$  we find 
     \begin{align*}
       &\frac{1}{\varepsilon^2}\int_{\Omega}W_{\text{\rm
          el}}\left((I{+}\varepsilon \nabla u)(I{+}\varepsilon \nabla
        v)^{-1}\right)\, {\rm d} X = \int_{\Omega}\frac{1}{2}{\rm D}^2 W_{\text{\rm el}}(I)\nabla (u{-}v): \nabla (u{-}v)  \, {\rm d} X+{\rm o}(\varepsilon)\\
      &\qquad  \to  \frac{1}{2}\int_{\Omega}\nabla (u{-}v): \C_{\text{\rm
          el}}\nabla (u{-}v) \, {\rm d} X ,\\
     &\frac{1}{\varepsilon^2}\int_{\Omega}W_{\text{\rm
       vi}}(I{+}\varepsilon \nabla v) \, {\rm d} X= \int_{\Omega}\frac{1}{2}{\rm
       D}^2 W_{\text{\rm vi}}(I)\nabla v: \nabla v \, {\rm d} X+{\rm
       o}(\varepsilon) \to  
   \frac{1}{2} \int_{\Omega}\nabla v: \C_{\text{\rm vi}}\nabla v \, {\rm d} X
   ,\\
    &\frac{1}{\varepsilon^2}\int_{\Omega}\psi\left(
                                    \varepsilon \nabla
                                    \dot{v}(I{+}\varepsilon
                                    \nabla{v})^{-1}
                                    \right) \, {\rm d} X=\int_{\Omega}\frac{1}{2}{\rm
                                    D}^2\psi(0)\nabla \dot {v}:\nabla
      \dot {v}\, {\rm d} X+{\rm o}(\varepsilon)  \to   \frac{1}{2}\int_{\Omega}\mathbb{D}\nabla \dot {v}:\nabla \dot{v}\, {\rm d} X.
     \end{align*}
      Here, we have assumed $W_{\text{\rm
         el}}(I)=W_{\text{\rm vi}}(I)=0$, ${\rm  D}W_{\text{\rm
         el}}(I)={\rm  D}W_{\text{\rm vi}}(I)=0$, and have defined   $\C_{\text{\rm el}}:={\rm D}^2 W_{\text{\rm el}}(I)$,
     $\C_{\text{\rm vi}}:={\rm D}^2 W_{\text{\rm vi}}(I)$, and
     $\mathbb{D}:={\rm D}^2\psi(0)$. 
   Moreover, we assume that the force $f$ is 
  small, i.e., $ f  = f^{\varepsilon}=\varepsilon
  f^{0}$.  Hence,
   by  neglecting the term  $f^0 \cdot \operatorname{id}_{\Omega}
   $,  which is independent of
  the displacement,  the rescaled  loading term  reads 
  \begin{equation*}    
   - \frac{1}{\varepsilon^2}\int_\Omega  f^\varepsilon  \cdot y_{\rm
      el} \circ y_{\rm vi} \, {\rm d} X =-
    \frac{1}{\varepsilon^2}\int_\Omega  \varepsilon f^0
    \cdot \varepsilon u \, {\rm d} X =-\int_\Omega   f^0
    \cdot  u \, {\rm d} X. 
   \end{equation*}
   
The above  pointwise  convergences are the classical 
heuristic linearization procedure. Still, one is left with actually
checking   that the finite-strain trajectories indeed converge to a
solution of the linearized system. This is the aim of our second
main result, Theorem \ref{thm limit lin}, where we prove that,
given a sequence of approximable solutions $(y_{\rm
  vi,\varepsilon},y_{\rm el,\varepsilon})_{\varepsilon}$ and upon defining
$y_{ \varepsilon}=y_{\rm el,\varepsilon}\circ y_{\rm
  vi,\varepsilon}$ and
the corresponding rescaled displacements $u_\varepsilon=(y_{ \varepsilon}{-} {\rm
  id}_\Omega)/\varepsilon$ and $v_\varepsilon=(y_{\rm
  vi,\varepsilon}{-} {\rm
  id}_\Omega)/\varepsilon$, the sequence
$(u_\varepsilon,v_\varepsilon)_\varepsilon$ converges pointwise in
time (up to
subsequences) to  
$(u,v):[0,T]\rightarrow H^1(\Omega;\R^d)\times H^1(\Omega;\R^d)$  with
$(u(0),v(0)) = (u^{0}{,}v^{0}) :=\lim_{\varepsilon \to 0}
(u_\varepsilon(0),v_\varepsilon(0))$
and satisfying, for all $t\in[0,T]$,
\begin{align}
  &\text{\it Linearized energy inequality:}\nonumber\\[2mm]
  & \frac{1}{2}\int_{\Omega}\nabla (u(t){-}v(t)): \C_{\text{\rm
          el}}\nabla (u(t){-}v(t))\,{\rm d} X +  \frac{1}{2} \int_{\Omega}\nabla v(t): \C_{\text{\rm vi}}\nabla v(t) \, {\rm d} X- \int_\Omega   f^0(t)
    \cdot  u(t) \, {\rm d} X \nonumber\\
  &\qquad  +  \int_0^t\!\!\!\int_{\Omega}\mathbb{D}\nabla \dot
    {v}(s):\nabla \dot{v}(s) \, {\rm d} X\, {\rm d}s\nonumber\\
  & \quad 
  \leq \frac{1}{2}\int_{\Omega}\nabla (u_0{-}v_0): \C_{\text{\rm
          el}}\nabla (u_0{-}v_0)\,{\rm d} X +  \frac{1}{2} \int_{\Omega}\nabla v_0: \C_{\text{\rm vi}}\nabla v_0 \, {\rm d} X - \int_\Omega   f^0(0)
    \cdot  u_0 \, {\rm d} X \nonumber\\
  &\qquad  -
   \int_{0}^{t}\!\!\! \int_\Omega \partial_s f^0(s)\cdot  u(s) \, {\rm
    d } X\, {\rm d}s,\label{eneq0}\\[4mm]
&\text{\it Linearized semistability:}\nonumber\\[3mm] 
  & \frac{1}{2}\int_{\Omega}\nabla (u(t){-}v(t)): \C_{\text{\rm
          el}}\nabla (u(t){-}v(t))\,{\rm d} X   - \int_\Omega   f^0(t)
    \cdot  u(t) \, {\rm d} X \nonumber\\
  &\quad  \leq \frac{1}{2}\int_{\Omega}\nabla (\hat{u}{-}v(t)): \C_{\text{\rm
          el}}\nabla (\hat{u}{-}v(t))\,{\rm d} X  - \int_\Omega   f^0(t)
    \cdot  \hat{u} \, {\rm d} X \quad \forall \hat{u} \ \ \text{admissible}. \label{elst0}
  \end{align}

 The linearized energy inequality and the linearized 
semistability deliver  a weak notion of solution for the
linearized problem  \eqref{linear PDE1}-\eqref{linear
  PDE2}.  Albeit  \eqref{eneq0}-\eqref{elst0} are too weak to fully
characterize the unique solution of linearized Poynting-Thomson
system \eqref{linear PDE1}-\eqref{linear
  PDE2}, the equilibrium system \eqref{linear PDE1} is fully
recovered. In particular, $u$ is uniquely determined at all times,
given $v$. Moreover, the linearized energy equality \eqref{eneq0} is
sharp and turns  out  to be an equality in specific cases.

To conclude, let us note that one could alternatively perform the linearization at the time-discrete level
and then pass to the time-continuous limit. This way one recovers the unique
strong solution of the linearized  Poynting-Thomson system \eqref{linear PDE1}- \eqref{linear
  PDE2}. This fact provides
some additional justification of the finite-strain model. Still, we do
not follow hier this alternative path, which could be easily treated
along the lines of the analysis in 
 of Sections \rm{\ref{sec:Proof of Existence}} and
\rm{\ref{sec:Proof of Linearization}}.

\section{Preliminaries}
	\label{sec:pre}
We devote this section to  setting notation and presenting some
preliminary results. 

\subsection{Notation}

  In what follows, we denote by $\R^{d\times d}$ the Euclidean space of $d\times d$ real matrices,  $d\ge 2$. Given $A\in \R^{d\times d}$ we define its (Frobenius)
  norm as $
    |A|^2 := A:A$,
  where the contraction product between $2$-tensors is defined as $A {:}
   B  := \sum_{i,j}A_{ij}B_{ij}$.
  Moreover, given a symmetric positive definite 4-tensor $\C\in\R^{d\times d\times d\times d}$, the corresponding induced matrix norm is defined as $
  |A|_{\C}^2 :=  \C A:A/2 $. We denote by $A^{\rm sym}:=(A+A^{\top})/2$ the symmetric part of a matrix $A\in \R^{d\times d}$.
  We shall use the following matrix sets
  \begin{align*}
   & SL(d):=\{ A \in \R^{d\times d} \;|\; \det  A
   =1\},\\
   & SO(d):=\{ A \in SL(d) \;|\; A A^{\top}
   =I\},\\
  & GL(d):=\{ A  \in \R^{d\times d} \;|\; \det  A \neq 0 \},\\
  & GL_{+}(d):=\{ A  \in \R^{d\times d} \;|\; \det  A >0 \}.
  \nonumber
  \end{align*}
  The symbol $B_r(A)\subset \R^{d\times d}$ denotes the open ball of radius $r>0$ and center $A\in\R^{d\times d}$.
  We make use of the function spaces
  \begin{align*}
    &H^{1}_{\sharp}(\Omega;\R^d):=\Big\{ u\in H^{1}(\Omega;\R^d)\;\Big|\; \int_{\Omega} u \,{\rm d}X=0 \Big\},\\
      &H^{1}_{\Gamma}(\Omega;\R^d):=\{ u\in H^{1}(\Omega;\R^d)\;|\; u=0 \text{ on } \Gamma\subset\partial\Omega \},
  \end{align*}
  where $\Gamma$ is a nonempty, open, and measurable subset of $\partial \Omega$.

  Moreover we denote by $\mathcal{H}^{d-1}$ the $(d-1)$-dimensional Hausdorff measure
  and by $|\omega|$ the $d$-dimensional Lebesgue measure of the measurable set $\omega$ in $\R^d$.

\subsection{Deformations and admissible states}\label{sec:admin}

Let  fix the reference configuration $\Omega$ of the body to  be a nonempty, open, bounded, and connected Lipschitz
subset of $\R^d$.
We assume without loss of generality that $\Omega$ is such that $\int_{\Omega}X\,{\rm d} X=0$.
 We let $\Gamma_D,\Gamma_N$ be open subsets of $\partial \Omega$ (in
 the relative topology  of $\partial\Omega$)  such that $\overline{\Gamma}_D\cup \overline{\Gamma}_N=\partial \Omega$, $\mathring{\Gamma}_D\cap \mathring{\Gamma}_N=\emptyset$, and $\mathcal{H}^{d-1}(\Gamma_D)>0$.

The viscous deformation is required to fulfill
\begin{equation*}
    y_{\text{\rm vi}}\in W^{1,p_{\text{\rm vi}}}(\Omega;\R^d) \quad \text{ for some }\quad p_{\text{\rm vi}}>d(d-1)
\end{equation*}
and to be locally volume-preserving, i.e., $\det \nabla y_{\text{\rm vi}}=1$ almost everywhere in $\Omega$. In the following, $y_{\text{\rm vi}}$ is tacitly identified with its H\"older-continuous representative. More precisely, $y_{\text{\rm vi}}\in C^{0,1-d/p_{\text{\rm vi}}}(\Omega;\R^d)$ and is almost everywhere differentiable (see \cite{fonseca_gangbo_2018}).
In addition, since we will use the change-of-variables formula to pass from Lagrangian to Eulerian variables, we require $y_{\text{\rm vi}}$ to be injective almost everywhere. Equivalently, we ask for the {\it Ciarlet-Ne\v{c}as condition} \cite{Ciarlet Necas}
\begin{equation}\label{CN}
    |\Omega|=\int_{\Omega}\det \nabla y_{\text{\rm vi}} \,{\rm d}X = |y_{\text{\rm vi}}(\Omega)|
\end{equation}
to hold.
As a consequence we have the change-of-variables formula
\begin{equation*}
    \int_{\omega}\varphi (y_{\text{\rm vi}}(X)) \,{\rm d}X=\int_{y_{\text{\rm vi}}(\omega)}\varphi(\xi)\,{\rm d}\xi
\end{equation*}
for every measurable set $\omega\subseteq \Omega$ and every measurable function $\varphi:y_{\text{\rm vi}}(\omega)\rightarrow \R^d$.
Note that $y_{\rm vi}\in W^{1,p}(\Omega;\R^d)$ has {\it distortion}
$K := |\nabla y_{\rm vi}|^d/\det \nabla y_{\rm vi}=|\nabla y_{\rm
  vi}|^d \in L^ {p_{\text{\rm vi}}/d}(\Omega; \R)$, since it is
locally volume preserving.  As $p_{\text{\rm vi}}/d > d - 1$, this
bound on the distortion $K$  implies that $y_{\rm vi}$ is either constant or open \cite[Theorem 3.4]{distorsion}. By the Ciarlet-
Ne\v cas condition \eqref{CN}, $y_{\rm vi}$ cannot be constant, and hence $y_{\rm vi}$ is open.
In particular $y_{\rm vi}(\Omega)$ is an open set.
 Moreover, we also have that $y_{\rm vi}$ is (globally) injective \cite[Lemma 3.3]{Grandi}, and that $y_{\rm vi}$
is actually a homeomorphism with inverse $y_{\text{\rm vi}}^{-1}\in W^{1,p_{\text{\rm vi}}/(d-1)} (y_{\rm vi}(\Omega);\R^d)$ (see \cite{fonseca_gangbo_2018}).

In order to make the statement of the model precise, we need to
require some regularity of the intermediate configuration $y_{\rm
  vi}(\Omega)$. We recall the following definition.
\medskip

\begin{definition}[$(\eta_1,\eta_2)$-Jones domain \cite{Jonesdomain}]\label{def:Jones}
Let $\eta_1,\eta_2>0$. A bounded open set $\omega \subset \mathbb{R}^{d}$ is said to be a \emph{$(\eta_1,\eta_2)$-Jones domain}, if for every $x, y \in \omega$ with $|x-y|<\eta_2$ there exists a Lipschitz curve $\gamma \in W^{1, \infty}([0,1] ; \omega)$ with $\gamma(0)=x$
and $\gamma(1)=y$ satisfying the following two conditions:
$$
l(\gamma):=\int_{0}^{1}|\dot{\gamma}(s)| \mathrm{d} s \leq \frac{1}{\eta_1}|x-y|
$$
and
$$
d(\gamma(t), \partial \omega) \geq \eta_1 \frac{|x-\gamma(t)||\gamma(t)-y|}{|x-y|} \quad \text{ for every } \quad t \in[0,1].
$$
The set of $(\eta_1,\eta_2)$-Jones domains will be  denoted by  $\mathcal{J}_{\eta_1,\eta_2}$.
\end{definition}
In the following, we will exploit the fact that $(\eta_1,\eta_2)$-Jones domains are {\it Sobolev extension domains}: for all $\eta_1,\eta_2>0$,  $p\in [1,\infty)$, and all $\omega\in \mathcal{J}_{\eta_1,\eta_2}$ there exists a positive constant $C=C(\eta_1,\eta_2, p, \omega, d)$ and a linear operator
  $E:W^{1,p}(\omega;\R^d)\rightarrow W^{1,p}(\R^d;\R^d)$
such that $Ey=y$ on $\omega$ and
\begin{equation*}
  \|Ey\|_{W^{1,p}(\R^d)}\leq C\|y\|_{W^{1,p}(\omega)} \quad \text{ for every } y\in W^{1,p}(\omega;\R^d).
\end{equation*}
Note that the class of $(\eta_1,\eta_2)$-Jones domains is closed under Hausdorff convergence \cite{Melching}.
 In the following, we will need to consider extensions and we then
ask for the regularity   
\begin{equation*}
   y_{\text{\rm vi}}(\Omega)\in \mathcal{J}_{\eta_1,\eta_2}.
\end{equation*}

Finally, since the problem will be formulated only in terms of the gradient of $y_{\text{\rm vi}}$, we impose the normalisation condition
\begin{equation}\label{av0}
    \int_{\Omega}y_{\text{\rm vi}} \,{\rm d}X=0.
\end{equation}

Given a viscous deformation $y_{\text{\rm vi}}$, we assume the elastic deformation to fulfill
\begin{equation*}
    y_{\text{\rm el}}\in W^{1,p_{\text{\rm el}}}(y_{\text{\rm vi}}(\Omega);\R^d) \quad \text{ for some } \quad p_{\text{\rm el}}>d
\end{equation*}
and we tacitly identify $y_{\rm el}$ with its H\"older-continuous representative.

 For all given   viscous deformation $y_{\text{\rm vi} }:\Omega \rightarrow \R^d$ and elastic deformation $y_{\text{\rm el}}:y_{\text{\rm vi}}(\Omega)^{\circ}\rightarrow \R^d$, we define the total deformation as the composition of the two, i.e.,
\begin{equation*}
   y:=y_{\text{\rm el}}\circ y_{\text{\rm vi}}:\Omega \rightarrow  \R^d.
\end{equation*}
We assume that $y$ satisfies  a  Dirichlet boundary condition
on $\Gamma_D$,  namely, 
\begin{equation}\label{DirBC}
    y=\operatorname{id} \quad \text{ on } \Gamma_D.
\end{equation}
Since $y_{\text{\rm vi}}$ is invertible and both $y_{\rm vi}$ and $y_{\rm el}$ are almost everywhere differentiable, the following chain rule
\begin{equation*}
    \nabla y(X)=\nabla y_{\text{\rm el}}(y_{\text{\rm vi}}(X)) \nabla y_{\text{\rm vi}}(X)
\end{equation*}
holds for almost every $X\in \Omega$. 
Hence, $y$ satisfies
\begin{equation*}
    \| \nabla y \|_{L^q(\Omega)}\leq \| \nabla y_{\text{\rm el}} \|_{L^{p_{\text{\rm el}}}(y_{\text{\rm vi}}(\Omega))}\| \nabla y_{\text{\rm vi}} \|_{L^{p_{\text{\rm vi}}}(\Omega)} \quad \text{ where }\quad \frac{1}{q}:=\frac{1}{p_{\text{\rm el}}}+\frac{1}{p_{\text{\rm vi}}},
\end{equation*}
as can be readily checked by a change of variables and by the H\"older inequality. In particular, the boundary condition \eqref{DirBC} should be understood in the classical trace sense.

To sum up, the set of {\it admissible states} is defined as
\begin{align}
       \mathcal{A}:=\Bigg\{(y_{\text{\rm el}},y_{\text{\rm vi}})&\in W^{1,p_{\text{\rm el}}}(y_{\text{\rm vi}}(\Omega);\R^d)\times W^{1,p_{\text{\rm vi}}}(\Omega;\R^d) \;\Bigg|\;
      \det \nabla y_{\text{\rm vi}}=1 \text{ a.e. in } \Omega,\notag \\
      &\int_{\Omega}y_{\text{\rm vi}}\, {\rm d}X=0, \;|\Omega|= |y_{\text{\rm vi}}(\Omega)|,\; y_{\text{\rm vi}}(\Omega)\in \mathcal{J}_{\eta_1,\eta_2},\; y=y_{\text{\rm el}}\circ y_{\text{\rm vi}}=\operatorname{id} \text{ on } \Gamma_{D}\Bigg\}.\notag
\end{align}
 Viscoelastic states are naturally depending on time. 
From now on, we are  hence  interested in {\it trajectories}
$(y_{\rm el}, y_{\rm vi}):[0,T]\rightarrow \mathcal{A}$  in the
set of admissible states.

\section{Main results}\label{sec:main}

We devote this section to the statements of our assumptions and our main results.

\subsection{Assumptions for the existence theory}\label{sec:hp existence}
In this section we specify the assumptions needed for the existence
results,  namely, Proposition \ref{thm: existence} and Theorem \ref{thm: existence limit}. 

The {\it total energy} of the system  at time $t\in[0,T]$ and state $(y_{\text{\rm el}},y_{\text{\rm vi}})\in \mathcal{A}$  is given by
\begin{equation*}
  \mathcal{E}(t,y_{\text{\rm el}},y_{\text{\rm vi}}):= \mathcal{W}(y_{\text{\rm el}},y_{\text{\rm vi}})-\langle \ell(t),y_{\text{\rm el}}\circ y_{\text{\rm vi}} \rangle,
\end{equation*}
where $\mathcal{W}(y_{\text{\rm el}},y_{\text{\rm vi}})$ is the {\it
  stored energy} and the pairing $\langle \ell(t),y_{\text{\rm
    el}}\circ y_{\text{\rm vi}} \rangle$ represents the work of {\it external mechanical actions}.

More precisely, the stored energy is defined as
\begin{equation*}
    \mathcal{W}(y_{\text{\rm el}},y_{\text{\rm vi}}):=\int_{y_{\text{\rm vi}}(\Omega)} W_{\text{\rm el}}(\nabla y_{\text{\rm el}}(\xi) )\,{\rm d}\xi + \int_{\Omega}  W_{\text{\rm vi}}(\nabla y_{\text{\rm vi}}(X))\,{\rm d}X
\end{equation*}
where $W_{\text{\rm el}}:\R^{d\times d}\rightarrow \R $ and
$W_{\text{\rm vi}}:\R^{d\times d}\rightarrow \R \cup \{\infty\}$  are
the stored {\it elastic} and  the  stored {\it viscous} energy densities, respectively. On the energy densities we assume that:
\begin{enumerate}[label=(E\arabic*)]
  \item\label{E1} there exist positive constants $c_1,c_2$ such that
  \begin{equation}\label{growthel}
    c_1|A|^{p_{\text{\rm el}}}\leq W_{\text{\rm el}}(A)\leq \frac{1}{c_1}(1+|A|^{p_{\text{\rm el}}}) \quad \text{ for every }  A\in GL(d),
  \end{equation}
  \begin{equation}\label{growthvi}
    W_{\text{\rm vi}}(A)\geq \begin{cases}
  c_2|A|^{p_{\text{\rm vi}}}- \frac{1}{c_2} \quad &\text{ for every } A\in SL(d)\\
  \infty &\text{ otherwise,}
  \end{cases}
  \end{equation}
  for $p_{\rm el}>d$ and $p_{\rm vi}>d(d-1)$.
  \item\label{E2} $ W_{\text{\rm el}}, W_{\text{\rm vi}}$ are polyconvex, i.e.,
   there exist two  convex functions $\hat{W}_{\text{\rm el}}, \hat{W}_{\text{\rm vi}}:\R^{\zeta (d)}\rightarrow \R \cup \{\infty \}$ such that
  \begin{equation*}
    W_{\text{\rm el}}(A)=\hat{W}_{\text{\rm el}}(T(A))\quad \text{ and } \quad W_{\text{\rm vi}}(A)=\hat{W}_{\text{\rm vi}}(T(A))
  \end{equation*}
  where the minors $T(A)$ of $A$ are given by  $T:\R^{d\times d}\rightarrow \R^{\zeta(d)}$
  \begin{equation*}
    T(A):=(A, \operatorname{adj}_2 A, \cdots ,\operatorname{adj}_{d} A ).
  \end{equation*}
  Here, $\operatorname{adj}_s A$ denotes the matrix of all minors $s\times s$ of the matrix $A \in \R^{d\times d}$, for $s=2,\cdots ,d$ and  $\zeta(d):=\sum_{s=1}^d\binom{s}{d}^2$.
\end{enumerate}

Notice that, since $p_{\text{\rm vi}}>d$, the mapping $ y_{\text{\rm vi}}\mapsto \operatorname{adj}_{s}\nabla y_{\text{\rm vi}}$ is $(W^{1,p_{\text{\rm vi}}}, L^{ p_{\text{\rm vi}}/s})$-weakly sequentially  continuous.
Hence, given $y_{\text{\rm vi},n}\rightharpoonup y_{\text{\rm vi}} $ in $W^{1,p_{\text{\rm vi}}}(\Omega;\R^d)$ with $\det \nabla y_{\text{\rm vi},n}=1$ almost everywhere in $\Omega$, we have that
\begin{equation*}
  1=\det\nabla y_{\text{\rm vi},n}\rightharpoonup \det\nabla
  y_{\text{\rm vi}}= 1 \quad  \text{in} \ L^{ p_{\text{\rm vi}}/d}
   (\Omega).
\end{equation*}
As  $ \nabla y_{\text{\rm vi}}(X)\in SL(d)$ a.e. in $\Omega$, we have that
\begin{equation*}
  \int_{\Omega}  W_{\text{\rm vi}}(\nabla y_{\text{\rm vi}}(X))\,{\rm d}X\leq \liminf_{n} \int_{\Omega}  W_{\text{\rm vi}}(\nabla y_{\text{\rm vi},n}(X))\,{\rm d}X
\end{equation*}
by polyconvexity of $W_{\rm vi}$. In particular,
$y_{\text{\rm vi}}\mapsto\int_{\Omega}  W_{\text{\rm vi}}(\nabla y_{\text{\rm vi}}(X))\,{\rm d}X$ is weakly lower semicontinuous in $W^{1,p_{\text{\rm vi}}}(\Omega;\R^d)$.

The growth condition \eqref{growthvi} ensures that all viscous
deformations $y_{\rm vi}$ of finite energy are incompressible. 
Local elastic incompressibility $\det \nabla
y_{\text{\rm el}}=1$ or even the weaker $\det \nabla
y_{\text{\rm el}}>0$   cannot be required, however. This is due to the
fact that we later need  to consider the Sobolev extension of
$y_{\text{\rm el}}$ from the moving domain $y_{\text{\rm vi}}(\Omega)$
to $\R^d$ in order to compute the limit of an infimizing
sequence.  As it is well-known, such extensions  may not
preserve the positivity of    $\det \nabla
y_{\text{\rm el}}$. 

 On the other hand,  our assumptions on the elastic energy
density  are compatible with {\it frame indifference}. In
particular, we could ask  $W_{\rm el}(RA)=W_{\rm el}(A)$ for every
rotation $R\in SO(d) $ and every $A\in \R^{d\times d}$. Note
nonetheless that this  property,  although fundamental from the
mechanical standpoint,  is  actually not  needed for the
analysis.  The above assumptions would be compatible with
requiring that $W_{\rm vi}$ is invariant by left multiplication with
special rotations, as well. Still, such an
invariance would be little relevant from the modeling viewpoint, for 
 the viscous energy density is defined on viscous
deformations, which take values in the intermediate
configuration. 

Eventually, the work of external mechanical actions is assumed to
result from a  given  time-dependent {\it body force}  $f:[0,T]\times
\Omega\rightarrow \R^d$ and a  given   time-dependent boundary
{\it traction} $g:[0,T]\times \Gamma_{N}\rightarrow \R^d$ as follows
\begin{equation}\label{elle}
    \langle \ell(t),y \rangle := \int_{\Omega} f(t,X)\cdot y(X)\,{\rm d}X + \int_{\Gamma_{N}}g(t,X)\cdot y(X)\,{\rm d}\mathcal{H}^{d-1}(X).
\end{equation}
We assume
\begin{enumerate}[label=(E\arabic*)]\setcounter{enumi}{2}
  \item\label{E3}
  $f\in W^{1,\infty}(0,T;L^{(q^*)'}(\Omega;\R^d))$ and $g\in W^{1,\infty}(0,T;L^{(q^\#)'}(\Gamma_{N};\R^d))$ where $q^*$ and $q^\#$ are the Sobolev and trace exponent related to $W^{1,q}(\Omega;\R^d)$,  respectively (see \cite{Roubicek}) and the prime denotes conjugation.
\end{enumerate}
Consequently, we have
\begin{equation*}
  \ell \in W^{1,\infty}{\left(0,T;(W^{1,q}(\Omega;\R^d))^*\right)}{,}
\end{equation*}
where $(W^{1,q}(\Omega;\R^d))^*$ is the dual space of $W^{1,q}(\Omega;\R^d)$.


Given a time-dependent viscous trajectory $y_{\rm
  vi}\,{:}\,[0,T]\rightarrow W^{1,p_{\rm vi}}(\Omega;\R^d)$, we define
the {\it total instantaneous dissipation} of the system \cite{MielkeRossiSavare} as
\begin{equation}\label{dissipation}
  \Psi( y_{\text{\rm vi}}, \dot{y}_{\text{\rm vi}} ):=\int_{\Omega} \psi(\nabla \dot y_{\text{\rm vi}}(\nabla y_{\text{\rm vi}})^{-1})\,{\rm d}X.
\end{equation}
Here and in the following, the dot represent a partial derivative with
respect to time. 
Above, the dissipation density $\psi:\R^{d\times d}\rightarrow [0,\infty)$ is assumed to be:
\begin{enumerate}[label=(E\arabic*)]\setcounter{enumi}{3}
  \item\label{E4}  convex and differentiable at $0$ with $\psi(0)=0$;
  \item\label{E5}  fulfilling
  \begin{equation}\label{dissipation density}
    \psi(A)\geq c_3|A|^{p_{\psi}} \quad \text{ for every } A\in \R^{d\times d}
  \end{equation}
  for some positive constant $c_3$;
 \item\label{E6} positively $p_{\psi}$-homogeneous, namely 
 \begin{equation}\label{p-homogeneous}
   \psi(\lambda A)=\lambda^{p_{\psi}}\psi(A) \quad \text{ for every } A\in \R^{d\times d}, \lambda \geq 0.
 \end{equation}
\end{enumerate}
 The form of the instantaneous dissipation is
 parallel to  the analogous definition in 
elastoplasticity,  where
nonetheless $\psi$ is assumed to be positively $1$-homogeneous, namely
$p_\psi=1$ \cite{Mielke02,Mielke03}. In particular, let us explicitly
point out that it does not fall within the frame-indifferent setting 
from \cite{Antman}. Indeed, in this case  viscous deformations take values in the
  intermediate configuration only and frame-indifference should not
  necessarily be imposed  there.

In the following, we ask
\begin{equation}\label{p_psi>2}
  p_{\psi}\geq 2 \geq \frac{d(d-1)}{d(d-1)-1} 
\end{equation}
 where we have used  $d\geq 2$.  In particular, we have
that  $p_{\psi}'<p_{\psi}$ and, by defining $p_r$ by $1/p_r:={1}/{p_{\psi}}+{1}/{p_{\text{\rm vi}}}$, one has that $p_r>1$. Again by H\"older's Inequality, this entails that
\begin{equation*}
  \|\nabla \dot {y}_{\text{\rm vi}}\|_{L^{p_r}(\Omega)}\leq
\|\nabla \dot {y}_{\text{\rm vi}}(\nabla y_{\text{\rm vi}})^{-1}\|_{ L^{p_{\psi}}(\Omega)} \|\nabla y_{\text{\rm vi}}\|_{L^{p_{\rm vi}}(\Omega)}\leq
c\Psi(y_{\text{\rm vi}}, \dot{y}_{\text{\rm vi}})\left(\mathcal{W}(y_{\text{\rm el}},y_{\text{\rm vi}})^{1/p_{\rm vi}}+1\right).
\end{equation*}
In particular, $\nabla \dot{y}_{\text{\rm vi}}$ belongs to
$L^{p_r}(\Omega;\R^{d\times d})$ with $p_r>1$ whenever energy and
dissipation are finite.

Here and in the following, the symbol $c$ denotes a generic positive constant, possibly depending on data and changing from line to line.


\subsection{Existence results}\label{sec:res existence}

Before presenting the statements of our main results, we make the notion of solution to the problem precise. To this aim, let
$\Pi_{\tau}:=\{ 0=t_0<t_1<...<t_{N}=T  \}$ denote the uniform partition of the time interval $[0,T]$ with time step $t_{i}-t_{i-1}=\tau>0$ for every $i=1,\dots, N_:=T/\tau$.
From now on, let $(y_{\text{\rm el},0},y_{\text{\rm vi},0})$ be a
compatible initial condition, i.e.,
\begin{equation}\label{initial}
  (y_{\text{\rm
    el},0},y_{\text{\rm vi},0})\in \mathcal{A} \ \  \text{with} \ \
\E(0,y_{\text{\rm el},0},y_{\text{\rm vi},0})<\infty.
\end{equation}
Given $(y_{\text{\rm el}}^0,y_{\text{\rm vi}}^0):=(y_{\text{\rm
    el},0},y_{\text{\rm vi},0})$, for  all  $i=1,...,N$ we define the incremental minimization problems
\begin{equation}\label{i min probl}
  \min_{(y_{\text{\rm el}},y_{\text{\rm vi}})\in \mathcal{A}} \left\{ \mathcal{E}(t_i,y_{\text{\rm el}},y_{\text{\rm vi}})
  + \tau \Psi\left(y^{i-1}_{\rm vi},\frac{y_{\rm vi}-y^{i-1}_{\rm vi}}{\tau}   \right) \right\}.
\end{equation}
We call a sequence of minimizers $(y_{\rm el}^i,y_{\rm vi}^i)_{i=0}^{N}$ of \eqref{i min probl} an {\it incremental solution} of the problem corresponding to time step $\tau$.

Note that incremental solutions exist. In particular, we have the
following.
\medskip

\begin{proposition}[Existence of incremental solutions]\label{thm: existence}

  Under assumptions {\rm \ref{E1}, \ref{E2}, \ref{E3}, \ref{E4}}, and 
  {\rm \ref{E5}} of Section {\rm \ref{sec:hp existence}}  and
  \eqref{initial} 
  the incremental minimization problem \eqref{i min probl} admits an incremental solution $(y_{\text{\rm el}}^{i},y_{\text{\rm vi}}^{i})_{i=0}^{N}\subset \mathcal{A} $.
\end{proposition}

 The proof of Proposition \ref{thm: existence} is given in Section
\ref{sec:time discr}. 

In the following, we make use of the following notation for interpolations. Given a vector $(u_{0},...,u_{N})$,
we define its backward-constant interpolant $\overline{u}_{\tau}$, its forward-constant interpolant $\underline{u}_{\tau}$, and its piecewise-affine interpolant $\hat{u}_{\tau}$ on the partition $\Pi_{\tau}$ as
  \begin{align*}
    &\overline{u}_\tau(0):=u_0, \quad \quad \overline{u}_\tau(t):=u_{i}  &\text{ if } t\in(t_{i-1},t_i] \quad\text{ for } i=1,\dots, N,\\
    &\underline{u}_{\tau}(T):=u_{N}, \quad \; \underline{u}_\tau(t):=u_{i-1} &\text{ if } t\in[t_{i-1},t_i) \quad\text{ for } i=1,\dots, N,\\
    &\hat{u}_\tau(0):=u_0, \quad \quad \hat{u}_\tau(t):=\frac{u_i-u_{i-1}}{t_i-t_{i-1}}(t-t_{i-1})+u_{i-1}
    &\text{ if } t\in(t_{i-1},t_i] \quad\text{ for } i=1,\dots, N.
  \end{align*}
We are now in the position of  introducing  our notion of
solution to the large-strain Poynting-Thomson model.

\medskip

\begin{definition}[Approximable solution]\label{def:approximable sol}
  We call $(y_{\rm el}, y_{\rm vi}):[0,T]\rightarrow \mathcal{A}$ an \emph{approximable solution} if
  there exist a sequence of uniform partitions 
  of the interval $[0,T]$ with mesh size $\tau\to 0$, corresponding incremental solutions $(y_{\rm el}^i,y_{\rm vi}^i)_{i=0}^{N}$, and a nondecreasing function $\delta:[0,T]\rightarrow[0,\infty )$
 such that, for every $0\leq s \leq t \leq T$,
 \begin{align}
   &\text{\rm Approximation:}\notag\\[1mm]
   &(\overline{y}_{\text{\rm
                                  el},\tau}(t),\overline{y}_{\text{\rm
                                  vi},\tau}(t))\rightharpoonup
                                  (y_{\text{\rm el}}(t),y_{\text{\rm
                                  vi}}(t)) \  \text{in} \
 W_{\rm loc}^{1,p_{\text{\rm el}}}( y_{\text{\rm vi}}(t,\Omega);\R^d)\times W^{1,p_{\text{\rm vi}}}(\Omega;\R^d),\notag\\
     &\int_{0}^{t}\Psi\left(\underline{y}_{{\rm vi},\tau},\dot{\hat{y}}_{{\rm vi},\tau}   \right) \rightarrow \delta(t), \notag \\
     &\int_{s}^{t}\Psi\left(\underline{y}_{{\rm vi},\tau},\dot{\hat{y}}_{{\rm vi},\tau}   \right) \leq \delta(t)-\delta(s),\notag\\[3mm]
   &\text{\rm Energy inequality:}\notag\\
   \quad&\mathcal{E}(t,y_{\text{\rm el}},y_{\text{\rm vi}})+
   p_{\psi}\delta(t)\leq \E(0,y_{\text{\rm el},0},y_{\text{\rm vi},0})-
   \int_{0}^{t} \langle\dot{\ell}(s),y \rangle, \label{energy inequality} \\[2mm]
  & \text{\rm Semistability:}\notag\\[2mm]
   \quad&\mathcal{E}\left(t, y_{\text{\rm el}}(t), y_{\text{\rm vi}}(t)\right) \leq \mathcal{E}\left(t,
   \tilde{y}_{\text{\rm el}}, y_{\text{\rm vi}}(t)\right) \quad \forall \tilde{y}_{\text{\rm el}}
   \  \text {with} \ \left(\tilde{y}_{\text{\rm el}}, y_{\text{\rm vi}}(t)\right) \in \mathcal{A}.\label{semistability}
 \end{align}
\end{definition}

Our first main result concerns the existence of approximable
solutions.

\medskip

\begin{theorem}[Existence of approximable solutions]\label{thm: existence limit}
    Under the assumptions {\rm \ref{E1}, \ref{E2}, \ref{E3}, \ref{E4},
      \ref{E5}}, and {\rm \ref{E6}} of Section {\rm \ref{sec:hp
        existence}}  and
  \eqref{initial} 
there exists an approximable solution $(y_{\text{\rm el}},y_{\text{\rm vi}}):[0,T]\rightarrow \mathcal{A} $.
\end{theorem}


 The proof of Theorem \ref{thm: existence limit} is detailed in Section \ref{sec:Proof of Existence}.

\begin{figure}[h!]  
  \centering
  \includegraphics[width=0.5\textwidth]{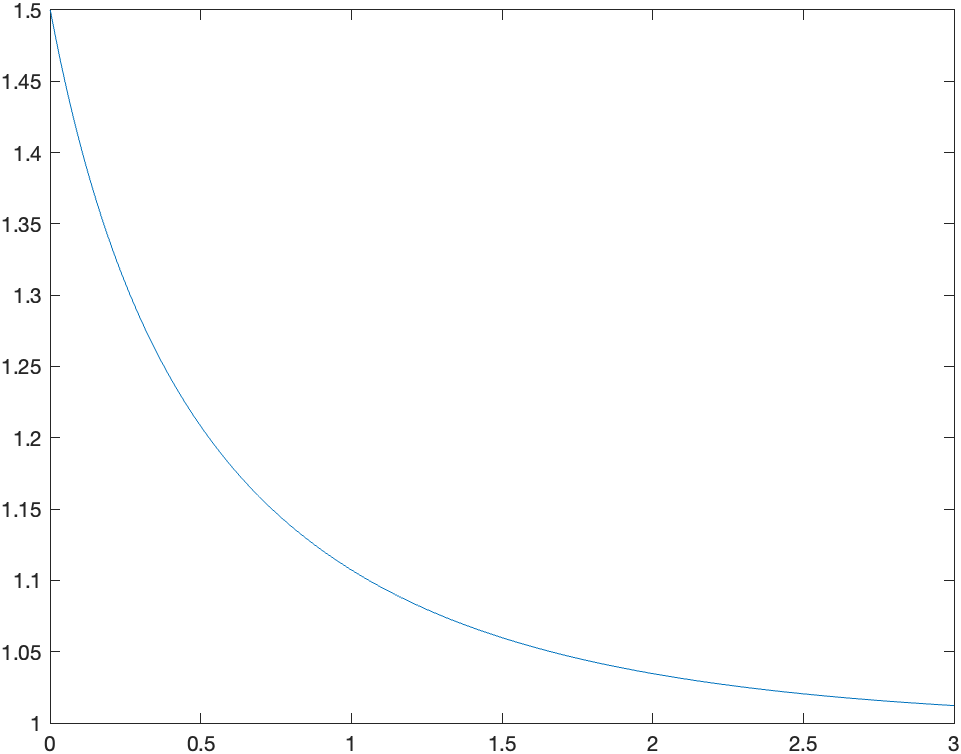}
  \caption{Evolution of the viscous strain  $t \in [0,3] \mapsto
    F_{\rm vi}(t)$ in the limit $\tau \to 0$ from problem
    \eqref{example},  starting from  $F_{\rm
      vi}^0=1.5$.  }
      \label{fig:example1}
\end{figure}

 As already mentioned in the introduction, the fact that solutions
are {\it approximable} ensures that viscous evolution actually
occurs,  even in absence of applied loads. 
We show this fact by resorting to the simplest, scalar model at  a
single material point. We consider the energy densities, the
dissipation to be quadratic, and that no loading is present. 
More precisely, we let $F\in \R$  and $F_{\rm vi}>0$
represent the  total and viscous (scalar) strains, respectively,
we define $W_{\rm el}(F_{\rm el})=W_{\rm el}(FF_{\rm
  vi}^{-1}):=\frac{1}{2}|FF_{\rm vi}^{-1}-1|^2$,  $W_{\rm vi}(F_{\rm
  vi}):=\frac{1}{2}|F_{\rm vi}-1|^2$, $\psi(\dot{F}_{\rm vi}F_{\rm
  vi}^{-1}):=\frac{1}{2}|\dot{F}_{\rm vi}F_{\rm vi}^{-1}|^2$, and 
we let
$\ell(t)\equiv 0$ for every $t\in [0,T]$. In this setting, the
discrete incremental problem \eqref{i min probl} is  specified as
 
  \begin{equation}\label{example}
    \min_{F\in \R, F_{\rm vi}>0}\left(\frac{1}{2}\left|F F_{\rm vi}^{-1}{-}1\right|^2+\frac{1}{2} |F_{\rm vi}{-}1|^2+\frac{1}{2\tau}\left| (F_{\rm vi}{-}F_{\rm vi}^{i-1})(F_{\rm vi}^{i-1})^{-1}\right|^{ 2}\right) \quad \text { for } i=1, \ldots, N,
  \end{equation}
   Take now  initial values $(F^0,F^0_{\rm vi})$ with
  $F^0_{\rm vi}\neq 1$,  so that some nonvanishing viscous stress
  present at time $0$. In this case,  
it is easy to  check  that the constant in time
  solution $(F^0,F^0_{\rm vi})$ satisfies the energy inequality and
  semistability, but it is not approximable.  This implies that
  the viscous strain $F_{\rm vi}$ corresponding to an approximable
  solution must evolve with time, see Figure
  \ref{fig:example1}. In this simple setting, asking the solution
  of the continuous problem to be approximable indeed implies uniqueness, as
  all discrete trajectories converge to the unique solution of the limiting
  differential problem.

\subsection{Assumptions for the linearization theory}\label{sec:hp linearization}

In addition to the assumptions stated in Section \ref{sec:hp existence}, we will require the following conditions in order to prove the linearization result.

On the stored elastic energy density $W_{\rm el}$ we assume that:
\begin{enumerate}[label=(L\arabic*)]
  \item\label{L1} $W_{\rm el}$ is locally Lipschitz;
  \item\label{L2} $W_{\rm el}$ satisfies the growth condition
  \begin{equation}\label{b_el}
    W_{\rm el}(A)\geq c_4 \operatorname{dist}^2 (A,SO(d))
  \end{equation}
  for some $c_4>0$;
  \item\label{L3} there exists a positive definite tensor $\C_{\rm el}$ such that, for every $\delta>0$, there exists $c_{\rm el}(\delta)>0$ satisfying
  \begin{equation}
    \left|W_{\rm el}(I+A)-|A|^2_{\C_{\rm el}}\right|\leq  \delta |A|^2_{\C_{\rm el}} \quad \text{ for every } A\in B_{c_{\rm el}(\delta)}(0).\label{L3a}
  \end{equation}
\end{enumerate}
In particular these conditions imply that $\C_{\rm el}$ is symmetric and
\begin{equation*}
  c_{4} |A^{\rm {sym}}|^2\leq |A|^2_{\C_{\rm el}} \quad \text{ for every } A\in \R^{d\times d}.
\end{equation*}
We can also equivalently state  inequality \eqref{L3a}  as follows:
\begin{equation}\label{c_el}
  (1-\delta)|A|^2_{\C_{\rm el}}\leq W_{\rm el}(I+A)\leq  (1+\delta) |A|^2_{\C_{\rm el}} \quad \text{ for every } A\in B_{c_{\rm el}(\delta)}(0).
\end{equation}

Concerning the viscous stored energy density $W_{\rm vi}$ we ask that
\begin{enumerate}[label=(L\arabic*)]\setcounter{enumi}{3}
  \item\label{L4}
  \begin{equation*}
    W_{\rm vi}(A)= \begin{cases*}
      \widetilde{W}_{\rm vi}(A) & if $A \in K$ \\
      \infty        & otherwise,
    \end{cases*}
\end{equation*}
where $K\subset \subset SL(d)$ contains a neighbourhood of the identity;
  \item\label{L5} $\widetilde{W}_{\rm vi}$ is locally Lipschitz continuous in a neighbourhood of the identity and
  \begin{equation}\label{b_vi}
    \widetilde{W}_{\rm vi}(I+A)\geq c_5|A|^2 \quad \text{ for every } A\in \R^{d\times d} \text{ with } I+A \in K
  \end{equation}
  for some $c_5>0$;
  \item\label{L6} there exists a positive definite tensor $\C_{\rm vi}$ such that, for every $\delta>0$, there exists $c_{\rm vi}(\delta)>0$ satisfying
  \begin{equation*}
    \left|\widetilde{W}_{\rm vi}(I+A)-|A|^2_{\C_{\rm vi}}\right|\leq  \delta |A|^2_{\C_{\rm vi}} \quad \text{ for every } A\in B_{c_{\rm vi}(\delta)}(0),
  \end{equation*}
  \end{enumerate}
  or, equivalently,
  \begin{equation}\label{c_vi}
    (1-\delta)|A|^2_{\C_{\rm vi}}\leq \widetilde{W}_{\rm vi}(I+A)\leq  (1+\delta) |A|^2_{\C_{\rm vi}} \quad \text{ for every } A\in B_{c_{\rm vi}(\delta)}(0).
  \end{equation}
As above we have that
\begin{equation*}
  c_{5} |A^{\rm {sym}}|^2\leq |A|^2_{\C_{\rm vi}} \quad \text{ for every } A\in \R^{d\times d}.
\end{equation*}
Moreover, there exists a constant $c_{K}>0$ (depending only on the compact set $K$) such that
\begin{equation}\label{ck}
    |A|+|A^{-1}|\leq c_K \quad \text{ for every } A\in K
\end{equation}
and
\begin{equation*}
    |A-I|\geq \frac{1}{c_K} \quad \text{ for every } A\in SL(d)\setminus K.
\end{equation*}
These last two inequalities  will provide  $L^{\infty}$-bounds
on the terms $\varepsilon \nabla v$ and $(I+\varepsilon \nabla
v)^{-1}$ later on. Note  however that the effect of the 
constraint $K$ will disappear as $\varepsilon\rightarrow 0$. In particular, the limiting linearized problem is independent of $K$.

On the forcing term $\ell^0$ we assume that
\begin{enumerate}[label=(L\arabic*)]\setcounter{enumi}{6}
  \item\label{L7} $  \ell^0 \in W^{1,1}{\left(0,T;(H^{1}(\Omega;\R^d))^*\right)}$.
\end{enumerate}

Finally, on the dissipation density $\psi$ we assume that
\begin{enumerate}[label=(L\arabic*)]\setcounter{enumi}{7}
  \item\label{L8}
  $\psi$ satisfies the growth condition
  \begin{equation}\label{a_psi}
    \psi(A)\geq c_6|A|^2\quad \text{ for every } A\in \R^{d \times d}
  \end{equation}
  for some $c_6>0$;
  \item\label{L9} there exists a positive definite tensor $\D$ such that, for every $\delta>0$, there exists $c_{\psi}(\delta)>0$ satisfying
  \begin{equation}\label{b_psi}
    \left|\psi(A)-|A|^2_{\D}\right|\leq  \delta |A|^2_{\D} \quad \text{ for every } A\in B_{c_{\psi}(\delta)}(0);
  \end{equation}
  \item\label{L10} $\psi$ is positively $2$-homogeneous, i.e., 
  \begin{equation*}
    \psi(\lambda A)=\lambda^{2}\psi(A) \quad \text{ for every } A\in \R^{d\times d}, \lambda \geq 0.
  \end{equation*}
\end{enumerate}

The specification $p_{\psi}=2$ of assumption \ref{L10} (compare with
the more general $p_{\psi}\geq 2$ from \ref{E6}) is just needed 
in the linearization setting to recover the linearized energy
inequality  \eqref{energy inequality 0} below.



\subsection{Linearization result}\label{sec:res linearization}

 Before moving on,  let us reformulate the setting and the existence results of Proposition \ref{thm: existence} and Theorem \ref{thm: existence limit} in terms of the linearization variables $u$ and $v$.
For all $\varepsilon>0$ fixed, the admissible set $\mathcal{A}$ is equivalently rewritten as
\begin{equation*}
  \begin{aligned}
  {\widetilde{\mathcal{A}}}_{\varepsilon}:=&
         \Bigg\{(u,v)\in W^{1,q}(\Omega;\R^d)\times W^{1,p_{\text{\rm
               vi}}}(\Omega;\R^d) \;\Bigg|\;   u=0  \text{ on }
         \Gamma_{D},\;
          \det (I+ \varepsilon \nabla v)=1, \\  &\qquad \int_{\Omega}v\,{\rm
            d}X=0,\;|\Omega|= |(\operatorname{id}+\varepsilon
          v)(\Omega)|, \; (\operatorname{id}+\varepsilon v)(\Omega)\in \mathcal{J}_{\eta_1,\eta_2}\Bigg\},
  \end{aligned}
\end{equation*}
where we recall that $\Omega$ is chosen to be such that $\int_{\Omega}\operatorname X\,{\rm d} X=0$ so that
\begin{equation*}
  0\stackrel{\eqref{av0}}{=}\int_{\Omega}y_{\rm vi}\,{\rm d} X=\int_{\Omega}(\operatorname{id}+\varepsilon v)\,{\rm d} X=\varepsilon \int_{\Omega}v\,{\rm d} X.
\end{equation*}

We use the following notation for the rescaled energies and dissipation
   \begin{align*}
      \mathcal{W}^{\varepsilon}_{\text{\rm
          el}}(u,v)&:=\frac{1}{\varepsilon^2}\int_{\Omega}W_{\text{\rm
          el}}\left((I{+}\varepsilon \nabla u)(I{+}\varepsilon \nabla
        v)^{-1}\right), \\
    \mathcal{W}^{\varepsilon}_{\text{\rm
     vi}}(v)&:=\frac{1}{\varepsilon^2}\int_{\Omega}W_{\text{\rm
              vi}}(I{+}\varepsilon \nabla v) ,\\
     \Psi^{\varepsilon}(v,\dot{v})&:=\frac{1}{\varepsilon^2}\int_{\Omega}\psi\left( \varepsilon \nabla
                                    \dot{v}(I{+}\varepsilon
                                    \nabla{v})^{-1}
                                    \right).
   \end{align*}
   Their corresponding linearized counterparts read
    \begin{align*}
    \mathcal{W}^0 _{\text{\rm el}}(u,v)&:=  \frac{1}{2}\int_{\Omega}\nabla (u{-}v): \C_{\text{\rm
          el}}\nabla (u{-}v),\\
 \mathcal{W}^0 _{\text{\rm vi}}(v)&:= \frac{1}{2} \int_{\Omega}\nabla v: \C_{\text{\rm vi}}\nabla v
   ,\\
 \Psi^0(\dot{v})&:=  \frac{1}{2}\int_{\Omega}\mathbb{D}\nabla \dot
        {v}:\nabla \dot{v}.
   \end{align*}
   
We also define for brevity
\begin{equation*}
 \E^{\varepsilon}(u,v):=\mathcal{W}^\varepsilon _{\text{\rm vi}}(v)+\mathcal{W}^\varepsilon _{\text{\rm el}}(u,v)-\langle \ell^{0} ,u \rangle \quad \text{ and } \quad \E^{0}(u,v):=\mathcal{W}^0 _{\text{\rm vi}}(v)+\mathcal{W}^0 _{\text{\rm el}}(u,v)-\langle \ell^0 ,u \rangle.
\end{equation*}
Finally, let $(u^{0}_{\varepsilon},v^{0}_{\varepsilon})\in
{\widetilde{\mathcal{A}}}_{\varepsilon} $ be a {\it well-prepared} sequence of initial data, namely
\begin{equation}\label{initial condition epsilon}
  (u^{0}_{\varepsilon},v^{0}_{\varepsilon})\rightharpoonup (u^{0},v^{0}) \text{ in } H^{1}(\Omega)\times H^{1}(\Omega) \quad \text{ and } \quad \lim_{\varepsilon \rightarrow 0} \E^{\varepsilon}(u^{0}_{\varepsilon},v^{0}_{\varepsilon})=\E^{0}(u^{0}
,v^{0})
\end{equation}

Proposition \ref{thm: existence} and Theorem \ref{thm: existence
  limit} can therefore be rewritten in terms of the new variables
$(u,v)$ and in the presence of the rescaling prefactor
$1/\varepsilon^2$ as follows.
\medskip

\begin{corollary}[Existence in terms of $(u_\varepsilon,v_\varepsilon)$]\label{thm: existence epsilon}
  Under the assumptions {\rm \ref{E1}, \ref{E2}, \ref{E3}, \ref{E4},
    \ref{E5}}, and {\rm \ref{L10}} of Section {\rm \ref{sec:hp
      existence}}  and \eqref{initial condition epsilon}  for every $\varepsilon>0$ there exist a sequence of partitions $(\Pi_{\tau^{\varepsilon}})_{\tau^{\varepsilon}}$
  of the interval $[0,T]$ with mesh size
  $\tau^{\varepsilon}\rightarrow 0$ and functions $(u_{\varepsilon},v_{\varepsilon}):[0,T]\rightarrow {\widetilde{\mathcal{A}}}_{\varepsilon} $
  such that for every $t\in [0,T]$
\begin{align}
  &\text{\rm Approximation:}\notag\\[1mm]
  &(\overline{u}_{\tau^\varepsilon}(t),\overline{v}_{\tau^\varepsilon}(t))\rightharpoonup (u_{\varepsilon}(t),v_{\varepsilon}(t))  \text{ in } W^{1,q}(\Omega;\R^d)\times W^{1,p_{\text{\rm vi}}}(\Omega;\R^d)\notag\\[4mm]
 & \text{\rm Energy inequality:}\notag\\
  \quad&\mathcal{W}^{\varepsilon}_{\text{\rm vi}}(v_{\varepsilon}(t))+\mathcal{W}^{\varepsilon}_{\text{\rm el}}(u_{\varepsilon}(t),v_{\varepsilon}(t))-\langle \ell^{0}, u_{\varepsilon}(t)\rangle +
       2 \int_{0}^t \Psi^{\varepsilon}(v_{\varepsilon},\dot{v}_{\varepsilon})\notag \\
      &\quad\leq \mathcal{W}^{\varepsilon}_{\text{\rm vi}}(
      v^{0}_{\varepsilon})+\mathcal{W}^{\varepsilon}_{\text{\rm el}}(u^{0}_{\varepsilon},v^{0}_{\varepsilon})-
      \int_{0}^{t} \langle\dot{\ell}^{0},u_{\varepsilon} \rangle \label{energy inequality epsilon}\\[2mm]
  %
  &\text{\rm Semistability:}\notag\\[2mm]
  \quad&\mathcal{W}^{\varepsilon}_{\text{\rm el}}(u_{\varepsilon}(t),v_{\varepsilon}(t))
    {-}\langle \ell^{\varepsilon}(t), u_{\varepsilon}(t)\rangle
    \leq \mathcal{W}^{\varepsilon}_{\text{\rm el}}(\tilde{u}_{\varepsilon},v_{\varepsilon}(t))
    {-}\langle \ell^{\varepsilon}(t), \tilde{u}_{\varepsilon}(t)\rangle \notag\\
    &\quad\quad \forall \tilde{u}_{\varepsilon} \text{ with } \left(\tilde{u}_{\varepsilon}, v_{\varepsilon}(t)\right) \in {\widetilde{\mathcal{A}}}_{\varepsilon}.\label{semistability epsilon}
  \end{align}
\end{corollary}

In the following result, we show that a sequence
$(u_{\varepsilon},v_{\varepsilon})_{\varepsilon}$ of  
approximable solutions at level $\varepsilon$ converges weakly to
$(u,v)$ satisfying the linearized energy and the linearized
semistability inequalities.
\medskip

\begin{theorem}[Linearization]\label{thm limit lin}
  For every $\varepsilon>0$ let $(u_{\varepsilon},v_{\varepsilon})$ be an approximable solutions given as in Corollary {\rm \ref{thm: existence epsilon}}.
Then, under the assumptions {\rm \ref{L1}, \ref{L2}, \ref{L3},
  \ref{L4}, \ref{L5}, \ref{L6}, \ref{L7}, \ref{L8}, \ref{L9}}, and
{\rm \ref{L10}} of Section {\rm \ref{sec:hp linearization}}  and \eqref{initial condition epsilon}  there exist functions $(u,v):[0,T]\rightarrow H^{1}_{\Gamma_D}(\Omega;\R^d)\times H^{1}_{\sharp}(\Omega;\R^d)$ such that, for every $t\in [0,T]$,  up to a not relabeled subsequence,
\begin{align*}
  u_{\varepsilon}(t)\rightharpoonup u(t), \; v_{\varepsilon}(t)\rightharpoonup v(t) \; &\text{ weakly in } H^{1}(\Omega;\R^{d}), \\
   \nabla \dot  v_{\varepsilon}(t)\rightharpoonup \nabla \dot  v(t) \; &\text{ weakly in } L^{2}(\Omega;\R^{d\times d}).
\end{align*}
Moreover, for every $t\in [0,T]$, we have
\begin{align}
  &\text{\rm Linearized energy inequality:}\notag\\
  \quad&\mathcal{W}^{0}_{\text{\rm vi}}(v(t))+\mathcal{W}^{0}_{\text{\rm el}}(u(t),v(t))-\langle \ell^0(t), u(t)\rangle +
       2\int_{0}^t \Psi^{0}(\dot{v}(s))\notag\\
      &\quad\leq \mathcal{W}^{0}_{\text{\rm vi}}(
      v^{0})+\mathcal{W}^{0}_{\text{\rm el}}(u^{0},v^{0}) -
        \langle \ell^0(0),u^0
        \rangle  -
      \int_{0}^{t} \langle\dot{\ell}^0(s),u(s) \rangle,\label{energy inequality 0}\\[2mm]
  &\text{\rm Linearized semistability:}\notag\\[2mm]
  \quad&\mathcal{W}^{0}_{\text{\rm el}}(u(t),v(t))-\langle \ell^0(t), u(t)\rangle \leq \mathcal{W}^{0}_{\text{\rm el}}(\hat{u},v(t))-\langle \ell^0(t),
    \hat u\rangle \quad\forall \hat{u}\in H^{1}_{\Gamma_D}(\Omega;\R^d).\label{semistability 0}
\end{align}
\end{theorem}

 The proof of Theorem \ref{thm limit lin} is to be found in Section \ref{sec:Proof of Linearization}
below.

Before moving on, let us remark that the linearized 
  energy inequality \eqref{energy inequality 0} and the linearized semistability
  \eqref{semistability 0} cannot be expected to uniquely determine
  solutions of the linearized problem \eqref{linear PDE1}-\eqref{linear PDE2}. On the other
  hand, inequalities \eqref{energy inequality 0}-\eqref{semistability
    0} would uniquely characterize solutions $(u,v)$ to \eqref{linear PDE1}-\eqref{linear PDE2} if in addition one assumes that $(u,v)$ are {\it
    approximable}, namely, are limits of time discretizations of \eqref{linear PDE1}-\eqref{linear PDE2}. Although the trajectories $(u,v)$ are  limits of approximable solutions
  $(u_{\varepsilon},v_{\varepsilon})$, we are not able
  to prove that $(u,v)$ are approximable themselves, for the property
  of being approximable seems not guaranteed to pass to the
  linearization limit.

\section{Time-discretization scheme: Proof of Proposition \ref{thm: existence}}\label{sec:time discr}

 To start with,  notice that the infimum in  the
incremental problems  \eqref{i min probl} is finite for every $i=1,\dots, N_\tau$. Indeed, since the initial condition satisfies $\E(0,y_{\text{\rm el}}^0,y_{\text{\rm vi}}^0)<\infty$, by arguing by induction and
choosing $(y_{\text{\rm el}},y_{\text{\rm vi}})=(y_{\text{\rm el}}^{i-1},y_{\text{\rm vi}}^{i-1})$, we get that
\begin{equation*}
  \mathcal{E}(t_i, y_{\text{\rm el}}, y_{\text{\rm vi}})+\tau \Psi\left(y^{i-1}_{\rm vi},\frac{y_{\rm vi}-y^{i-1}_{\rm vi}}{\tau}   \right)=\mathcal{E}(t_i, y_{\text{\rm el}}^{i-1}, y_{\text{\rm vi}}^{i-1})< \infty.
\end{equation*}
Fix now $1\le i\le N$ and let $(y^i_{\text{\rm el},m},y^i_{\text{\rm vi},m} )_{m\in \N}=(y_{\text{\rm el},m},y_{\text{\rm vi},m} )_{m\in \N}\subset \mathcal{A}$
be an infimizing sequence for problem \eqref{i min probl} at time step $i$.

\subsection{Coercivity}

Let us first show that $(y_{\text{\rm el},m},y_{{\rm vi},m} )_{m\in \N}$ is bounded in $W^{1,p_{\text{\rm el}}}(y_{{\rm vi}, m}(\Omega);\R^d)$
$\times$ $ W^{1,p_{\text{\rm vi}}}(\Omega;\R^d)$. This requires some care since $y_{\text{\rm el},m}$ is defined on the moving domain $y_{\text{\rm vi},m}(\Omega)$. Since the infimum is finite, we have by \eqref{growthel} and \eqref{growthvi}
\begin{equation*}
  c_1\int_{y_{{\rm vi},m}(\Omega)} |\nabla y_{{\rm el},m}|^{p_{\rm el}}+c_2 \int_{\Omega} |\nabla y_{{\rm vi},m}|^{p_{\rm vi}}-\frac{|\Omega|}{c_2}\leq
  \mathcal{W}(y_{\text{\rm el},m},y_{\text{\rm vi},m}) \leq c -\langle\ell (t_i),y_m \rangle
\end{equation*}
where we have posed $y_m:=y_{{\rm el}, m}\circ y_{{\rm vi}, m}$.
The loading term can be controlled as follows
\begin{equation*}
  \begin{aligned}
    |\langle\ell (t_i),y_m \rangle|&\leq \|\ell (t_i)\|_{(W^{1,q}(\Omega))^*}\|y_m\|_{W^{1,q}(\Omega)}\leq c\|\ell (t_i)\|_{(W^{1,q}(\Omega))^*}\|\nabla y_m\|_{L^{q}(\Omega)}\\
    &\stackrel{\text{H\"older}}{\leq} c \|\ell (t_i)\|_{(W^{1,q}(\Omega))^*}\|\nabla y_{\text{\rm el},m}\|_{L^{p_{\text{\rm el}}}(y_{\text{\rm vi},m}(\Omega))}\|\nabla y_{\text{\rm vi},m}\|_{L^{p_{\text{\rm vi}}}(\Omega)}\\
    &\stackrel{\text{Young}}{\leq} c \|\ell (t_i)\|^{1/q'}_{(W^{1,q}(\Omega))^*}+\frac{c_1}{2}\|\nabla y_{\text{\rm el},m}\|^{p_{\text{\rm el}}}_{L^{p_{\text{\rm el}}}(y_{\text{\rm vi},m}(\Omega))}
    +\frac{c_2}{2}\|\nabla y_{\text{\rm vi},m}\|^{p_{\text{\rm vi}}}_{L^{p_{\text{\rm vi}}}(\Omega)}.
  \end{aligned}
\end{equation*}
This entails that
\begin{equation*}
  \|\nabla y_{\text{\rm el},m}\|^{p_{\text{\rm el}}}_{L^{p_{\text{\rm el}}}(y_{\text{\rm vi},m}(\Omega))}
  +\|\nabla y_{\text{\rm vi},m}\|^{p_{\text{\rm vi}}}_{L^{p_{\text{\rm vi}}}(\Omega)}\leq c,
\end{equation*} which in turn guarantees that
\begin{equation*}
  \|\nabla y_{m}\|^{q}_{L^{q}(\Omega)}\leq c.
\end{equation*}

Now,  using the  growth condition \eqref{growthvi} and the
Poincar\'e-Wirtinger inequality, recalling that  $y_{\text{\rm vi}}$
has  zero mean,  we have that
\begin{equation*}
  \|y_{\text{\rm vi},m}\|_{W^{1,p_{\text{\rm vi}}}(\Omega)}\leq c.
\end{equation*}

Recalling that $y_m$ satisfies the Dirichlet  boundary condition \eqref{DirBC}, by the Poincar\'e inequality we obtain $$\| y_m \|_{W^{1,q}(\Omega)}\leq c.$$
A change of variables ensures that
\begin{equation*}
  \int_{y_{\text{\rm vi},m}(\Omega)}|y_{\text{\rm el},m}|^q \,{\rm d}\xi=\int_{\Omega}|y_m|^q \,{\rm d}X\leq c
\end{equation*}
so that $\| y_{\text{\rm el},m} \|_{L^{q}(y_{{\rm vi},m}(\Omega))}\leq c$, as well.
Again the Poincar\'e
 inequality guarantees that
\begin{equation}\label{ineqyel}
  \| y_{\text{\rm el},m} \|_{W^{1,p_{\rm el}}(y_{{\rm
        vi},m}(\Omega) ) }\leq c.
\end{equation}

Up to a not relabeled subsequence we hence have that
\begin{align}
  y_{\text{\rm vi},m}=y^i_{\text{\rm vi},m}\rightharpoonup y^i_{\text{\rm vi}} \quad &\text{ in } W^{1,p_{\text{\rm vi}}}(\Omega;\R^d) \label{lim y_vi m}\\
  y_{m}=y^i_{m}=y^i_{\text{\rm vi},m}\circ y^i_{\text{\rm el},m}\rightharpoonup y^i \quad &\text{ in } W^{1,q}(\Omega;\R^d).  \notag
\end{align}

We now want to extract a converging subsequence   from the elastic
deformations  $y_{\text{\rm el},m}$, which are however defined on the moving domains $y_{\text{\rm vi},m}(\Omega)$. Consider the trivial extensions $\overline{y_{\text{\rm el},m}}$ and $\overline{\nabla y_{\text{\rm el},m}}$ on the whole $\R^d$
by setting ${y_{\text{\rm el},m}}$ and ${\nabla y_{\text{\rm el},m}}$ to be zero on $\R^d\setminus y_{\text{\rm vi},m}(\Omega) $, respectively.
Recalling the bound \eqref{ineqyel}, we have (up to a subsequence)
\begin{align}
  \overline{y_{\text{\rm el},m}}\rightharpoonup  y^i_{\text{\rm el}}\quad& \text{ in } L^{p_{\text{\rm el}}}(\R^d;\R^{d}) \notag\\
 \overline{\nabla y_{\text{\rm el},m}}\rightharpoonup G\quad& \text{ in } L^{p_{\text{\rm el}}}(\R^d;\R^{d\times d}).\label{lim y_el m}
\end{align}
 We want to show that $G=\nabla y^i_{\text{\rm el}}$ on the limiting
 set $y^i_{\text{\rm vi}}(\Omega)$. By Sobolev embedding, possibly by
 extracting a further subsequence, we have that  $y_{\text{\rm
     vi},m} \rightarrow y^i_{\text{\rm vi}} $  uniformly. Letting $\omega \subset \subset y^i_{\text{\rm vi}}(\Omega)$, for $m$ large enough we eventually have that $\omega \subset \subset y_{\text{\rm vi},m}(\Omega)$.
 By uniqueness of the limit we have $\overline{y_{\text{\rm el},m}}\rightharpoonup y^i_{\text{\rm el}} \text{ in } L^{p_{\text{\rm el}}}(\omega;\R^{d})$ and $\overline{\nabla y_{\text{\rm el},m}}\rightharpoonup \nabla y^i_{\text{\rm el}}=G \text{ in } L^{p_{\text{\rm el}}}(\omega;\R^{d\times d})$.
 Hence $G=\nabla y^i_{\text{\rm el}}$ in every $\omega \subset \subset y^i_{\text{\rm vi}}(\Omega)$. An exhaustion argument ensures that $G=\nabla y^i_{\text{\rm el}}$ in $ y^i_{\text{\rm vi}}(\Omega)$.


\subsection{Closure of the set of admissible deformations}\label{Closure adm def}
Let us now check that the weak limit $(y^i_{\text{\rm el}},y^i_{\text{\rm vi}} )$ belongs to the admissible set $\mathcal{A}$.
First of all, since $p_{\text{\rm vi}}>d$ we have that
\begin{equation*}
  1=\det \nabla y^i_{\text{\rm vi},m}\rightharpoonup \det\nabla y^i_{\text{\rm vi}} \quad \text{ in } L^{p_{\rm vi}/d}(\Omega)
\end{equation*}
and hence $\det \nabla y^i_{\text{\rm vi}}=1$ almost everywhere. 
On the other hand,  Lemmas $\rm 3.1-3.2$  in \cite{Melching} 
imply  that $y^i_{\text{\rm vi}}(\Omega)\in \mathcal{J}_{\eta_1,
  \eta_2}$. By the linearity of the  mean  and trace operators and by the weak convergence of $y^i_{{\rm vi},m}$, we find $\int_{\Omega}y^i_{\text{\rm vi}}\,{\rm d}X=0$ and $y^i=\operatorname{id}$ on $\Gamma_{D}$.
Moreover, by \cite[Lemma $\rm 5.2 (i)$]{Grandi} we have that
\begin{equation*}
  |y^i_{{\rm vi},m}(\Omega)\Delta y^i_{{\rm vi}}(\Omega)|\rightarrow 0,
\end{equation*}
where the symbol $\Delta$ denotes the symmetric difference,
and, for every $\omega \subset \Omega$, that $\mathbbm{1}_{y_{\text{\rm vi},m}^i(\omega)}\rightarrow \mathbbm{1}_{y_{\text{\rm vi}}^i(\omega)}$ almost everywhere in $\Omega$.
This implies that $y^i_{{\rm vi}}$ satisfies the Ciarlet-Ne\v cas condition, since
\begin{equation*}
  |\Omega|=|y^i_{{\rm vi},m}(\Omega)|\rightarrow |y^i_{{\rm vi}}(\Omega)|.
\end{equation*}

It remains to show that $y^i=y^i_{\text{\rm el}}\circ y_{\text{\rm
    vi}}^i$. Let us take  any measurable  $\omega \subset \Omega$ and consider, by changing variables,
\begin{equation*}
  \begin{aligned}
    \int_{\omega}y^i(X)\,{\rm d}X&\leftarrow\int_{\omega}y_m(X)\,{\rm d}X=\int_{\omega}y_{\text{\rm el},m} (y_{\text{\rm vi},m}(X))\,{\rm d}X=\int_{y_{\text{\rm vi},m}(\omega)}y_{\text{\rm el},m}(\xi)\,{\rm d}\xi\\
    &=\int_{\R^d}y_{\text{\rm el},m}(\xi)\mathbbm{1}_{y_{\text{\rm vi},m}(\omega)}(\xi)\,{\rm d}\xi\rightarrow \int_{\R^d}y_{\text{\rm el}}^i(\xi)\mathbbm{1}_{y_{\text{\rm vi}}^i(\omega)}(\xi)\,{\rm d}\xi
    =\int_{\omega}y_{\text{\rm el}}^i( y_{\text{\rm vi}}^i(X))\,{\rm d}X,
  \end{aligned}
\end{equation*}
where in the last limit we used the weak convergence of $ y_{\text{\rm el},m}$ and the strong convergence of $\mathbbm{1}_{y_{\text{\rm vi}}^i(\omega)}$.
Since $\omega \subset \Omega$ is arbitrary we conclude  that $y^i=y^i_{\text{\rm el}}\circ y_{\text{\rm
    vi}}^i$. In particular, we have that  $(y_{\rm
  el}^i,y_{\rm vi}^i)\in \mathcal{A}$. 

\subsection{Weak lower semicontinuity}

We aim to show that the functional in \eqref{i min probl} is weakly
lower semicontinuous with respect to the above convergences. 

By polyconvexity of the viscous energy density $W_{\text{\rm vi}}$ and \eqref{lim y_vi m}, we have
\begin{equation*}
\int_{\Omega} W_{\text{\rm vi}}(\nabla y^{i}_{\text{\rm vi}})\,{\rm d}X \leq \liminf_{m\rightarrow \infty}\int_{\Omega} W_{\text{\rm vi}}(\nabla y_{\text{\rm vi},m})\,{\rm d}X.
\end{equation*}
For what concerns the dissipation, from the weak convergence of $y_{\text{\rm vi},m} $ in $W^{1,p_{\text{\rm vi}}}(\Omega)$, we also have
\begin{equation*}
  \frac{\nabla (y_{\text{\rm vi},m}- y^{i-1}_{\text{\rm vi}})}{\tau}(\nabla y^{i-1}_{\text{\rm vi}})^{-1}\rightharpoonup   \frac{\nabla (y^i_{\text{\rm vi}}- y^{i-1}_{\text{\rm vi}})}{\tau}(\nabla y^{i-1}_{\text{\rm vi}})^{-1} \quad \text{ in }  L^{p_\psi}(\Omega;\R^{d\times d}).
\end{equation*}
Hence, by the weak lower semicontinuity of $\Psi$, it follows that
\begin{equation*}
   \Psi\left(y^{i-1}_{\rm vi},\frac{y^i_{\rm vi}-y^{i-1}_{\rm vi}}{\tau}   \right)\leq \liminf_{m\rightarrow \infty} \Psi\left(y^{i-1}_{{\rm vi}},\frac{y_{{\rm vi},m}-y^{i-1}_{\rm vi}}{\tau}   \right).
\end{equation*}

As the loading term is linear, we have
\begin{equation*}
  \langle \ell (t_i),y^i \rangle = \lim_{m\rightarrow \infty}\langle \ell (t_i),y_m \rangle
\end{equation*}
by weak convergence of $y_m$.

Finally, for any $\omega \subset \subset y^i_{\text{\rm vi}}(\Omega)$ we can treat the elastic energy as follows
 \begin{equation*}
   \int_{\omega} W_{\text{\rm el}}(\nabla y^i_{\text{\rm el}})\,{\rm d}\xi\leq \liminf_{m\rightarrow \infty}\int_{\omega} W_{\text{\rm el}}(\nabla y_{\text{\rm el},m})\,{\rm d}\xi \stackrel{\eqref{growthel}}{\leq} \liminf_{m\rightarrow \infty}\int_{y_{\text{\rm vi},m}(\Omega)} W_{\text{\rm el}}(\nabla y_{\text{\rm el},m})\,{\rm d}\xi
   ,
 \end{equation*}
where we have used the polyconvexity of $W_{\text{\rm el}}$ and convergence \eqref{lim y_el m}. Taking the supremum over $\omega \subset \subset y^i_{\text{\rm vi}}(\Omega)$ we conclude via an exhaustion argument that
\begin{equation*}
  \int_{y_{\text{\rm vi},m}(\Omega)} W_{\text{\rm el}}(\nabla y_{\text{\rm el}})\,{\rm d}\xi\leq \liminf_{m\rightarrow \infty}\int_{y_{\text{\rm vi},m}(\Omega)} W_{\text{\rm el}}(\nabla y_{\text{\rm el},m})\,{\rm d}\xi.
\end{equation*}

All in all, we have proved that $(y_{\rm
  el}^i,y_{\rm vi}^i)\in \mathcal{A}$ and
\begin{equation*}
  \mathcal{E}(t_i,y_{\text{\rm el}}^i,y_{\text{\rm vi}}^i)
  + \tau \Psi\left(y^{i-1}_{\rm vi},\frac{y_{\rm vi}^i-y^{i-1}_{\rm vi}}{\tau}   \right)=\min_{(y_{\text{\rm el}},y_{\text{\rm vi}})\in \mathcal{A}} \left\{ \mathcal{E}(t_i,y_{\text{\rm el}},y_{\text{\rm vi}})
  + \tau \Psi\left(y^{i-1}_{\rm vi},\frac{y_{\rm vi}-y^{i-1}_{\rm vi}}{\tau}   \right) \right\}
\end{equation*}
so that the assertion of Proposition \ref{thm: existence} follows.

\section{Existence of approximable solutions: Proof of Theorem
  \ref{thm: existence limit}}\label{sec:Proof of Existence}

 We split the proof in subsequent steps. The basic energy
estimate and its consequences are presented in Subsection \ref{sec:
  estimates}. The energy estimate is then sharpened in Subsection
\ref{sec: de giorgi estimates}, leading to the discrete energy
inequality. By taking limits as the time step $\tau$ goes to $0$, the
time-continuous energy inequality \eqref{energy inequality} and the time-continuous
semistability \eqref{semistability} are proved in Subsections \ref{sec:1}
and \ref{sec:2}, respectively. 

\subsection{ Energy estimate and its consequences}\label{sec: estimates}

Let $(y^i_{\rm el}  y^i_{\rm vi})_{i=0}^{N}$ be a solution to \eqref{i min probl}.
By minimality we have, for every $i=1,...,N$,
\begin{equation*}
  \begin{aligned}
    \mathcal{E}(t_i,y^i_{\text{\rm el}},y^i_{\text{\rm vi}})+\tau\Psi\left(y^{i-1}_{\rm vi},\frac{y^i_{\rm vi}-y^{i-1}_{\rm vi}}{\tau}   \right)\leq \mathcal{E}(t_i,y^{i-1}_{\text{\rm el}},y^{i-1}_{\text{\rm vi}})\\
    \quad=
    \mathcal{E}(t_{i-1},y^{i-1}_{\text{\rm el}},y^{i-1}_{\text{\rm vi}})-\int_{t_{i-1}}^{t_i}\langle\dot{\ell},y^{i-1} \rangle .
  \end{aligned}
\end{equation*}
Summing up over $i=1,...,n\leq N$ we get
\begin{equation}\label{discrete energy ineq}
  \mathcal{E}(t_n,y^n_{\text{\rm el}},y^n_{\text{\rm vi}})+\sum_{i=1}^n \tau\Psi\left(y^{i-1}_{\rm vi},\frac{y^i_{\rm vi}-y^{i-1}_{\rm vi}}{\tau}   \right)
  \leq \mathcal{E}(0,y_{\text{\rm el},0},y_{\text{\rm vi},0})-\sum_{i=1}^n \int_{t_{i-1}}^{t_i}\langle\dot{\ell},y^{i-1} \rangle .
\end{equation}
By using the notation for the interpolants, we have, for all $t\in [0,T],$
  \begin{align*}
    &\mathcal{E}(\overline{t}_{\tau}(t),\overline{y}_{\text{\rm el},\tau}(t),\overline{y}_{\text{\rm vi},\tau}(t))+
    \int_{0}^{\overline{t}_{\tau}(t)} \Psi\left(\underline{y}_{{\rm vi},\tau},\dot{\hat{y}}_{{\rm vi},\tau}   \right)\leq \mathcal{E}(0,y_{{\rm el},0},y_{{\rm vi},0})-
    \int_{0}^{\overline{t}_{\tau}(t)} \langle\dot{\ell},\underline{y}_{\tau} \rangle
    \\ &\quad\leq \mathcal{E}(0,y_{{\rm el},0},y_{{\rm vi},0})+
    \int_{0}^{\overline{t}_{\tau}(t)}\|\dot{\ell}\|_{(W^{1,q}(\Omega))^*}\|\underline{y}_{\tau}\|_{W^{1,q}(\Omega)}\\
    &\quad\stackrel{\text{Poincar\'e}}{\leq} \mathcal{E}(0,y_{{\rm el},0},y_{{\rm vi},0})+
    c\int_{0}^{\overline{t}_{\tau}(t)}\|\dot{\ell}\|_{(W^{1,q}(\Omega))^*}\|\nabla \underline{y}_{\tau}\|_{L^{q}(\Omega)}\\
    &\quad\leq \mathcal{E}(0,y_{{\rm el},0},y_{{\rm vi},0})+
    c\int_{0}^{\overline{t}_{\tau}}\|\nabla \underline{y}_{\text{\rm el},\tau}\|_{L^{p_{\text{\rm el}}}(\underline{y}_{{\rm vi},\tau}(t,\Omega))}\|\nabla \underline{y}_{\text{\rm vi},\tau}\|_{L^{p_{\text{\rm vi}}}(\Omega)}.
  \end{align*}
On the other hand, by the growth assumptions \eqref{growthel}, \eqref{growthvi}, and \eqref{dissipation density}, we also have
  \begin{align}
    \mathcal{E}(\overline{t}_{\tau}(t),&\overline{y}_{\text{\rm el},\tau}(t),\overline{y}_{\text{\rm vi},\tau}(t))+
    \int_{0}^{\overline{t}_{\tau}(t)} \Psi\left(\underline{y}_{{\rm vi},\tau},\dot{\hat{y}}_{{\rm vi},\tau}   \right)
    \geq c\|\nabla \overline{y}_{\text{\rm vi},\tau}(t)\|_{L^{p_{\text{\rm vi}}}(\Omega)}^{p_{\rm vi}} \notag\\
     &+c\|\nabla \overline{y}_{\text{\rm el},\tau}(t)\|_{L^{p_{\text{\rm el}}}(\overline{y}_{{\rm vi},\tau}(t,\Omega))}^{p_{\rm el}}
    +
    c\int_0 ^{\overline{t}_{\tau}(t)}\|\nabla \dot
       {\hat{y}}_{\text{\rm vi},\tau}(\nabla \underline{y}_{\text{\rm
       vi},\tau})^{-1}\|^{p_{\psi}}_{L^{p_{\psi}}(\Omega)}- c . \label{ineq 3}
  \end{align}
In particular, combining the two inequalities above we get
\begin{equation*}
  \begin{aligned}
    &c\|\nabla \overline{y}_{\text{\rm vi},\tau}(t)\|^{p_{\rm vi}}_{L^{p_{\text{\rm vi}}}(\Omega)}+ c\|\nabla \overline{y}_{\text{\rm el},\tau}(t)\|^{p_{\rm el}}_{L^{p_{\text{\rm el}}}(\overline{y}_{{\rm vi},\tau}(t,\Omega))}
    +\int_0 ^{\overline{t}_{\tau}(t)}\|\nabla \dot {\hat{y}}_{\text{\rm vi},\tau}(\nabla \underline{y}_{\text{\rm vi},\tau})^{-1}\|^{p_{\psi}}_{L^{p_{\psi}}(\Omega)}
    \\ &\quad\leq
    \mathcal{E}(0,y_{{\rm el},0},y_{{\rm vi},0})+
    c\int_{0}^{\overline{t}_{\tau}(t)}\|\nabla \underline{y}_{\text{\rm el},\tau}\|_{L^{p_{\text{\rm el}}}(\underline{y}_{\text{\rm vi},\tau}(t,\Omega))}\|\nabla \underline{y}_{\text{\rm vi},\tau}\|_{L^{p_{\text{\rm vi}}}(\Omega)}+c
    \\&\quad \leq c+
    c\int_{0}^{\overline{t}_{\tau}(t)}
    \left(\|\nabla \underline{y}_{\text{\rm el},\tau}\|^{p_{\text{\rm el}}}_{L^{p_{\text{\rm el}}}(\underline{y}_{{\rm vi},\tau}(\Omega))}+\|\nabla \underline{y}_{\text{\rm vi},\tau}\|^{p_{\text{\rm vi}}}_{L^{p_{\text{\rm vi}}}(\Omega)}\right).
  \end{aligned}
\end{equation*}
We can apply the Discrete Gronwall Lemma \cite[(C.2.6), p. 534]{Kruzik Roubicek} to find
\begin{equation*}
  \|\nabla \overline{y}_{\text{\rm vi},\tau}(t)\|_{L^{p_{\text{\rm vi}}}(\Omega)}+\|\nabla \overline{y}_{\text{\rm el},\tau}(t)\|_{L^{p_{\text{\rm el}}}(\overline{y}_{{\rm vi},\tau}(\Omega))}\leq c.
\end{equation*}
Thus, for every $t\in [0,T]$ we have 
\begin{equation*}
  \| \overline{y}_{\text{\rm vi},\tau}(t)\|_{
  W^{1,p_{\text{\rm vi}}}(\Omega)}\leq c \quad \text{ and } \quad \| \nabla \overline{y}_{\text{\rm el},\tau}(t)\|_{
  L^{p_{\text{\rm el}}}(\overline{y}_{\text{\rm vi},\tau}(t,\Omega))}\leq c.
\end{equation*}
Then, using Poincar\'e inequality on the total deformation, we find $\| \overline{y}_{\tau}(t)\|_{
W^{1,q}(\Omega)}\leq c $ and hence, as before, for every $t\in [0,T]$
\begin{equation*}
  \| \overline{y}_{\text{\rm el},\tau}(t)\|_{
  W^{1,p_{\text{\rm el}}}(\overline{y}_{\text{\rm vi},\tau}(t,\Omega))}\leq c.
\end{equation*}
Moreover, thanks to \eqref{ineq 3}, we also have for every $t\in[0,T]$
\begin{equation}\label{coercivity diss p_psi}
  \int_0 ^{\overline{t}_{\tau}(t)}\|\nabla \dot {\hat{y}}_{\text{\rm vi},\tau}(\nabla \underline{y}_{\text{\rm vi},\tau})^{-1}\|^{p_{\psi}}_{L^{p_{\psi}}(\Omega)}\leq c.
\end{equation}
By recalling that $1/p_r=1/p_{\psi}+1/p_{\text{\rm vi}}$  this
implies that 
\begin{equation*}
  \int_0^T\|\nabla \dot {\hat{y}}_{\text{\rm vi},\tau}\|_{
  L^{p_r}(\Omega)}\leq
  \int_0^T \|\nabla \dot {\hat{y}}_{\text{\rm vi},\tau}(\nabla \underline{y}_{\text{\rm vi},\tau})^{-1}\|_{
   L^{p_{\psi}}(\Omega)}
\|\nabla \underline{y}_{\text{\rm vi},\tau}\|_{
L^{p_{\text{\rm vi}}}(\Omega)}\leq c.
\end{equation*}


\subsection{ Energy inequality, sharp version}\label{sec: de giorgi estimates}

In the previous section, we have found  the energy estimate 
\eqref{discrete energy ineq},  which features  the dissipation
with a prefactor $1$. In order to prove the  sharp version of the
energy inequality  \eqref{energy inequality} with the prefactor $p_\psi$, we need a finer argument, mutated from \cite{AGS}.

First, we introduce some notation. Let 
$$\mathcal{V} = \Big\{y_{\rm vi}\in W^{1,p_{\text{\rm vi}}}(\Omega;\R^d) \;\Big|\;
      \det \nabla y_{\text{\rm vi}}=1 \text{ a.e. in } \Omega\Big\}$$
and, for all $i=1,\dots,N$, define the functionals $\Phi^i:
[0,T]\times   \mathcal{V}\times
     \mathcal{A}   \rightarrow \R$ as 
\begin{align*} 
  &\Phi^i(\tau;y_{\rm old}  ,y_{\rm el}, y_{\rm vi}):=\mathcal{E}(t_i,y_{\text{\rm el}},y_{\text{\rm vi}})
  + \tau \Psi\left(y_{\rm old},\frac{y_{\rm vi}-y_{\rm old}}{\tau}   \right).
\end{align*}
Recall that, by definition \eqref{dissipation} of $\Psi$ and by the $p_{\psi}$-homogeneity \eqref{p-homogeneous}, we have 
\begin{equation}\label{phi equiv hom}
  \tau \Psi\left(y_{\rm old},\frac{y_{\rm vi}-y_{\rm old}}{\tau}\right)=\frac{1}{\tau^{p_{\psi}-1}}\int_{\Omega}\psi\left((\nabla y_{\rm old})^{-1}(\nabla y_{\rm vi}-\nabla y_{\rm old})\right).
\end{equation}
 For all $(t, y_{\rm old})\in  [0,T]\times  \mathcal{V}$ we  also define the  minimal value of
the latter functional as 
\begin{equation*}
  \phi_{\tau}^i(y_{\rm old}) :=\inf_{(y_{\text{\rm el}},y_{\text{\rm
        vi}})\in \mathcal{A}} \Phi^i(\tau;y_{\rm old}  ,y_{\rm el},
  y_{\rm vi}) 
\end{equation*}
 and denote  the set of minimizers by $J_{\tau}^i(y_{\rm old}):= \argmin\{ \Phi^i(\tau;y_{\rm old}  ,y_{\rm el}, y_{\rm vi}) \,|\, (y_{\text{\rm el}},y_{\text{\rm vi}}) \in \mathcal A \} $, which is nonempty by Proposition \ref{thm: existence}. Finally, introduce 
\begin{align*}
  \Psi_{\tau}^{+,i}(y_{\rm old}):= \sup_{(y_{\text{\rm el},\tau},y_{\text{\rm vi},\tau})\in J_{\tau}^i(y_{\rm old})}\int_{\Omega}\psi\left((\nabla y_{\rm old})^{-1}\frac{(\nabla y_{{\rm vi},\tau}-\nabla{y_{\rm old}})}{\tau}  \right),
  \\
  \Psi_{\tau}^{- ,i}(y_{\rm old}):= \inf_{(y_{\text{\rm el},\tau},y_{\text{\rm vi},\tau})\in J^i_{\tau}(y_{\rm old})}\int_{\Omega}\psi\left((\nabla y_{\rm old})^{-1}\frac{(\nabla y_{{\rm vi},\tau}-\nabla{y_{\rm old}})}{\tau}  \right).
\end{align*}

We start by stating an auxiliary result, providing the continuity
property of the map $\tau\mapsto \phi^i_{\tau}(y_{\rm old})$ in $0$
and the monotonicity of $\tau\mapsto\Psi\left(y_{\rm old}, y_{{\rm
      vi},\tau}{-}y_{\rm old}\right)$.
\medskip

\begin{lemma}
  For every $i=1,\dots, N$ and every $y_{\rm old}\in 
  \mathcal{V}$,  we have
  \begin{equation}\label{lim phi tau to zero}
    \lim_{\tau\searrow 0}\phi^i_{\tau}(y_{\rm old})=\E( t_{i},y_{{\rm el}}, y_{\rm old})
  \end{equation}
  where $y_{{\rm el}}\in\argmin \{\E(t_i,\widetilde{y}_{\rm el}, y_{\rm old})\,|\, \widetilde{y}_{\rm el}\in W^{1,p_{\rm el}}(y_{\rm old}(\Omega);\R^d)\}$.

  Moreover, if  $0<\tau_0<\tau_1$, then 
  \begin{equation}\label{monotonicity diss}
    \Psi\left(y_{\rm old}, y_{{\rm vi},\tau_0}{-}y_{\rm old}\right)\leq \Psi\left(y_{\rm old}, y_{{\rm vi},\tau_1}{-}y_{\rm old}\right) \quad\text{ for every } (y_{{\rm el},\tau_j},y_{{\rm vi},\tau_j})\in J_{\tau_j}^{i}(y_{\rm old}),\; j=0,1.
  \end{equation}
\end{lemma}

\begin{proof}
  We start by proving the continuity property of $\tau\mapsto
  \phi^i_{\tau}(y_{\rm old})$. Let $(y_{{\rm el},\tau },y_{{\rm
      vi},\tau })\in J^{i}_{\tau}(y_{\rm old})$. By the growth
  condition \eqref{dissipation density},  the  $p_{\psi}$-homogeneity \eqref{p-homogeneous}, and coercivity, we have 
  \begin{equation*}
   \|(\nabla y_{\rm old})^{-1}(\nabla y_{{\rm vi},\tau}-\nabla y_{\rm old})\|_{L^{p_{\rm vi}}(\Omega)}\leq c\Psi\left(y_{\rm old},y_{{\rm vi},\tau}{-}y_{\rm old} \right)= c\tau^{p_{\psi}}\Psi\left(y_{\rm old},\frac{y_{{\rm vi},\tau}{-}y_{\rm old}}{\tau}\right)\leq c \tau^{p_{\psi}{-}1}.
  \end{equation*}
   The proves that  $\nabla y_{{\rm vi},\tau} \to
  \nabla y_{\rm old}$ in $L^{p_{\rm vi}}(\Omega; \R^{d\times
    d})$  as $\tau \to 0$.  Moreover, by \eqref{ineqyel}, we have $ y_{{\rm el},\tau} \rightharpoonup  y_{{\rm el}}$ weakly in $W^{1,p_{\rm el}}(\Omega; \R^{d})$. Thus, by weak lower semicontinuity, we have
  \begin{equation*}
    \lim_{\tau \searrow 0}\phi^{i}_{\tau}(y_{\rm old})=\lim_{\tau \searrow 0} \Phi^i(\tau;y_{\rm old},y_{{\rm el},\tau}
    ,y_{{\rm vi},\tau})\geq \liminf_{\tau \searrow 0}\E(t_i,y_{{\rm el},\tau},y_{{\rm vi},\tau})\geq \E(t_i,y_{{\rm el}},y_{{\rm old}}).
  \end{equation*}
  On the other hand,  from   minimality  we get 
  $\E(t_i,y_{{\rm el}},y_{{\rm old}})\geq \phi^i_{\tau}(y_{{\rm old}})$.
  This implies that 
  \begin{equation*}
    \lim_{\tau \searrow 0}\phi^{i}_{\tau}(y_{\rm old})=\E(t_i,y_{{\rm el}},y_{{\rm old}}).
  \end{equation*}
  The fact that $y_{{\rm el}}\in\argmin \{\E(t_i,\widetilde{y}_{\rm el}, y_{\rm old})\,|\, \widetilde{y}_{\rm el}\in W^{1,p_{\rm el}}(y_{\rm old}(\Omega);\R^d)\}$ follows  from   minimality since
  \begin{equation*}
    \E(t_i,y_{\rm el},y_{\rm old})=\lim_{\tau \searrow 0}\phi^{i}_{\tau}(y_{\rm old})\leq \lim_{\tau \searrow 0} \Phi^i(\tau;y_{\rm old}, \widetilde{y}_{\rm el}, y_{\rm old})=\E(t_i,\widetilde{y}_{\rm el},y_{\rm old})
  \end{equation*}
  for every $\widetilde{y}_{\rm el}\in W^{1,p_{\rm el}}(y_{\rm old}(\Omega);\R^d)$.

  Let us now prove the monotonicity of $\tau\mapsto\Psi\left(y_{\rm old}, y_{{\rm vi},\tau}{-}y_{\rm old}\right)$. 
  Let $0<\tau_0<\tau_1$ and $y_{{\rm el},\tau_j}, y_{{\rm vi},\tau_j}\in J^i_{\tau_j}(y_{\rm old})$, $j=0,1$.  From   minimality, we have that
  \begin{align*}
    \phi^i_{\tau_0}&= \E(t_i, y_{{\rm el},\tau_0}, y_{{\rm vi},\tau_0})+\frac{1}{\tau_0^{p_{\psi}-1}} \Psi\left( y_{\rm old}, y_{{\rm el},\tau_0}-y_{\rm old}\right)\\
    &\leq \E(t_i, y_{{\rm el},\tau_1}, y_{{\rm vi},\tau_1})+\frac{1}{\tau_0^{p_{\psi}-1}}\Psi\left( y_{\rm old}, y_{{\rm el},\tau_1}-y_{\rm old}\right)\\
    &= \E(t_i, y_{{\rm el},\tau_1}, y_{{\rm vi},\tau_1})+\frac{1}{\tau_1^{p_{\psi}-1}} \Psi\left( y_{\rm old}, y_{{\rm el},\tau_1}-y_{\rm old}\right)+\left(\frac{1}{\tau_0^{p_{\psi}-1}}- \frac{1}{\tau_1^{p_{\psi}-1}} \right)\Psi\left( y_{\rm old}, y_{{\rm el},\tau_1}-y_{\rm old}\right)\\
    &\leq  \E(t_i, y_{{\rm el},\tau_0}, y_{{\rm vi},\tau_0})+\frac{1}{\tau_1^{p_{\psi}-1}} \Psi\left( y_{\rm old}, y_{{\rm el},\tau_0}-y_{\rm old}\right)+\left(\frac{1}{\tau_0^{p_{\psi}-1}}- \frac{1}{\tau_1^{p_{\psi}-1}} \right)\Psi\left( y_{\rm old}, y_{{\rm el},\tau_1}-y_{\rm old}\right).
  \end{align*}
  This implies that 
  \begin{equation*}
    \left(\frac{1}{\tau_0^{p_{\psi}-1}}- \frac{1}{\tau_1^{p_{\psi}-1}} \right)\Psi\left( y_{\rm old}, y_{{\rm el},\tau_1}-y_{\rm old}\right)\leq \left(\frac{1}{\tau_0^{p_{\psi}-1}}- \frac{1}{\tau_1^{p_{\psi}-1}} \right)\Psi\left( y_{\rm old}, y_{{\rm el},\tau_1}-y_{\rm old}\right),
  \end{equation*}
which concludes the proof.
\end{proof}

In the following Lemma, we calculate the derivative with respect to
$\tau$ of the minimal incremental energy $\phi^i_{\tau}$ and provide a
crucial estimate. \medskip

\begin{lemma}
  For every $y_{\rm old}\in  \mathcal{V}$  and $i=1,\dots N$, the map $\tau \mapsto \phi^i_{\tau}(y_{\rm old})$ is locally Lipschitz on $(0,1)$. Moreover, we have
  \begin{equation}\label{phi tau derivative}
    \frac{d}{d\tau}\phi^i_{\tau}(y_{\rm old})=-(p_{\psi}-1)\Psi_{\tau}^{\pm,i}(y_{\rm old})
  \end{equation}
  for almost every $\tau\in (0,1).$ In particular, for almost every $\tau\in (0,1)$ we have
  \begin{equation}\label{improved estimate}
    \tau\Psi\left(y_{\rm old}, \frac{y_{{\rm vi},\tau}-y_{\rm old}}{\tau}\right)+(p_{\psi}-1)\int_0^\tau \Psi_{r}^{\pm,i}(y_{\rm old}) \, {\rm d}r=\mathcal{E}(t_{i},{y}_{\rm el},y_{\rm old})-\mathcal{E}(t_i,y_{\text{\rm el},\tau},y_{\text{\rm vi},\tau})
  \end{equation}
  for every $(y_{\text{\rm el},\tau},y_{\text{\rm vi},\tau})\in J^i_{\tau}(y_{\rm old})$, for some $y_{{\rm el}}=\argmin \{\E(t_i,\widetilde{y}_{\rm el}, y_{\rm old})\,|\, \widetilde{y}_{\rm el}\in W^{1,p_{\rm el}}(y_{\rm old}(\Omega);\R^d)\}$.

\end{lemma}

\begin{proof}
  For every $\tau_0\neq \tau_1$ and $(y_{\text{\rm el},{\tau_j}},y_{\text{\rm vi},{\tau_j}})\in J^i_{\tau_j}(y_{\rm old})$, $j=0,1$, by minimality we have
  \begin{align*}
    \phi_{\tau_0}(y_{\rm old})-\phi_{\tau_1}(y_{\rm old})&\leq \Phi^i(\tau_0;y_{\rm old},y_{\text{\rm el},{\tau_1}},y_{\text{\rm vi},{\tau_1}})-\Phi^i(\tau_1;y_{\rm old},y_{\text{\rm el},{\tau_1}},y_{\text{\rm vi},{\tau_1}})\\
    &=\frac{1}{\tau_0^{p_{\psi}-1}}\Psi\left(y_{\rm old},y_{{\rm vi},{\tau_1}}-y_{\rm old}\right)-\frac{1}{\tau_1^{p_{\psi}-1}}\Psi\left(y_{\rm old},y_{{\rm vi},{\tau_1}}-y_{\rm old}\right)
    \\
    &=\frac{\tau_1^{{p_{\psi}-1}}-\tau_0^{{p_{\psi}-1}}}{(\tau_1\tau_0)^{p_{\psi}-1}}
    \int_{\Omega}\psi\left((\nabla y_{\rm old})^{-1}(\nabla y_{{\rm vi},{\tau_1}}-\nabla y_{\rm old})\right),
  \end{align*}
where we used \eqref{phi equiv hom}.
We can perform an analogous calculation for \begin{equation*}
  \phi_{\tau_0}(y_{\rm old})-\phi_{\tau_1}(y_{\rm old})\geq \Phi^i(\tau_0;y_{\rm old},y_{\text{\rm el},\tau_0},y_{\text{\rm vi},\tau_0})-\Phi^i(\tau_1;y_{\rm old},y_{\text{\rm el},\tau_0},y_{\text{\rm vi},\tau_0})
\end{equation*} so that, by combining the two above inequalities, for $\tau_0<\tau_1$ we find 
\begin{align*}
  \frac{\tau_1^{{p_{\psi}-1}}-\tau_0^{{p_{\psi}-1}}}{(\tau_1\tau_0)^{p_{\psi}-1}(\tau_1-\tau_0)}
    \int_{\Omega}\psi\left((\nabla y_{\rm old})^{-1}(\nabla y_{{\rm vi},{\tau_0}}-\nabla y_{\rm old})\right)\leq \frac{\phi_{\tau_0}(y_{\rm old})-\phi_{\tau_1}(y_{\rm old})}{\tau_1-\tau_0} \\
    \quad\leq \frac{\tau_1^{{p_{\psi}-1}}-\tau_0^{{p_{\psi}-1}}}{(\tau_1\tau_0)^{p_{\psi}-1}(\tau_1-\tau_0)}
    \int_{\Omega}\psi\left((\nabla y_{\rm old})^{-1}(\nabla y_{{\rm vi},{\tau_1}}-\nabla y_{\rm old})\right).
\end{align*}
Taking the supremum over $(y_{\text{\rm el},{\tau_0}},y_{\text{\rm vi},{\tau_0}})\in J^i_{\tau_0}(y_{\rm old})$ in the left-hand side and the infimum over $(y_{\text{\rm el},{\tau_1}},y_{\text{\rm vi},{\tau_1}})\in J^i_{\tau_1}(y_{\rm old})$ in the right hand side, we find
\begin{equation*}
  \frac{\tau_0(\tau_1^{{p_{\psi}-1}}-\tau_0^{{p_{\psi}-1}})}{\tau_1^{p_{\psi}-1}(\tau_1-\tau_0)}\Psi_{\tau_0}^{+,i}(y_{\rm old})\leq \frac{\phi_{\tau_0}(y_{\rm old})-\phi_{\tau_1}(y_{\rm old})}{\tau_1-\tau_0}\leq \frac{\tau_1(\tau_1^{{p_{\psi}-1}}-\tau_0^{{p_{\psi}-1}})}{\tau_0^{p_{\psi}-1}(\tau_1-\tau_0)}\Psi_{\tau_1}^{-,i}(y_{\rm old}),
\end{equation*}
which implies that $\tau \mapsto \phi^i_{\tau}(y_{\rm old})$ is locally Lipschitz. Then, passing to the limit for $\tau_1{\searrow} \tau$ and $\tau_0{\nearrow }\tau$, we get \eqref{phi tau derivative}.

Integrating \eqref{phi tau derivative} from $\tau_0>0$ to $\tau$, we have
\begin{equation*}
  \phi^i_{\tau}(y_{\rm old})-\phi^i_{\tau_0}(y_{\rm old})=-(p_{\psi}-1)\int_{\tau_0}^{\tau}\Psi^{\pm,i}_r(y_{\rm old})\, {\rm d}r.
\end{equation*}
Letting $\tau_0 {\searrow} 0$, recalling \eqref{lim phi tau to zero}, and the definition of $(y_{\text{\rm el},\tau},y_{\text{\rm vi},{\tau}})\in J^i_{\tau}(y_{\rm old})$, we get \eqref{improved estimate}.
\end{proof}

We now state the definition of {\it De Giorgi variational
  interpolation} \cite[Definition 3.2.1]{AGS}, which in our setting
refers  to  the viscous deformation $y_{\rm vi}$  only. 
\medskip

\begin{definition}[De Giorgi variational interpolation]
  Let $(y_{{\rm el},\tau}^i,y_{{\rm vi},\tau}^i)_{i=0}^{N}$ be an {\it
    incremental solution} of the problem of \eqref{i min probl}. We
  call \emph{De Giorgi variational interpolation} of $(y_{{\rm
      vi},\tau}^i)_{i=0}^{N}$ any interpolation $\widetilde{y}_{{\rm
      vi},\tau}$ of the discrete values  with 
  $(\overline{y}_{{\rm el},\tau},\widetilde{y}_{{\rm
      vi},\tau}):[0,T]\rightarrow \mathcal{A}$  that satisfies 
  \begin{equation*}
    \widetilde{y}_{{\rm vi},\tau}(t)=\widetilde{y}_{{\rm vi},\tau}(
    t_{i-1}+ r  )\in J_{ r }(\widetilde{y}_{{\rm vi},\tau}^{i-1}) \quad \text{ if } t_{i-1}+ r  \in (t_{i-1},t_{i}]
  \end{equation*}
  for every $i=1,\dots, N$.
\end{definition}

The following Proposition provides the sharp energy estimate on the discrete level, providing an equality instead of an inequality.

\begin{proposition}[Discrete energy equality]
  Let $(y_{{\rm el},\tau}^i,y_{{\rm vi},\tau}^i)_{i=0}^{N}$ be an {\it incremental solution} of the problem of \eqref{i min probl}. Then, for every $1\leq n\leq N$ we have 
  \begin{align}
    \tau \sum_{i=1}^n \psi\left(\left(\nabla y_{{\rm vi}, \tau}^{i-1}\right)^{-1}\left(\nabla y_{{\rm vi},\tau}^i-\nabla y_{{\rm vi},\tau}^{i-1}\right)\right)+(p_{\psi}-1)\sum_{i=1}^n \int_{t_{i-1}}^{t_{i}} G_\tau^{p_{\psi}}( r )\, {\rm d}  r +\E\left(t_{i},y_{{\rm el},\tau}^n, y_{{\rm vi},\tau}^n\right)\notag
    \\
    =\E\left(0,y_{{\rm el},0}, y_{{\rm vi},0}\right)-\sum_{i=1}^{n}\int_{t_{i-1}}^{t_{i}}\langle \dot{\ell}, y^{i-1}_{{\rm el},\tau}\circ y^{i-1}_{{\rm el},\tau}\rangle \label{energy equality deiscrete}
  \end{align}
  where 
  \begin{equation*}
    G_\tau(t):=\left({\Psi^{\pm,i}_{ r }(y_{{\rm vi},\tau}^{i-1})}\right)^{1/p_{\psi}} \quad \text { for } t= t_{i-1} + r  \in (t_{i-1},t_{i}].
  \end{equation*}
\end{proposition}

\begin{proof}
  By \eqref{improved estimate} for $y_{\rm old}= y_{{\rm vi},\tau}^{i-1}$, $y_{{\rm vi},\tau}= y_{{\rm vi},\tau}^i$, and $y_{{\rm el},\tau}= y_{{\rm el},\tau}^i$, we find 
  \begin{align*}
    \tau\Psi\left(y_{{\rm vi},\tau}^{i-1}, \frac{y_{{\rm vi},\tau}^i-y_{{\rm vi},\tau}^{i-1}}{\tau}\right)+(p_{\psi}-1)\int_0^\tau |G_{\tau}( r )|^{p_{\psi}} \, {\rm d} r + \mathcal{E}(t_i,y_{\text{\rm el},\tau}^i,y_{\text{\rm vi},\tau}^i) =\mathcal{E}(t_{i},y_{\text{\rm el},\tau}^{i-1},y_{{\rm vi},\tau}^{i-1})\\
    =\mathcal{E}(t_{i-1},y_{\text{\rm el},\tau}^{i-1},y_{{\rm vi},\tau}^{i-1})-\int_{t_{i-1}}^{t_i}\langle \dot{\ell}, y^{i-1}_{{\rm el},\tau}\circ y^{i-1}_{{\rm el},\tau}\rangle,
  \end{align*}
  where we used the definition of $G_{\tau}$ and of De Giorgi variational interpolation. Then, summing from $i=1$ to $i=n$ we get \eqref{energy equality deiscrete}.
\end{proof}

Before passing to the limit for $\tau \rightarrow 0$ in the energy
equality \eqref{energy equality deiscrete}, we need to characterize
the limit of the De Giorgi variational interpolation. In the following
Lemma we show that such limit  coincides with that  of the
backward interpolants.
\medskip

\begin{lemma}\label{lemma limit de giorgi}
  If $\overline{y}_{{\rm vi}, \tau}(t)\rightharpoonup y_{\rm vi}(t)$ in $W^{1, p_{\rm vi}}(\Omega;\R^d)$, then $\widetilde{y}_{{\rm vi}, \tau}(t)\rightharpoonup y_{\rm vi}(t)$ in $W^{1, p_{\rm vi}}(\Omega;\R^d)$.
\end{lemma}

\begin{proof}
  First, let us show that, for $\tau>0$ and $t\in (t^\tau_{i-1}, t^\tau_{i}]$ fixed, $\| \widetilde{y}_{{\rm vi}, \tau}(t)- \overline{y}_{{\rm vi}, \tau}(t)\|_{L^1(\Omega)}\leq c\tau^{p_{\psi}-1}$. We have, by definition of $\nabla \overline{y}_{{\rm vi}, \tau}$ and H\"older inequality,
  \begin{align*}
    \|\nabla \widetilde{y}_{{\rm vi}, \tau}(t)-\nabla \overline{y}_{{\rm vi}, \tau}(t)\|_{L^1(\Omega)}&\leq \|\nabla y_{{\rm vi},\tau}^{i-1}(\nabla y_{{\rm vi},\tau}^{i-1})^{-1}\left(\nabla \widetilde{y}_{{\rm vi}, \tau}(t)-\nabla {y}_{{\rm vi}, \tau}^i\right)\|_{L^1(\Omega)}\\
    &\leq \|\nabla y_{{\rm vi},\tau}^{i-1}\|^{p_{\psi}'}_{L^{p_{\psi}'}(\Omega)}\|(\nabla y_{{\rm vi},\tau}^{i-1})^{-1}\left(\nabla \widetilde{y}_{{\rm vi}, \tau}(t)-\nabla {y}_{{\rm vi}, \tau}^i\right)\|^{p_{\psi}}_{L^{p_{\psi}}(\Omega)}.
  \end{align*}
Since $p_{\psi}\geq 2$ by \eqref{p_psi>2}  we have that  $p_{\psi}'\leq p_{\psi}$. Hence, by the boundedness of $\nabla y_{{\rm vi},\tau}^{i-1}$ in $L^{p_{\psi}}(\Omega;\R^{d\times d})$ and the fact that $\Omega$ is bounded, we have that $\|\nabla y_{{\rm vi},\tau}^{i-1}\|^{p_{\psi}'}_{L^{p_{\psi}'}(\Omega)}\leq c$ uniformly in $i$ and $\tau$. Thus, by growth condition \eqref{dissipation density}, we have 
\begin{align*}
  \|\nabla \widetilde{y}_{{\rm vi}, \tau}(t)-\nabla \overline{y}_{{\rm vi}, \tau}(t)\|_{L^1(\Omega)}
  &\leq c \|(\nabla y_{{\rm vi},\tau}^{i-1})^{-1}\left(\nabla {y}_{{\rm vi}, \tau}^i-\nabla {y}_{{\rm vi}, \tau}^{i-1}\right)\|_{L^{p_{\psi}}(\Omega)}^{p_{\psi}}\\
  &\quad+c\|(\nabla y_{{\rm vi},\tau}^{i-1})^{-1}\left(\nabla \widetilde{y}_{{\rm vi}, \tau}(t)-\nabla {y}_{{\rm vi}, \tau}^{i-1}\right)\|_{L^{p_{\psi}}(\Omega)}^{p_{\psi}}\\
  &\leq c \Psi\left( y_{{\rm vi},\tau}^{i-1},  {y}_{{\rm vi}, \tau}^i-{y}_{{\rm vi}, \tau}^{i-1} \right)+c \Psi\left(y_{{\rm vi},\tau}^{i-1},  \widetilde{y}_{{\rm vi}, \tau}-{y}_{{\rm vi}, \tau}^{i-1} \right)\\
  &\leq c \Psi\left( y_{{\rm vi},\tau}^{i-1},  {y}_{{\rm vi}, \tau}^i-{y}_{{\rm vi}, \tau}^{i-1} \right),
\end{align*}
where in the last inequality we used the definition of $\widetilde{y}_{{\rm vi}, \tau}$ and the monotonicity property \eqref{monotonicity diss}. Using the $p_{\psi}$-homogeneity \eqref{p-homogeneous} and the boundedness of the dissipation, we get
\begin{equation*}
  \|\nabla \widetilde{y}_{{\rm vi}, \tau}(t)-\nabla \overline{y}_{{\rm vi}, \tau}(t)\|_{L^1(\Omega)}
  \leq c\tau^{p_{\psi}} \Psi\left( y_{{\rm vi},\tau}^{i-1},  \frac{{y}_{{\rm vi}, \tau}^i-{y}_{{\rm vi}, \tau}^{i-1}}{\tau} \right)\leq c \tau^{p_{\psi}-1}.
\end{equation*}
Then $\| \widetilde{y}_{{\rm vi}, \tau}(t)- \overline{y}_{{\rm vi},
  \tau}(t)\|_{L^1(\Omega)}\leq c\tau^{p_{\psi}-1}$ follows since
$\widetilde{y}_{{\rm vi}, \tau}(t)$ and $ \overline{y}_{{\rm vi},
  \tau}(t)$ have zero  mean. 

The assertion follows as $\Omega$ is bounded, by assumption
$\overline{y}_{{\rm vi}, \tau}(t)\rightharpoonup y_{\rm vi}(t)$ in $W^{1, p_{\rm vi}}(\Omega;\R^d)$, and $\overline{y}_{{\rm vi}, \tau}(t)$ is bounded in $W^{1,p_{\rm vi}}(\Omega;\R^d)$ by coercivity, as shown in Section \ref{sec: estimates} for $\overline{y}_{{\rm vi}, \tau}(t)$.
\end{proof}

\subsection{Proof of the energy inequality}\label{sec:1}

In the following, we extract further subsequences without relabeling whenever necessary.

Assume to be given a sequence of partitions $(\Pi_{\tau})_{\tau}$ with $\tau\rightarrow 0$ and denote by $(y_{{\rm el}}^i,y_{{\rm vi}}^i)_{i=0}^{N}$ the corresponding incremental solutions.
The estimates in Section \ref{sec: estimates} and Lemma \ref{lemma limit de giorgi} ensure that for every $t\in [0,T]$
\begin{equation*}
  \begin{aligned}
    \overline{y}_{{\rm vi},\tau}(t)\rightharpoonup y_{{\rm vi}}(t)\quad &\text{ in }
 W^{1,p_{\rm vi}}(\Omega;\R^d), \\
 \widetilde{y}_{{\rm vi},\tau}(t)\rightharpoonup y_{{\rm vi}}(t)\quad &\text{ in }
 W^{1,p_{\rm vi}}(\Omega;\R^d), \\
    \overline{y}_{\tau}(t)\rightharpoonup y(t)\quad &\text{ in }
    W^{1,q}(\Omega;\R^d).
  \end{aligned}
\end{equation*}
Moreover, by Sobolev embedding
we have that $\hat y_{\text{\rm vi},\tau}\rightharpoonup y_{\text{\rm vi}}$ weakly in $C([0,T];W^{1,p_r}(\Omega;\R^d))$.

As regards the elastic deformation, given $t\in [0,T]$ by extracting a subsequence $(\tau^t_k)_{k\in \N}$ possibly depending on $t$ we get
\begin{equation*}
  \overline{y}_{{\rm el},\tau_k^t}(t)\rightharpoonup y_{{\rm el}}(t)\quad \text{ in } 
  W^{1,p_{\rm el}}( y_{\rm vi}(t,\Omega);\R^d).
\end{equation*}
Note that here we have to implement an exhaustion argument for dealing
with the moving domains $\overline{y}_{\rm vi}(t,\Omega)$, exactly as
in Section \ref{sec:time discr}. Moreover, the total deformation 
$y$ can be proved to fulfill  $y  =y_{\text{\rm el}}\circ y_{\text{\rm
    vi}} $ by  arguing as in Section \ref{Closure adm def}.

We aim at passing to the limit in the energy equality \eqref{energy
  equality deiscrete}, which can be rewritten, thanks to the
definition of $G_{\tau}(t)$ and of $\Psi^{+,i}_{ r }$, in the weaker form
\begin{align}
  \mathcal{E}(\overline{t}_{\tau}(t),\overline{y}_{\text{\rm el},\tau}(t),\overline{y}_{\text{\rm vi},\tau}(t))\,{+}
  \int_{0}^{\overline{t}_{\tau}(t)} \Psi\left(\underline{y}_{{\rm vi},\tau},\dot{\hat{y}}_{{\rm vi},\tau}   \right)
  +
  (p_{\psi}{-}1)\int_{0}^{\overline{t}_{\tau}(t)} \Psi\left(\underline{y}_{{\rm vi},\tau},\frac{\widetilde{y}_{{\rm vi},\tau}{-}\underline{y}_{{\rm vi},\tau}}{\tau} \right)
  \notag
  \\
  \leq \mathcal{E}(0,y_{{\rm el},0},y_{{\rm vi},0})-
  \int_{0}^{\overline{t}_{\tau}(t)} \langle\dot{\ell},\underline{y}_{\tau} \rangle.\label{discrete energy ineq 2}
\end{align}

Passing to the $\liminf$ in the left-hand side of inequality \eqref{discrete energy ineq 2}, we find by lower semicontinuity
\begin{equation*}
  \mathcal{E}(t,y_{\text{\rm el}}(t),y_{\text{\rm vi}}(t))\leq \liminf_{\tau\rightarrow 0}\mathcal{E}(\overline{t}_{\tau}(t),\overline{y}_{\text{\rm el},\tau}(t),\overline{y}_{\text{\rm vi},\tau}(t)).
\end{equation*}
Let us now study the first dissipation term in \eqref{discrete energy ineq 2}. The calculations for the second one are analogous by Lemma \ref{lemma limit de giorgi}.
Recalling that by definition $\overline{t}_{\tau}(t)\geq t $ and that $\psi\geq 0$, we have that
\begin{equation*}
  \liminf_{\tau\rightarrow 0}\int_{0}^{\overline{t}_{\tau}(t)}\Psi\left(\underline{y}_{{\rm vi},\tau},\dot{\hat{y}}_{{\rm vi},\tau}   \right)
  \geq \liminf_{\tau\rightarrow 0}\int_{0}^{t}\Psi{\left(\underline{y}_{{\rm vi},\tau},\dot{\hat{y}}_{{\rm vi},\tau}   \right)}.
\end{equation*}
Moreover, up to a subsequence, $\nabla \dot {\hat{y}}_{\text{\rm vi},\tau}(\nabla \underline{y}_{\text{\rm vi},\tau})^{-1}\rightharpoonup l$ weakly in
$L^{p_{\psi}}(\Omega;\R^{d\times d})$ by \eqref{coercivity diss p_psi}.
It hence remains to identify the limit $l$. To this end, let us define $$(t,\xi)\in [0,T]\times \underline{y}_{{\rm vi},\tau}(t,\Omega)\mapsto v_{\tau}(t,\xi):=\dot{\hat{y}}_{\text{\rm vi},\tau}(t,\underline{y}_{\text{\rm vi},\tau}^{-1}(\xi))\in \R^d.$$
By a pointwise-in-time change of variables we have
\begin{equation*}
  \int_0^{T}\!\!\int_{\underline{y}_{\text{\rm vi},\tau}(t,\Omega)}|v_{\tau}(t,\xi)|^{p_{\psi}}\,{\rm d}\xi{\rm d}t=  \int_0^{T}\!\!\int_{\Omega}|\dot{\hat{y}}_{\text{\rm vi},\tau}(t,X)|^{p_{\psi}}\,{\rm d}X{\rm d}t\leq c.
\end{equation*}
In order to obtain a bound on the gradient $\nabla v_{\tau}$, let us consider
\begin{equation*}
  \begin{aligned}
    \int_0^{T}\!\!\int_{\underline{y}_{\text{\rm vi},\tau}(t,\Omega)}|\nabla v_{\tau}(t,\xi)|^{p_{\psi}}\,{\rm d}\xi{\rm d}t&=\int_0^{T}\!\!\int_{\Omega}\left|\nabla \dot{\hat{y}}_{\text{\rm vi},\tau}(\nabla \underline{y}_{\text{\rm vi},\tau})^{-1}(t,X)\right|^{p_{\psi}}\,{\rm d}X{\rm d}t
    \\&\stackrel{\eqref{dissipation density}}{\leq} \int_0^{T}\Psi\left(\underline{y}_{{\rm vi},\tau},\dot{\hat{y}}_{{\rm vi},\tau}   \right)\leq c.
  \end{aligned}
\end{equation*}
For given $t_0\in (0,T)$, let us show that $\cap_{t\in [t_0,t_0+\delta]}y_{\text{\rm vi}}(t,\Omega)$ is not empty for small $\delta>0$.
Notice that, by Sobolev embedding, $\underline{y}_{\text{\rm vi},\tau}\rightarrow y_{\text{\rm vi}}$ in $C([0,T]\times \overline{\Omega})$. Hence, for every $\epsilon>0$, there exists $\overline{\tau}=\overline{\tau}(\epsilon)$ such that, for every $\tau\leq \overline{\tau}$, we have
\begin{equation*}
\sup_{X\in \overline{\Omega}}\Big| \underline{y}_{\text{\rm vi},\tau}(t,X)-y_{\text{\rm vi}}(t,X) \Big|\leq \frac{\epsilon}{2}.
\end{equation*}
Moreover, since $y_{\text{\rm vi}}$ is absolutely continuous in time, 
for $|t-s|<\nu$ and $\nu>0$ small we also have
\begin{equation*}
    \sup_{X\in \overline{\Omega}}\left| y_{\text{\rm vi}}(t,X)-y_{\text{\rm vi}}(s,X)  \right|\leq \frac{\epsilon}{2}.
\end{equation*}
Combining these two inequalities we get
\begin{equation*}
  \sup_{X\in \overline{\Omega}}\Big| \underline{y}_{\text{\rm vi},\tau}(t,X)-y_{\text{\rm vi}}(s,X)  \Big|\leq \epsilon
\end{equation*}
for $\tau$ and $\nu$ small enough.
We can hence fix $\omega\subset\subset \cap_{t\in [t_0,t_0+\nu]}\underline{y}_{\text{\rm vi}}(t,\Omega)$ and trivially extend $v_{\tau}$ on $\R^d\setminus \omega$.
Then, thanks to the bounds above we have that $v_{\tau}\rightharpoonup v$ weakly in  $L^{p_{\psi}}([0,T];W^{1,p_{\psi}}(\omega))$. We have to show that $v=\dot{y}_{\text{\rm vi}}\circ y_{\text{\rm vi}}^{-1} $.
In fact, we have
  \begin{align}
  \int_{t_0}^{t_0+\nu}\!\!\int_{\omega}v(t,\xi) \,{\rm d}\xi{\rm d}t \leftarrow  \int_{t_0}^{t_0+\nu}\!\!\int_{\omega}v_{\tau}(t,\xi)\,{\rm d}\xi{\rm d}t=&
  \int_{t_0}^{t_0+\nu}\!\!\int_{\underline{y}_{\text{\rm vi},\tau}^{-1}(t,\omega)}\dot{\hat{y}}_{\text{\rm vi},\tau}(t,X)\,{\rm d}X{\rm d}t\notag
  \\
  =\int_{t_0}^{t_0+\nu}\!\!\int_{\R^d}\dot{\hat{y}}_{\text{\rm vi},\tau}(t,X)\mathbbm{1}_{\underline{y}_{\text{\rm vi},\tau}^{-1}(t,\omega)}(t,X)\,{\rm d}X{\rm d}t
  \rightarrow&
    \int_{t_0}^{t_0+\nu}\!\!\int_{\R^d}\dot{y}_{\text{\rm vi}}(t,X)\mathbbm{1}_{y_{\text{\rm vi}}^{-1}(t,\omega)}(t,X)\,{\rm d}X{\rm d}t\notag\\
    &=  \int_{t_0}^{t_0+\nu}\!\!\int_{\omega}\dot{y}_{\text{\rm vi}}(t,y_{\text{\rm vi}}^{-1}(t,\xi))\,{\rm d}\xi{\rm d}t,\label{v=comp}
  \end{align}
where we have used that $\hat{y}_{\text{\rm vi},\tau}\rightharpoonup y_{\text{\rm vi}}$ weakly in $C([0,T];W^{1,p_{\text{\rm vi}}}(\Omega))$,  $\mathbbm{1}_{\underline{y}_{\text{\rm vi},\tau}^{-1}(\omega)}\rightarrow \mathbbm{1}_{y_{\text{\rm vi}}^{-1}(\omega)}$ strongly in $L^1(\omega)$,
and the fact that
$\mathbbm{1}_{\underline{y}_{\text{\rm vi},\tau}^{-1}(t,\omega)}$ is bounded.
Since in \eqref{v=comp} $t_0,\,\nu,
$ and $ \omega$ are arbitrary, we have that $v=\dot{y}_{\text{\rm vi}}\circ y_{\text{\rm vi}}^{-1}$ and  we have hence identified $l=\nabla v$.
By weak lower semicontinuity, we thus have that
\begin{equation*}
  \begin{aligned}
    \liminf_{\tau\rightarrow 0}\int_{0}^{\overline{t}_{\tau}(t)}\Psi\left(\underline{y}_{{\rm vi},\tau},\dot{\hat{y}}_{{\rm vi},\tau}   \right)&=
    \liminf_{\tau\rightarrow 0}\int_{0}^{\overline{t}_{\tau}(t)}\!\!
    \int_{\overline{y}_{{\rm vi},\tau}(s,\Omega)}\psi(\nabla v_{\tau}(s,\xi))\,{\rm d}\xi{\rm d}s\\
    &\geq \int_{0}^{t}\!\! \int_{{y}_{{\rm vi}}(s,\Omega)}\psi(\nabla v(s,\xi))\,{\rm d}\xi{\rm d}s=
    \int_{0}^{t}\Psi\left(y_{{\rm vi},\tau},\dot{y}_{{\rm vi},\tau}   \right).
  \end{aligned}
\end{equation*}
Thanks to the boundedness and to the weak lower semicontinuity of the energy and of the dissipation we can apply Helly's Selection Principle  \cite[Thm. B.5.13, p. 611]{mielke2015rate} and find a nondecreasing function $\delta:[0,T]\rightarrow [0,\infty)$ such that
\begin{subequations}
  \begin{align}
&\int_{0}^{t}\Psi\left(\underline{y}_{{\rm vi},\tau},\dot{\hat{y}}_{{\rm vi},\tau}   \right) \rightarrow \delta(t),  \label{Helly2}\\
&\int_{s}^{t}\Psi\left(\underline{y}_{{\rm vi},\tau},\dot{\hat{y}}_{{\rm vi},\tau}   \right) \leq \delta(t)-\delta(s) \label{Helly3}
\end{align}
\end{subequations}
for every $s,t\in [0,T]$.
Then, fixing $t\in [0,T]$, we have
\begin{equation*}
  \delta(t)\stackrel{\eqref{Helly2}}{=}\lim_{k\rightarrow \infty} \int_{0}^{t}\Psi\left(\underline{y}_{{\rm vi},\tau},\dot{\hat{y}}_{{\rm vi},\tau}   \right)\leq
  \liminf_{k\rightarrow \infty} \int_{0}^{\underline{t}_{\tau}(t)}\Psi\left(\underline{y}_{{\rm vi},\tau},\dot{\hat{y}}_{{\rm vi},\tau}   \right).
\end{equation*}

Setting $\theta_{\tau}(s):=-\langle\dot{\ell}(s),\underline{y}_{\tau} (s)\rangle$, by the regularity of $\ell$ and the boundedness of $(\underline{y}_{\tau}(t) )_{\tau}$ for almost every $t\in [0,T]$ in $
W^{1,q}(\Omega;\R^d)$,
we have that $(\theta_{\tau})_{\tau}$ is equiintegrable. Hence, we can apply the Dunford-Pettis Theorem (see, e.g., \cite[Thm. B.3.8, p. 598]{mielke2015rate}) to find a subsequence such that
\begin{equation*}
  \theta_{\tau}\rightharpoonup \theta \quad \text{ in } L^1(0,T).
\end{equation*}
Furthermore, thanks to the boundedness of the energy and the dissipation, we are able to find  further  $t$-dependent subsequences $(\tau^t_k)_{k\in \N}$ such that
  \begin{equation*}
    \theta_{\tau^t_k}\rightarrow \limsup_{\tau\rightarrow 0}\theta_{\tau}(t)=:\overline{\theta}(t),
  \end{equation*}
and, by regularity of $\ell$, that
\begin{equation*}
  \overline{\theta}(t):=\lim_{k\rightarrow \infty}\theta_{\tau^t_k}=\lim_{k\rightarrow \infty}\langle\dot{\ell}(t),\underline{y}_{\tau^t_k}(t) \rangle=\langle\dot{\ell}(t),y(t) \rangle.
\end{equation*}
In conclusion, passing to the $\liminf$
in the left-hand side and to the limit in the right-hand side of \eqref{discrete energy ineq 2} we retrieve energy inequality \eqref{energy inequality}.

\subsection{Proof of the semistability condition}\label{sec:2}


Fix now $t\in[0,T]$ and
recall that $\overline{y}_{{\rm el},\tau_k^t}(t)\rightharpoonup y_{{\rm el}}(t)$ in 
$W^{1,p_{\rm el}}( y_{\rm vi}(t,\Omega);\R^d)$.

By minimality of the incremental solution we have
\begin{equation*}
  \mathcal{E}\left(\overline{t}_{\tau}(t), \overline{y}_{\text{\rm
        el},\tau}(t), \overline{y}_{\text{\rm vi},\tau}(t)\right) \leq
  \mathcal{E}\left(\overline{t}_{\tau}(t), 
    \tilde{\overline{y}}_{\text{\rm el}},  \overline{y}_{\text{\rm vi},\tau}(t)\right)
\end{equation*}
for every   $\tilde{\overline{y}}_{\text{\rm el} }$ with $\left(\tilde{\overline{y}}_{\text{\rm el}}, \overline{y}_{\text{\rm vi}}(t)\right) \in \mathcal{A}$.
Let $(\tilde{y}_{\rm el}, y_{\rm vi}(t))\in \mathcal{A}$  be
given.
We want to show that  one can choose
$\tilde{\overline{y}}_{\text{\rm el},\tau}$ with
$\left(\tilde{\overline{y}}_{\text{\rm el},\tau}, \overline{y}_{\text{\rm
      vi}}(t)\right) \in \mathcal{A}$ in such a way that 
\begin{align}
  0&\leq \limsup_{ \tau \to 0}\Big(\mathcal{E}\left(\overline{t}_{\tau}(t), \tilde{\overline{y}}_{\text{\rm el},\tau}, \overline{y}_{\text{\rm vi},\tau}(t)\right) - \mathcal{E}\left(\overline{t}_{\tau}(t), \overline{y}_{\text{\rm el},\tau}(t), \overline{y}_{\text{\rm vi},\tau}(t)\right)\Big)\notag\\
  &\leq
  \mathcal{E}\left(t, \tilde{y}_{\text{\rm el}}, y_{\text{\rm vi}}(t)\right) - \mathcal{E}\left(t, y_{\text{\rm el}}(t), y_{\text{\rm vi}}(t)\right),\label{semistability proof}
\end{align}
which  would then imply    \eqref{semistability}. 

Since $\left(\tilde{y}_{\text{\rm el}}, y_{\text{\rm vi}}(t)\right)
\in \mathcal{A} $, we have that $y_{\text{\rm vi}}(\Omega)\in
\mathcal{J}_{\eta_1,\eta_2}$ and  $y_{\text{\rm vi}}(\Omega)$ is a Sobolev extension domain.
Hence, there exists  a linear and bounded  extension operator
$E:W^{1,p_{\text{\rm el}}}(y_{\text{\rm vi}}(\Omega);\R^d)$
$\rightarrow$ $ W^{1,p_{\text{\rm el}}}(\R^d;\R^d)$. We thus 
define $\tilde{\overline{y}}_{\text{\rm el},\tau}\in W^{1,p_{\rm
    el}}(\overline{y}_{\text{\rm vi},\tau}(\Omega);\R^d)$ as the
restriction to $\overline{y}_{\text{\rm vi},\tau}(\Omega)$ of the extension
$E\tilde{y}_{\text{\rm el}}$, namely, 
\begin{equation*}
  \tilde{\overline{y}}_{\text{\rm el},\tau}:=E\tilde{y}_{\text{\rm el}}\Big|_{\overline{y}_{\text{\rm vi},\tau}(\Omega)}.
\end{equation*}
  In the
following, we just concentrate our attention on the stored
elastic energy part, since the treatment of the loading term is
immediate.   We write
\begin{align}
  \int_{\overline{y}_{\text{\rm vi},\tau}(t,\Omega)}W_{\text{\rm el}}(\nabla \tilde{\overline{y}}_{\text{\rm el},\tau})\,{\rm d}\xi-
  \int_{\overline{y}_{\text{\rm vi},\tau}(t,\Omega)}W_{\text{\rm el}}(\nabla {\overline{y}}_{\text{\rm el},\tau})\,{\rm d}\xi
  =\int_{\overline{y}_{\text{\rm vi},\tau}(t,\Omega)\cap y_{\rm vi} (t,\Omega)}W_{\text{\rm el}}(\nabla \tilde{\overline{y}}_{\text{\rm el},\tau})\,{\rm d}\xi\notag\\
  +\int_{\overline{y}_{\text{\rm vi},\tau}(t,\Omega)\setminus y_{\rm vi} (t,\Omega)}W_{\text{\rm el}}(\nabla \tilde{\overline{y}}_{\text{\rm el},\tau})\,{\rm d}\xi
  \quad-\int_{\overline{y}_{\text{\rm vi},\tau}(t,\Omega)}W_{\text{\rm el}}(\nabla {\overline{y}}_{\text{\rm el},\tau})\,{\rm d}\xi\label{ineq999}
\end{align}
By the growth condition \eqref{growthel} on $W_{\text{\rm el}}$ and the fact that on the set $\overline{y}_{\text{\rm vi},\tau}(t,\Omega)$
we have $\tilde{\overline{y}}_{\text{\rm el},\tau}=\tilde{y}_{\text{\rm el},\tau}$, which is uniformly bounded in $ W^{1,p_{\rm el}}(y_{{\rm vi},\tau}(t,{\Omega});\R^d)$, we find
\begin{equation*}
  \int_{\overline{y}_{\text{\rm vi},\tau}(t,\Omega)\setminus y_{\rm vi} (t,\Omega)}W_{\text{\rm el}}(\nabla \tilde{\overline{y}}_{\text{\rm el},\tau})\,{\rm d}\xi=
  \int_{\tilde{y}_{\text{\rm vi},\tau}(t,\Omega)\setminus y_{\rm vi} (t,\Omega)}W_{\text{\rm el}}(\nabla  {\tilde{y}}_{\text{\rm el},\tau})\,{\rm d}\xi
  \stackrel{\eqref{growthel}}{\leq} c\left|\overline{y}_{\text{\rm vi},\tau}(t,\Omega)\setminus y_{\rm vi} (t,\Omega)\right|.
\end{equation*}
Since the measure of the set $\overline{y}_{\text{\rm
    vi},\tau}(t,\Omega)\setminus y_{\text{\rm vi}}(t,\Omega) $
vanishes as $\tau$ goes to $0$ by the uniform convergence of
$\overline{y}_{\text{\rm vi},\tau}$  to ${y}_{\text{\rm vi}}$,  we have
\begin{equation*}
  \lim_{\tau\rightarrow 0}  \int_{\overline{y}_{\text{\rm vi},\tau}(t,\Omega)}W_{\text{\rm el}}(\nabla \tilde{\overline{y}}_{\text{\rm el},\tau})\,{\rm d}\xi
  =\int_{y_{\text{\rm vi}}(t,\Omega)}W_{\text{\rm el}}(\nabla \tilde{y}_{\text{\rm el}})\,{\rm d}\xi.
\end{equation*}
We can hence pass to the $\limsup$ in inequality \eqref{ineq999}
 as $\tau \to 0$  and obtain \eqref{semistability proof}, which is nothing but the semistability \eqref{semistability}.

\section{Linearization: Proof of Theorem \ref{thm limit lin}}\label{sec:Proof of Linearization}

 We first prove in Subsection \ref{sec:co} some coercivity results, uniform with respect to the
linearization parameter $\varepsilon$, which in turn provide  a priori
estimates on the sequence of approximable solutions
$(u_\varepsilon,v_\varepsilon)_\varepsilon$. Then, we check in
Subsection \ref{sec:ga} some $\Gamma$-$\liminf$ inequalities for the
energy and the dissipation. Eventually, in Subsection
\ref{sec:con} we show that the approximable solutions
$(u_\varepsilon,v_\varepsilon)_\varepsilon$ converge, up to subsequences, to solutions of
the linearized problem in the sense of Theorem  \ref{thm limit lin}.

In the following, we use the notation 
\begin{equation*}
  W^{\varepsilon}_{\rm el}(A):=\frac{1}{\varepsilon^2}W_{\rm el}(I+\varepsilon A), \quad \widetilde{W}^{\varepsilon}_{\rm vi}(A):=\frac{1}{\varepsilon^2}\widetilde{W}_{\rm vi}(I+\varepsilon A), \quad \psi^{\varepsilon}(A):=\frac{1}{\varepsilon^2}\psi(\varepsilon A)
\end{equation*}
for the rescaled energy and dissipation densities.

\subsection{Coercivity}\label{sec:co}

 We devote this subsection to the proof of the following. 
\medskip

\begin{lemma}[Coercivity]\label{lemma coercivity}
  For every $(u,v)\in \widetilde{\mathcal{A}}_{\varepsilon}$, it holds
  \begin{equation*}
    \|u\|^2_{H^1(\Omega)}+\|v\|^2_{H^1(\Omega)}+\|\nabla \dot  v\|^2_{L^2(\Omega)}+\|\varepsilon \nabla v\|_{L^{\infty}(\Omega)}\leq c{\left( 1+
    \mathcal{W}^\varepsilon _{\text{\rm vi}}(v)+\mathcal{W}^\varepsilon _{\text{\rm el}}(u,v)+\Psi^\varepsilon(v)   \right)}.
  \end{equation*}
\end{lemma}

  Notice   the bound on the term $\|\varepsilon \nabla
  v\|_{L^{\infty}(\Omega)}$, which follows from assumption {\rm
    \ref{L4}}. This bound will play an important role in passing to
  the limit for $\varepsilon \to 0$.  

  \begin{proof}[Proof of Lemma \ref{lemma coercivity}]
     With no loss of generality we  can assume $\mathcal{W}^\varepsilon _{\text{\rm vi}}(v)+\mathcal{W}^\varepsilon _{\text{\rm el}}(u,v)+\Psi^\varepsilon(v)<\infty$

 By assumption \ref{L4} we have that $I+\varepsilon \nabla v \in K
 \subset \subset SL(d)$  almost everywhere in  $ \Omega$.
By using \eqref{ck} we get that $|I+\varepsilon \nabla v|\leq c_K$,
 hence 
\begin{equation*}
\|\varepsilon \nabla v\|_{L^{\infty}(\Omega)}\leq c.
\end{equation*}

Since $v$ has zero  mean  by assumption, by applying the
Poincar\'e-Wirtinger inequality and  by taking into account the
  growth condition \eqref{b_vi} we get
\begin{equation*}
  \|v\|^2_{H^1(\Omega)}\leq c \|\nabla v\|^2_{L^2(\Omega)}=\frac{c}{\varepsilon^2} \|\varepsilon\nabla v\|^2_{L^2(\Omega)}\leq
   \frac{c}{\varepsilon^2}\int_{\Omega} W_{\rm vi}(I+\varepsilon \nabla v){\rm d}X =c\mathcal{W}^{\varepsilon}_{\rm vi }(v).
\end{equation*}

 Using condition \eqref{a_psi} and the fact that $|I+\varepsilon \nabla v|$ is bounded in $L^{\infty}$ we get
\begin{equation*}
  \begin{aligned}
    \|\nabla   \dot{v}\|^2_{L^2(\Omega)}=\frac{1}{\varepsilon^2}\int_{\Omega}|\varepsilon\nabla  \dot{v}|^2 {\rm d}X&\leq
    \frac{1}{\varepsilon^2}\int_{\Omega}|\varepsilon \nabla  \dot{v}|^2|I+\varepsilon \nabla v|^{-2}|I+\varepsilon \nabla v|^{2} {\rm d}X\\
    &\leq\frac{c}{\varepsilon^2}\int_{\Omega}\psi\left(\varepsilon\nabla \dot{v} (I+\varepsilon \nabla v)^{-1} \right){\rm d}X=c\Psi^{\varepsilon}(v).
  \end{aligned}
\end{equation*}

 In order to obtain the $H^1$-bound on $u$, we start by fixing $Q\in
 SO(d)$ and  define   $F_{\rm el}:=\nabla y (I+\varepsilon \nabla v)^{-1}$, where we recall that $y=\operatorname{id}+\varepsilon u$. We have
\begin{equation*}
  \begin{aligned}
    |\nabla y - Q|^2&=|\nabla y - Q(I+\varepsilon \nabla v)+\varepsilon Q\nabla v|^2=|(F_{\rm el}-Q)(I+\varepsilon \nabla v)+\varepsilon Q\nabla v|^2\\
    &\leq c(|F_{\rm el}-Q|^2|I+\varepsilon \nabla v|^2+\varepsilon^2 |\nabla v|^2)\leq c(|F_{\rm el}-Q|^2+\varepsilon^2 |\nabla v|^2).
  \end{aligned}
\end{equation*}
Taking the infimum over $Q\in SO(d)$ we get
\begin{equation*}
  \operatorname{dist}^2(\nabla y,SO(d))\leq c( \operatorname{dist}^2(F_{\rm el},SO(d))+\varepsilon^2|\nabla v|^2).
\end{equation*}
We  now   integrate over $\Omega$ and, thanks to assumption \eqref{b_el} and the estimate on $\|v\|^2_{H^{1}(\Omega)}$, we find
\begin{equation*}
  \int_{\Omega}\operatorname{dist}^2(\nabla y,SO(d)){\rm d}X \leq c \int_{\Omega}W_{\rm el}(F_{\rm el}){\rm d}X+c\varepsilon^2\|\nabla v\|^2_{L^2(\Omega)}\leq c\varepsilon^2\left(\mathcal{W}^{\varepsilon}_{\rm el}(u,v)+\mathcal{W}^{\varepsilon}_{\rm vi}(v)\right).
\end{equation*}
The classical Rigidity Estimate \cite[Theorem $\rm 3.1$]{Rigidity}
implies that there exists  a constant rotation  $\hat Q\in SO(d)$ such that
\begin{equation*}
  \|\nabla y-\hat Q\|^2_{L^2(\Omega)}\leq c\|\operatorname{dist}(\nabla y,SO(d))\|^2_{L^2(\Omega)}.
\end{equation*}
We hence have that
\begin{equation*}
  \|\nabla y-\hat Q\|^2_{L^2(\Omega)}\leq c\varepsilon^2\left(\mathcal{W}^{\varepsilon}_{\rm el}(u,v)+\mathcal{W}^{\varepsilon}_{\rm vi}(v)\right).
\end{equation*}
Recalling that $y=\operatorname{id}$ on $\Gamma_D$, by \cite[(3.14)]{DalMaso} we also deduce
\begin{equation*}
  \|I-\hat Q\|^2_{L^2(\Omega)}\leq c\varepsilon^2\left(\mathcal{W}^{\varepsilon}_{\rm el}(u,v)+\mathcal{W}^{\varepsilon}_{\rm vi}(v)\right).
\end{equation*}
In conclusion, we get that
\begin{align*}
  \|\nabla u\|^2_{L^2(\Omega)}&=\frac{1}{\varepsilon^2}\|\nabla y -I\|^2_{L^2(\Omega)}\leq \frac{2}{\varepsilon^2}\|\nabla y -\hat Q\|^2_{L^2(\Omega)} +\frac{2}{\varepsilon^2}\|\hat Q -I\|^2_{L^2(\Omega)}\\
  &\leq c\left(\mathcal{W}^{\varepsilon}_{\rm el}(u,v)+\mathcal{W}^{\varepsilon}_{\rm vi}(v)\right)
\end{align*}
 whence the assertion follows.  \end{proof}

\subsection{ \texorpdfstring{${\Gamma}$-$\liminf$}{} inequalities}\label{sec:ga}

In order to proceed with the linearization, we need to establish $\Gamma$-$\liminf$ inequalities.
 At first,  we prove the following Lemma on the convergence of
the densities.
\medskip

\begin{lemma}[Convergence of the densities]
Assume conditions \ref{L3}, \ref{L6}, and \ref{L9}. Then, we have
\begin{equation*}
  W^{\varepsilon}_{\rm el}\rightarrow |\cdot|^2_{\C_{\rm el}}, \quad\quad
  \widetilde{W}^{\varepsilon}_{\rm vi}\rightarrow |\cdot|^2_{\C_{\rm vi}},
  \quad\quad
  \psi^{\varepsilon}\rightarrow |\cdot|^2_{\D}
\end{equation*}
locally uniformly. Moreover, we have
  \begin{align}
    &|z|^2_{\C_{\rm el}}\leq \inf\left\{ \liminf_{\varepsilon\rightarrow 0 }W^{\varepsilon}_{\rm el}(z_{\varepsilon})\;|\; z_{\varepsilon} \rightarrow z \text{ in } \R^{d\times d}  \right\},
    \label{ineq inf densities el}\\
    &|z|^2_{\C_{\rm vi}}\leq \inf\left\{ \liminf_{\varepsilon\rightarrow
    0 }\widetilde{W}^{\varepsilon}_{\rm vi}(z_{\varepsilon})\;|\; z_{\varepsilon} \rightarrow z \text{ in } \R^{d\times d}  \right\},\label{ineq inf densities vi}\\
    &|z|^2_{\D}\leq \inf\left\{ \liminf_{\varepsilon\rightarrow 0 }\psi^{\varepsilon}(z_{\varepsilon})\;|\; z_{\varepsilon} \rightarrow z \text{ in } \R^{d\times d}  \right\}.\label{ineq inf densities D}
  \end{align} 
\end{lemma}

\begin{proof}
  Let $K_0\subset \subset \R^{d\times d}$  be given.  Fix
  $\delta >0$ and let $c_{\rm el}(\delta)$ be the corresponding
  constant  from  assumption \eqref{c_el}. Then, for sufficiently small $\varepsilon$, we have that $\varepsilon K_0 \subset B_{c_{\rm el}(\delta)}(0)$. Hence, by \eqref{c_el} we find
  \begin{equation*}
    \limsup_{\varepsilon \rightarrow 0}\sup_{K_0}\left|W^{\varepsilon}_{\rm el}(\cdot)-|\cdot|^2_{\C_{\rm el}}\right|\leq \delta \sup_{K_0}|\cdot|^2_{\C_{\rm el}}\leq \delta c.
  \end{equation*}
  Since $\delta$ is arbitrary, we get local uniform convergence for $W^{\varepsilon}_{\rm el}$. For $\widetilde{W}^{\varepsilon}_{\rm vi}$ and $\psi^{\varepsilon}$ proof of convergence is analogous, using the corresponding conditions \eqref{c_vi} and \eqref{b_psi}, respectively.

For the $\Gamma$-$\liminf$ inequalities \eqref{ineq inf densities
  el}-\eqref{ineq inf densities D}, let $(z_{\varepsilon})_{\varepsilon}\subset \R^{d\times d}$ be such that $z_{\varepsilon}\rightarrow z$ in $\R^{d\times d}$.
  Assume without loss of generality that
  $\sup_{\varepsilon}W^{\varepsilon}_{\rm
    el}(z_{\varepsilon})<\infty.$ Then, the inequality follows from
  local uniform convergence. The same applies  to   $\widetilde{W}^{\varepsilon}_{\rm vi}$
  and $\psi^{\varepsilon}$.
\end{proof}

We are now  in the position of proving the  $\Gamma $-$\liminf$ inequalities for the
functionals.
\medskip

\begin{lemma}[$\Gamma $-$\liminf$ inequalities]\label{gamma liminf ineq}
  For every $(u,v)\in {\widetilde{\mathcal{A}}}_{\varepsilon}$, we have
  \begin{equation*}
    \begin{aligned}
      \mathcal{W}^{0}_{\rm el}(u,v)+  \mathcal{W}^{0}_{\rm vi}(v)\leq \inf \Big\{\liminf_{\varepsilon \rightarrow 0}\big(&
      \mathcal{W}^{\varepsilon}_{\rm el}(u_{\varepsilon},v_{\varepsilon})+  \mathcal{W}^{\varepsilon}_{\rm vi}(v_{\varepsilon})\big)\;\\
      &\big|\; (u_{\varepsilon},v_{\varepsilon})\rightharpoonup (u,v) \text{ weakly in }H^1(\Omega;\R^{d})^2\Big\},
    \end{aligned}
  \end{equation*}
  \begin{equation*}
    \int_0^t\Psi^{0}(\dot{v})\leq \inf \left\{\liminf_{\varepsilon \rightarrow 0}
    \int_0^t\Psi^{\varepsilon}(v_{\varepsilon},\dot{v}_{\varepsilon})\;\big|\; v_{\varepsilon}\rightharpoonup v \text{ weakly in }H^1([0,t]; L^2(\Omega;\R^{d}))
    \right\}.
  \end{equation*}
\end{lemma}
\begin{proof}
  Let $(u_{\varepsilon},v_{\varepsilon})\rightharpoonup (u,v) \text{ weakly in }H^1(\Omega;\R^{d})^2$ and assume without loss of generality that
  $\sup_{\varepsilon}\left(\mathcal{W}^{\varepsilon}_{\rm el}(u_{\varepsilon},v_{\varepsilon})+  \mathcal{W}^{\varepsilon}_{\rm vi}(v_{\varepsilon})\right)<\infty$.

  Thanks to  inequality \eqref{ineq inf densities vi} and \cite[Lemma
  4.2]{MielkeStefanelli} we immediately  handle the stored viscous
  energy terms as  
  \begin{equation*}
    \int_{\Omega}|\nabla v|^2_{\C_{\rm vi}}{\rm d}X\leq \liminf_{\varepsilon \rightarrow 0}\int_{\Omega}W^{\varepsilon}_{\rm vi}(\nabla v_{\varepsilon}){\rm d} X= \liminf_{\varepsilon \rightarrow 0}\frac{1}{\varepsilon^2}
    \int_{\Omega}W_{\rm vi}(I+\varepsilon \nabla v_{\varepsilon}){\rm d} X.
  \end{equation*}


  The  treatment of the stored  elastic energy term requires some steps. First, notice that, since $\sup_{\varepsilon}\mathcal{W}^{\varepsilon}_{\rm vi}(v_{\varepsilon})<\infty$,
  we have that $I+\varepsilon \nabla v_{\varepsilon}\in K$ almost everywhere in $\Omega$.
  Hence, $\|\varepsilon \nabla
  v_{\varepsilon}\|_{L^{\infty}(\Omega)}\leq c $ uniformly in
  $\varepsilon$ and $(I+\varepsilon \nabla v_{\varepsilon})^{-1}$ is bounded in $L^{\infty}(\Omega;\R^{d\times d})$ by \eqref{ck} as well.

  Let us then define the auxiliary tensor $Z_{\varepsilon}$ as
  \begin{equation*}
    Z_{\varepsilon}:=\frac{1}{\varepsilon}\left( (I+\varepsilon \nabla v_{\varepsilon})^{-1}-I+\varepsilon \nabla v_{\varepsilon} \right)=\varepsilon (I+\varepsilon \nabla v_{\varepsilon})^{-1}(\nabla v_{\varepsilon})^2
  \end{equation*}
  so that $(I+\varepsilon \nabla v_{\varepsilon})^{-1}=I-\varepsilon \nabla v_{\varepsilon}+ \varepsilon Z_{\varepsilon}$.
  Notice that $\|\varepsilon Z_{\varepsilon}\|_{L^{\infty}(\Omega)}\leq c$ since $\|\varepsilon \nabla v_{\varepsilon}\|_{L^{\infty}(\Omega)}\leq c$. Furthermore,
  \begin{equation*}
    \| Z_{\varepsilon}\|_{L^{1}(\Omega)}\leq \varepsilon \|(I+\varepsilon \nabla v_{\varepsilon})^{-1}(\nabla v_{\varepsilon})^2\|_{L^1(\Omega)}\leq c\varepsilon \|\nabla v_{\varepsilon}\|_{L^2(\Omega)}\leq c\varepsilon.
  \end{equation*}
  Hence, $Z_{\varepsilon}$ is bounded in $L^2(\Omega;\R^{d\times d})$ by interpolation, namely,
  \begin{equation*}
    \| Z_{\varepsilon}\|_{L^{2}(\Omega)}\leq \|\varepsilon
    Z_{\varepsilon}\|_{L^{\infty}(\Omega)}^{ 1/2}\frac{\|Z_{\varepsilon}\|_{L^{1}(\Omega)}^{ 1/2}}{\varepsilon^{ 1/2}} \leq c.
  \end{equation*}
  We therefore  conclude that   $Z_{\varepsilon}\rightharpoonup 0$ weakly in $L^2(\Omega;\R^{d\times d})$. 
  
  Define now $F^{\varepsilon}_{\rm el}:=(I+\varepsilon \nabla u_{\varepsilon})(I+\varepsilon \nabla v_{\varepsilon})^{-1}$ and
  \begin{equation*}
    A^{\varepsilon}:= \frac{F^{\varepsilon}_{\rm el}-I}{\varepsilon}= \frac{1}{\varepsilon}\left((I+\varepsilon \nabla u_{\varepsilon})(I+\varepsilon \nabla v_{\varepsilon})^{-1}-I\right).
  \end{equation*}
  We want to show that $A^{\varepsilon}\rightharpoonup \nabla u-\nabla v$ weakly in $L^{2}(\Omega;\R^{d\times d})$.
 Let us compute
  \begin{equation*}
    A^{\varepsilon}= \frac{1}{\varepsilon}\left((I+\varepsilon \nabla u_{\varepsilon})(I-\varepsilon \nabla v_{\varepsilon}+\varepsilon Z_{\varepsilon})-I\right)=\nabla u_{\varepsilon}-\nabla v_{\varepsilon}+Z_{\varepsilon}-\varepsilon (\nabla u_{\varepsilon} \nabla v_{\varepsilon}- \nabla u_{\varepsilon}Z_{\varepsilon}).
  \end{equation*}
  Since $\nabla u_{\varepsilon}-\nabla v_{\varepsilon}\rightharpoonup \nabla u-\nabla v $ and $Z_{\varepsilon}\rightharpoonup 0 $ weakly in $L^2(\Omega;\R^{d\times d})$, it remains to show that $H_{\varepsilon}:=\varepsilon (\nabla u_{\varepsilon} \nabla v_{\varepsilon}- \nabla u_{\varepsilon}Z_{\varepsilon})\rightharpoonup 0$ weakly in $L^2(\Omega;\R^{d\times d})$.
  Notice that $\|H_{\varepsilon}\|_{L^{2}(\Omega)}\leq c$ since $\nabla u_{\varepsilon}$ is bounded in $L^2(\Omega;\R^{d\times d})$ and $\varepsilon \nabla v_{\varepsilon}$ and $ \varepsilon Z_{\varepsilon}$ are bounded in $L^{\infty}(\Omega;\R^{d\times d})$.
  Moreover, since $\nabla v_{\varepsilon}$ and $ Z_{\varepsilon}$ are bounded in $L^{2 }(\Omega;\R^{d\times d})$, then $\|H_{\varepsilon}\|_{L^{1}(\Omega)}\leq c\varepsilon $ so that $H_{\varepsilon}\rightharpoonup 0$ weakly in $L^2(\Omega;\R^{d\times d})$.

  Hence, we have by \eqref{ineq inf densities el} and \cite[Lemma 4.2]{MielkeStefanelli} that
  \begin{align*}
    \int_{\Omega}|\nabla u -\nabla v|^2_{\C_{\rm el}}{\rm d}X&\leq
                                                               \liminf_{\varepsilon
                                                               \rightarrow
                                                               0}\int_{\Omega}W^{\varepsilon}_{\rm
                                                               el}(A^{\varepsilon})\,{\rm
                                                               d} X\\
    &= \liminf_{\varepsilon \rightarrow 0}\frac{1}{\varepsilon^2}
    \int_{\Omega}W_{\rm el}\left((I{+}\varepsilon \nabla u_{\varepsilon})(I{+}\varepsilon \nabla v_{\varepsilon})^{-1}\right){\rm d} X.
  \end{align*}


  Let $(v_{\varepsilon})_\varepsilon$ be such that $ v_{\varepsilon} \rightharpoonup v$ weakly in $H^1([0,t];L^2(\Omega;\R^{d}))$, and $\sup_{\varepsilon}\Psi^{\varepsilon}(v_{\varepsilon},\dot{v}_{\varepsilon})<\infty$.
  By coercivity we have that (up to a non relabeled subsequence) $\nabla \dot  v_{\varepsilon}(I+\varepsilon \nabla v_{\varepsilon})^{-1}\rightharpoonup Y$ weakly in $L^2(\Omega)$. We want to identify the limit as $Y=\nabla \dot  v$. First, notice that
  \begin{equation*}
    \nabla \dot  v_{\varepsilon}(I+\varepsilon \nabla v_{\varepsilon})^{-1}-\nabla \dot  v_{\varepsilon}=-\varepsilon \nabla  \dot{v}_{\varepsilon}\nabla v_{\varepsilon}(I+\varepsilon \nabla v_{\varepsilon})^{-1}=:Y_{\varepsilon}\rightharpoonup 0 \,\text{ weakly in } L^2(\Omega;\R^{d\times d}).
  \end{equation*}
  Indeed, $\| Y_{\varepsilon}\|_{L^{1}(\Omega)}\leq c\varepsilon$ and $\| Y_{\varepsilon}\|_{L^{2}(\Omega)}\leq c$. 
  Then the $\Gamma$-$\liminf$ inequality for the dissipation term
  follows from \eqref{ineq inf densities D} and \cite[Lemma
  4.2]{MielkeStefanelli}  applied  on the
  domain  $[0,t]\times \Omega $. 
\end{proof}

\subsection{Convergence of approximable solutions}\label{sec:con}

  Thanks to Lemma \ref{lemma coercivity} and the energy inequality \eqref{energy inequality epsilon} we have
    \begin{align}
       \|u_{\varepsilon}(t)\|^2_{H^1(\Omega)}&\leq
                                               c\left(1+\E^{\varepsilon}(t,u_{\varepsilon},v_{\varepsilon})\right)
                                               \leq c\left(1+\E^{\varepsilon}(t,u_{\varepsilon},v_{\varepsilon})+\int_{0}^t\Psi^{\varepsilon}(v_{\varepsilon},\dot{v}_{\varepsilon})\right)\notag \\
      &\leq
      c\left(1+\E^{\varepsilon}(u^0_{\varepsilon},v^0_{\varepsilon})+\int_{0}^t \langle \dot{\ell}^{\varepsilon}, u_{\varepsilon} \rangle \right).\label{u epsilon bound}
    \end{align}
By the Gronwall Lemma \cite[Lemma C.2.1, p. 534]{Kruzik Roubicek} this implies that $\|u_{\varepsilon}(t)\|_{H^1(\Omega)}\leq c$ for every $t\in [0,T]$.

Concerning $v_{\varepsilon}$ we similarly deduce from Lemma \ref{lemma coercivity} that
\begin{equation*}
     \|v_{\varepsilon}(t)\|^2_{H^1(\Omega)}\leq  c\left(1+\E^{\varepsilon}(t,u_{\varepsilon},v_{\varepsilon})\right)\leq
    c\left(1+\E^{\varepsilon}(u^0_{\varepsilon},v^0_{\varepsilon})+\int_{0}^t \langle \dot{\ell}^{\varepsilon}, u_{\varepsilon} \rangle \right)
\end{equation*}
so that $\|v_{\varepsilon}(t)\|_{H^1(\Omega)}\leq c$ for every $t\in [0,T]$, as well.

Again, by Lemma \ref{lemma coercivity} we have
$\|\nabla \dot  v_{\varepsilon}(t)\|^2_{L^{2}(\Omega)}\leq
c\Psi^{\varepsilon}(v_{\varepsilon}(t),\dot{v}_{\varepsilon}(t))$ for
every $t\in [0,T]$.
 This yields 
\begin{equation*}
    \int_0^t \|\nabla \dot  v_{\varepsilon}\|^2_{L^{2}(\Omega)}
    \leq c\int_0^t \Psi^{\varepsilon}(v_{\varepsilon},\dot{v}_{\varepsilon}) 
   \leq
   c\left(1+\E^{\varepsilon}(u^0_{\varepsilon},v^0_{\varepsilon})+\int_{0}^t \langle \dot{\ell}^{\varepsilon}, u_{\varepsilon} \rangle \right),
\end{equation*}
hence, $\|\nabla \dot  v_{\varepsilon}\|_{L^{2}(0,T;L^2(\Omega))}\leq c$.
Therefore, up to a non relabeled subsequence, we find
\begin{equation*}
   v_{\varepsilon}(t)\rightharpoonup v(t) \text{ in } H^{1}(\Omega;\R^d), \quad
\nabla \dot  v_{\varepsilon}(t)\rightharpoonup \nabla \dot  v(t) \text{ in } L^{2}(\Omega;\R^{d\times d}) \quad\end{equation*}
for almost every $t\in [0,T]$.
    Notice that, since $\E^{\varepsilon}(t,u_{\varepsilon}(t),v_{\varepsilon}(t))<\infty$ for every $t\in [0,T]$, from assumption \ref{L4} it follows $
  I+\varepsilon \nabla v_{\varepsilon}\in K$ for almost every $x \in
  \Omega$ and for every $t\in [0,T]$. In particular, $\varepsilon
  \nabla v_{\varepsilon}$ are uniformly bounded. Since $v_\varepsilon
  \in \mathcal{A}_{\varepsilon}$ by developing the determinant
  as a third-order polynomial we get
  $$1=\det(I+\varepsilon \nabla v_{\varepsilon})=1+\varepsilon
  \operatorname{tr}\nabla v_{\varepsilon}+
  \varepsilon^2\operatorname{tr}(\operatorname{cof}\nabla
  v_{\varepsilon})+ \varepsilon^3 \det \nabla v_{\varepsilon}+{\rm o}(\varepsilon^4).$$
  By using $\|\nabla v_{\varepsilon}(t)\|_{L^2(\Omega)}\leq c$ and
  $\varepsilon\|\nabla v_{\varepsilon}(t)\|_{L^\infty(\Omega)}\leq c$
  for a.e.  $t\in (0,T)$ we hence conclude that
  \begin{align*}
    \|\operatorname{tr}\nabla v_\varepsilon(t) \|_{L^1(\Omega)} \leq  \varepsilon \|\operatorname{tr}(\operatorname{cof}\nabla
  v_{\varepsilon}(t)) \|_{L^1(\Omega)} +\varepsilon^2\| \det \nabla
    v_\varepsilon(t) \|_{L^1(\Omega)} \leq c\varepsilon
  \end{align*}
  for a.e. $t\in (0,T)$. By passing to the limit as $\varepsilon \to
  0$, this ensures that  $\operatorname{tr}\nabla v = 0$ a.e.

Fix now $t\in[0,T]$. By \eqref{u epsilon bound} we have
  \begin{equation}\label{u epsilon subseq}
    u_{\varepsilon}(t)\rightharpoonup u(t) \text{ in } H^{1}(\Omega;\R^d), \quad\end{equation}
  where at this point the subsequence above may in general depend on
  $t$. However, we shall see that this is not the case by uniqueness
  of the limit (see  below).


The  linearized  energy inequality \eqref{energy inequality 0}
follows immediately from the  energy inequality \eqref{energy
  inequality epsilon} at level $\varepsilon$,   thanks to the
$\liminf$-inequalities in Lemma \ref{gamma liminf ineq} and  to
 the continuity of $\dot{\ell}$.


The  linearized  semistability condition \eqref{semistability
  0} on the other hand is more delicate, since it requires passing to
the $\limsup$ on the right-hand side of the  semistability
condition  \eqref{semistability epsilon}  by  choosing a suitable recovery sequence $\widetilde{u}_{\varepsilon}$.
In the following, we will  drop the indication of the time
dependence (note that time is fixed in this statement) and  simply
denote
$u_{\varepsilon}(t)=u_{\varepsilon},v_{\varepsilon}(t)=v_{\varepsilon},
u(t)=u$, and $v(t)=v$, to simplify notation.

 We start by showing that, for all fixed $\hat{u}\in
H^{1}_{\Gamma_D}(\Omega;\R^d) $ one can choose a recovery sequence
$(\tilde u_{\varepsilon})_{\varepsilon}$ such that 
\begin{equation}
  0\stackrel{\eqref{semistability epsilon}}{\leq} \limsup_{\varepsilon \rightarrow 0}\left (\mathcal{W}^{\varepsilon}_{\text{\rm el}}(\widetilde{u}_{\varepsilon},v_{\varepsilon}) - \mathcal{W}^{\varepsilon}_{\text{\rm el}}(u_{\varepsilon},v_{\varepsilon}) \right)
  \leq
  \mathcal{W}^{0}_{\text{\rm el}}(\hat{u},v) -
  \mathcal{W}^{0}_{\text{\rm el}}(u,v) .\label{vis}
\end{equation}
 With no loss of generality we can assume by  density 
that  $\hat{u}$  has the form 
\begin{equation*}
  \hat{u}:= u+\widetilde{u} \quad \text{ where } \widetilde{u}\in C^{\infty}_c(\Omega;\R^d).
\end{equation*}
 As   inequality \eqref{semistability epsilon} holds for every
$\widetilde{u}_{\varepsilon}$ such that
$(\widetilde{u}_{\varepsilon},v_{\varepsilon})\in
\widetilde{A}_{\varepsilon}$, i.e., $\widetilde{u}_{\varepsilon}\in
H^{1}_{\Gamma_D}(\Omega;\R^d)$,  we can choose 
\begin{equation*}
  \widetilde{u}_{\varepsilon}:= \hat{u}+u_{\varepsilon}-u=\widetilde{u}+u_{\varepsilon}.
\end{equation*}
Notice that we have
\begin{equation}\label{diff sum u_epsilon}
  \widetilde{u}_{\varepsilon}-u_{\varepsilon}=\widetilde{u} \quad \text{ and } \quad \widetilde{u}_{\varepsilon}+u_{\varepsilon}=\widetilde{u}+2u_{\varepsilon}\rightharpoonup \widetilde{u}+2u \quad \text{ in } H^{1}(\Omega;\R^d).
\end{equation}
To check   inequality \eqref{vis}   we need to show that
  \begin{align}
    &\limsup_{\varepsilon \rightarrow 0}\frac{1}{\varepsilon^2}\left(\int_{\Omega}\left(W_{\text{\rm el}}((I+\varepsilon \nabla \widetilde{u}_{\varepsilon})(I+\varepsilon \nabla v_{\varepsilon})^{-1}) - W_{\text{\rm el}}((I+\varepsilon \nabla {u}_{\varepsilon})(I+\varepsilon \nabla v_{\varepsilon})^{-1})\right){\rm d}X \right)\notag \\
    &\qquad \leq
    \int_{\Omega}\left(|\nabla (\hat{u}-v)|^2_{\C_{\text{\rm el}}}-|\nabla (u-v)|^2_{\C_{\text{\rm el}}}  \right){\rm d}X. \label{Th semistability}
  \end{align}
Let us first study the limiting behaviour of the arguments of these energy densities. We define $(I+\varepsilon \nabla \widetilde{u}_{\varepsilon})(I+\varepsilon \nabla v_{\varepsilon})^{-1}=I+\varepsilon A_{\varepsilon}$, namely
\begin{equation*}
\begin{aligned}
  A_{\varepsilon}&:=\frac{1}{\varepsilon}\left((I+\varepsilon \nabla \widetilde{u}_{\varepsilon})(I+\varepsilon \nabla v_{\varepsilon})^{-1}-I\right)\\
  &=(\nabla \widetilde{u}_{\varepsilon}-\nabla v_{\varepsilon})-\varepsilon \nabla \widetilde{u}_{\varepsilon}\nabla v_{\varepsilon}
  +\varepsilon(I+\varepsilon \nabla \widetilde{u}_{\varepsilon})(\nabla v_{\varepsilon})^2(I+\varepsilon \nabla v_{\varepsilon})^{-1}\\
  &=(\nabla \widetilde{u}_{\varepsilon}-\nabla v_{\varepsilon})-\varepsilon \nabla \widetilde{u}_{\varepsilon}\nabla v_{\varepsilon}
  +  M_{\varepsilon} +\varepsilon \nabla \widetilde{u}_{\varepsilon}M_{\varepsilon}
  ,
\end{aligned}
\end{equation*}
where we have set $M_{\varepsilon}:=\varepsilon(\nabla v_{\varepsilon})^2(I+\varepsilon \nabla v_{\varepsilon})^{-1}$. Similarly, we can write $(I+\varepsilon \nabla u_{\varepsilon})(I+\varepsilon \nabla v_{\varepsilon})^{-1}=I+\varepsilon B_{\varepsilon}$ by letting
\begin{equation*}
\begin{aligned}
  B_{\varepsilon}&:=\frac{1}{\varepsilon}\left((I+\varepsilon \nabla u_{\varepsilon})(I+\varepsilon \nabla v_{\varepsilon})^{-1}-I\right)\\
  &=(\nabla u_{\varepsilon}-\nabla v_{\varepsilon})-\varepsilon \nabla u_{\varepsilon}\nabla v_{\varepsilon}
  +  M_{\varepsilon} +\varepsilon \nabla u_{\varepsilon}M_{\varepsilon}.
\end{aligned}
\end{equation*}
Notice that by definition of $M_{\varepsilon}$ and the fact that $I+\varepsilon \nabla v_{\varepsilon}\in K$ we have
\begin{equation*}
  \|\varepsilon M_{\varepsilon}\|_{L^{\infty}(\Omega)}\leq c, \quad \| M_{\varepsilon}\|_{L^{1}(\Omega)}\leq c \varepsilon\|( \nabla v_{\varepsilon})\|_{L^{2}(\Omega)}\leq c\varepsilon .
\end{equation*}
This implies by interpolation that $M_{\varepsilon}$ is also bounded in $L^2(\Omega;\R^{d\times d})$, hence $M_{\varepsilon}\rightharpoonup 0$ weakly in $L^{2}(\Omega;\R^{d\times d})$. Then, we have
\begin{equation*}
  A_{\varepsilon}-B_{\varepsilon}=(\nabla \widetilde{u}_{\varepsilon}-u_{\varepsilon})(I-\varepsilon \nabla v_{\varepsilon}+\varepsilon M_{\varepsilon}){\stackrel{\eqref{diff sum u_epsilon}}{=}}
  \nabla \widetilde{u}+\nabla \widetilde{u}(-\varepsilon \nabla v_{\varepsilon}+\varepsilon M_{\varepsilon})\rightarrow \nabla \widetilde{u} \, \text{ strongly in } L^2(\Omega;\R^{d\times d}),
\end{equation*}
since $\nabla \widetilde{u}\in C^{\infty}_c(\Omega;\R^{d\times d})$ is bounded in $L^{\infty}(\Omega;\R^{d\times d})$ and $(-\varepsilon \nabla v_{\varepsilon}+\varepsilon M_{\varepsilon})\rightarrow 0$ strongly in $L^2(\Omega;\R^{d\times d})$.
Moreover, by recalling \eqref{diff sum u_epsilon} we have that
\begin{equation*}
  A_{\varepsilon}+B_{\varepsilon}=(\nabla \widetilde{u}_{\varepsilon}+u_{\varepsilon})(I-\varepsilon \nabla v_{\varepsilon}+\varepsilon M_{\varepsilon})-2(\nabla v_{\varepsilon}-M_{\varepsilon})\rightharpoonup \nabla \widetilde{u}+2\nabla u-2\nabla v \, \text{ weakly in } L^2(\Omega;\R^{d\times d}).
\end{equation*}

Fix now $\delta>0$ and let $c_{\rm el}(\delta)$ be as in assumption
\eqref{c_el}. Let us define the set
\begin{equation*}
  \Omega_{\varepsilon}^{\delta}:=\left\{x\in \Omega \;|\; \varepsilon |A_{\varepsilon}|+\varepsilon|B_{\varepsilon}|\leq c_{\rm el}(\delta)  \right\}
\end{equation*}
 containing all points where $\varepsilon |A_{\varepsilon}|$ and
$\varepsilon|B_{\varepsilon}|$ are small. 
Notice that
\begin{equation}\label{bad set meas}
  |\Omega \setminus \Omega_{\varepsilon}^{\delta}|=\int_{\Omega \setminus \Omega_{\varepsilon}^{\delta}}1\, {\rm d}X\leq \frac{\varepsilon^2}{c^2_{\rm el}(\delta)}\int_{\Omega}\left( |A_{\varepsilon}|+|B_{\varepsilon}|\right)^2
  {\rm d}X\leq c\frac{\varepsilon^2}{c^2_{\rm el}(\delta)} \rightarrow 0 \quad \text{ as } \varepsilon\rightarrow 0
\end{equation}
since $A_{\varepsilon}$ and $ B_{\varepsilon} $ are bounded in
$L^2(\Omega;\R^{d\times d})$. We split the integrals in the left-hand
side of \eqref{Th semistability} in the sum of the integrals on the 
sets $\Omega_{\varepsilon}^{\delta}$ and on the complementary sets
$\Omega \setminus \Omega_{\varepsilon}^{\delta}$. By  using assumption \eqref{c_el}, on the 
sets $\Omega_{\varepsilon}^{\delta}$  we have
\begin{align}
  &\frac{1}{\varepsilon^2}\int_{\Omega_{\varepsilon}^\delta}\left(W_{\text{\rm
    el}}(I+\varepsilon A_{\varepsilon}) - W_{\text{\rm
    el}}(I+\varepsilon B_{\varepsilon})\right){\rm d}X\nonumber\\
  &\quad \leq
   \int_{\Omega}\left(|A_{\varepsilon}|^2_{\C_{\rm el}} -|B_{\varepsilon}|^2_{\C_{\rm el}}+\delta(|A_{\varepsilon}|^2_{\C_{\rm el}} +|B_{\varepsilon}|^2_{\C_{\rm el}} )  \right){\rm d}X.\label{second}
\end{align}
 The first term in the right-hand side above can be treated as
follows 
\begin{align*}
    &\int_{\Omega}\left(|A_{\varepsilon}|^2_{\C_{\rm el}}
      -|B_{\varepsilon}|^2_{\C_{\rm el}} )  \right){\rm d}X =\frac{1}{2}\int_{\Omega} \C_{\rm el} (A_{\varepsilon}+B_{\varepsilon}):(A_{\varepsilon}-B_{\varepsilon})){\rm d}X\\
     &\quad\rightarrow \frac{1}{2}\int_{\Omega}\C_{\rm el}((\nabla \hat{u}{-}\nabla v){+}(\nabla u{-}\nabla v)){:}((\nabla \hat{u}{-}\nabla v){-}(\nabla u{-}\nabla v)){\rm d}X\\
    &\quad=\int_{\Omega}\left(|\nabla (\hat{u}-v)|^2_{\C_{\text{\rm el}}}-|\nabla (u-v)|^2_{\C_{\text{\rm el}}}  \right){\rm d}X
  \end{align*}
by  means of  the strong convergence of $  A_{\varepsilon}-B_{\varepsilon}\rightarrow \nabla \widetilde{u}$
and the weak convergence of
 $A_{\varepsilon}+B_{\varepsilon}\rightharpoonup \nabla
 \widetilde{u}+2\nabla u-2\nabla v$ in $L^{2}(\Omega;\R^{d\times
   d})$. On the other hand, the second term  in the right-hand
 side of
 \eqref{second}   satisfies
\begin{equation*}
  \int_{\Omega}\delta(|A_{\varepsilon}|^2_{\C_{\rm el}} +|B_{\varepsilon}|^2_{\C_{\rm el}} )  {\rm d}X\leq \delta c
\end{equation*}
since $A_{\varepsilon} $ and $B_{\varepsilon}$ are bounded in
$L^2(\Omega;\R^{d\times d})$.

Hence, it remains to show that the integrals  in \eqref{Th
  semistability}   on the complements $\Omega\setminus
\Omega_{\varepsilon}^{\delta}$  converge to $0$    as $\varepsilon \rightarrow 0$. In order to do so, let us define
\begin{equation*}
  F_1:=(I+\varepsilon \nabla u_{\varepsilon})(I+\varepsilon \nabla v_{\varepsilon})^{-1} \quad F_2:= \nabla \widetilde{u}(I+\varepsilon \nabla v_{\varepsilon})^{-1}.
\end{equation*}
Since by definition $\nabla \widetilde{u}_{\varepsilon}=\nabla \widetilde{u}+\nabla u_{\varepsilon}$ and $W_{\rm}$ is locally Lipschitz, we can write
\begin{equation*}
  \begin{aligned}
    \frac{1}{\varepsilon^2} \int_{\Omega\setminus \Omega_{\varepsilon}^{\delta}}&\left(W_{\text{\rm el}}((I+\varepsilon \nabla \widetilde{u}_{\varepsilon})(I+\varepsilon \nabla v_{\varepsilon})^{-1}) - W_{\text{\rm el}}((I+\varepsilon \nabla u_{\varepsilon})(I+\varepsilon \nabla v_{\varepsilon})^{-1})\right){\rm d}X\\
    &=\frac{1}{\varepsilon^2}\int_{\Omega\setminus \Omega_{\varepsilon}^{\delta}} |W_{\rm el}(F_1+\varepsilon F_2)- W_{\rm el}(F_1)|\,  {\rm d}X\leq \frac{1}{\varepsilon^2} \int_{\Omega\setminus \Omega_{\varepsilon}^{\delta}} \varepsilon|F_2|\,{\rm d}X\\
    &\stackrel{\eqref{bad set meas}}{\leq}
    \frac{c}{\varepsilon^2}\frac{\varepsilon^2}{c^2_{\rm el}(\delta)} \varepsilon \rightarrow 0,
  \end{aligned}
\end{equation*}
where we used that $F_2$ is uniformly bounded in
$L^{\infty}(\Omega;\R^{d\times d})$. This concludes the proof of 
inequality \eqref{Th
  semistability}. The check of linearized semistability
\eqref{semistability 0} then follows as soon as one passes to the
limit in the loading terms, which is straightforward. 

  In particular, we have proved that $u$ solves the linear minimization problem
  \begin{equation*}
    \mathcal{W}^{0}_{\text{\rm el}}(u(t),v(t))-\langle \ell^0(t), u(t)\rangle =\argmin_{\hat{u}\in H^{1}_{\Gamma_D}(\Omega;\R^d)} \mathcal{W}^{0}_{\text{\rm el}}(\hat{u},v(t))-\langle \ell^0(t), \hat u\rangle
  \end{equation*}
  for given $v$,
  thanks to \eqref{semistability 0}. Hence, the limit $u$ is unique
  and measurable in time, since it is the image of $v$ through a
  linear operator. We also remark that this implies that subsequences
  in \eqref{u epsilon subseq} can be chosen independently of $t$.

	\section*{Acknowledgements}
 AC and US are supported by the Austrian Science Fund (FWF) through
 project F65. US and MK are partially funded by the FWF-GA\v CR
 project I5149 (GA\v{C}R 21-06569K), the OeAD-WTZ
 project CZ09/2023, and the project M\v{S}MT \v{C}R 8J23AT008. US is partially funded by projects I4354 and P32788. 

	\section*{Statement on data availability}
 Data sharing not applicable to this article as no datasets were
 generated or analysed during the current study.

\end{document}